%% file: paper.tex
\documentclass[preprint]{elsarticle}
\usepackage{amscd}
\usepackage{amssymb, latexsym}
\usepackage{amsmath}
\usepackage{amsfonts}
\usepackage{amscd}
\usepackage{longtable}
\usepackage{amsthm}
\usepackage{url}

\newtheorem{thm}{Theorem}[subsection]
\newtheorem{cor}[thm]{Corollary}
\newtheorem{lem}[thm]{Lemma}
\newtheorem{prop}[thm]{Proposition}

\newtheorem{defn}[thm]{Definition}
\newtheorem{rem}[thm]{Remark}

\numberwithin{equation}{section}

\def\Z{{\mathbb Z}}

\def\:{\colon}
\def\a{\alpha}
\def\b{\beta}
\def\g{\gamma}
\def\d{\delta}

\def\la{\lambda}
\def\v{\nu}
\def\w{\omega}

\def\p{\phi}
\def\t{\otimes}
\def\Gg{{\mathfrak{g}}}
\def\Uu{{\mathfrak{u}}}
\def\h{\operatorname{H}}

\begin{document}

\begin{frontmatter}

\title{Second Cohomology Groups for Algebraic Groups and their Frobenius Kernels}

\author{Caroline B. Wright \corref{cor1}}
\address{Department of Mathematics, University of Georgia, Athens, GA 30602, USA}

\cortext[cor1]{Corresponding Address: Department of Mathematics, University of Arizona, P.O. Box 210089,  Tucson, AZ 85721, USA}
\ead{cwright@math.arizona.edu}

\begin{abstract}
Let $G$ be a simple simply connected algebraic group scheme defined over an algebraically closed field of characteristic $p > 0$.  Let $T$ be a maximal split torus in $G$, $B \supset T$ be a Borel subgroup of $G$ and $U$ its unipotent radical.  Let $F: G \rightarrow G$ be the Frobenius morphism.  For $r \geq 1$ define the Frobenius kernel, $G_r$, to be the kernel of $F$ iterated with itself $r$ times.  Define $U_r$ (respectively $B_r$) to be the kernel of the Frobenius map restricted to $U$ (respectively $B$).  Let $X(T)$ be the integral weight lattice and $X(T)_+$ be the dominant integral weights.\\
\indent The computations of particular importance are $\h^2(U_1,k)$, $\h^2(B_r,\la)$ for $\la \in X(T)$, $\h^2(G_r,H^0(\la))$ for $\la \in X(T)_+$, and $\h^2(B,\la)$ for $\la \in X(T)$.  The above cohomology groups for the case when the field has characteristic 2 one computed in this paper.  These computations complete the picture started by Bendel, Nakano, and Pillen for $p \geq 3$ \cite{BNP2}.
\end{abstract}


\end{frontmatter}

\section{Introduction}

\subsection{} Let $G$ be a simple simply connected affine algebraic group scheme defined over ${\mathbb F}_p$ and $k$ be an algebraically closed field of characteristic $p>0$.  Let $T$ be a maximal split torus of $G$ and $B \supset T$ be the Borel subgroup of $G$.  Let $F:G \rightarrow G$ denote the Frobenius morphism and $G_1$, the first Frobenius kernel, denote the scheme-theoretic kernel of $F$.   More generally, higher Frobenius kernels, $G_r$, are defined by taking the kernel of the iteration of $F$ with itself $r$ times.\\
\indent It is a well known fact that the representation theory of $G_1$ is equivalent to the representation theory of the restricted Lie algebra $\Gg = \operatorname{Lie}(G)$. Knowledge about the second cohomology groups is important because of the information it gives us about central extensions of the underlying algebraic structures.  A central question in representation theory of algebraic groups is to understand the structure and the vanishing of the line bundle cohomology, $H^n(\la)={\mathcal H}^n(G/B,\mathcal{L}(\la))$ for $\la \in X(T)$, where $\mathcal{L}(\la)$ is the line bundle over the flag variety, $G/B$.  Fundamental to the understanding of the line bundle cohomology is the computation of the rational cohomology groups, in particular the calculation of $\h^\bullet(B, \la)$.   Listed below are the calculations of cohomology of Frobenius kernels that aid in the computation of the line bundle cohomology.

\begin{itemize}
\item[(1.1)] $\h^n(\Uu,k)$, where $\Uu = \operatorname{Lie}(U)$
\vspace{.15cm} 
\item[(1.2)] $\h^n(U_1,k)$
\vspace{.15cm} 
\item[(1.3)]  $\h^n(B_1,\la), \mbox{ for } \la \in X(T)$ 
\vspace{.15cm}
\item[(1.4)]  $\h^n(B_r,\la), \mbox{ for } \la \in X(T)$ 
\vspace{.15cm}
\item[(1.5)]  $\h^n(B,\la),  \mbox{ for } \la \in X(T)$ 
\vspace{.15cm}
\item[(1.6)]  $\h^n(G_1,H^0(\la)), \mbox{ for } \la \in X(T)_+$
\vspace{.15cm}
\item[(1.7)]  $\h^n(G_r,H^0(\la)), \mbox{ for } \la \in X(T)_+$

\end{itemize}

\indent The goal of this paper focuses on the calculations (1.1) - (1.7) for $n=2$ and $p=2$, which complete the picture for $n=2$.  Bendel, Nakano, and Pillen \cite{BNP2} computed the above groups for $n=2$ and $p \geq 3$.  There are many theorems that only hold for $p \geq 3$, which makes the computations for $p=2$ harder.  Furthermore, when calculating $\h^2(U_1,k)$, an even index of connection for many types of Lie Algebras creates interesting cases for $p=2$.  The even index of connection creates many case-by-case considerations.  However, in the end we found that $\h^2(G_r, H^0(\lambda))$ always has a good filtration, satisfying Donkin's conjecture for $V=H^0(\la)$.  More details about the good filtrations is found in Appendix D. \footnotemark
\footnotetext{The appendices are on my webpage: \small{\url{http://www.math.arizona.edu/~ cwright/Research.html}} and on the arXiv.}
\subsection{History}
\indent In 1983, Friedlander and Parshall \cite{FP1} calculated various cohomology groups of algebraic groups, starting with the special case when $G$ is the general linear group with coefficients in the adjoint representation; then extended the idea to general algebraic groups with coefficients in $V^{(r)}$, where $V$ is a $G$-module.    They also calculated (1.7) for $n=1,2$ and for $k=H^0(0)$ for $p \not= 2,3$.\\
\indent Andersen and Jantzen \cite{AJ}, determined (1.6) for $p \ge h$, where $h$ is the Coxeter number (i.e., $\displaystyle{h=\left< \rho, \a^\vee \right> +1}$, more precisely, $h$ is one more than the height of the highest root). Andersen and Jantzen also determined (1.3) for $\la = w \cdot 0 + p\v$ for $p > h$, where $w \in W$ and $\v \in X(T)$.  These results originally had restrictions on the type of root system involved, which were later removed by Kumar, Lauritzen, and Thomsen \cite{KLT}.\\
\indent In 1984, Andersen \cite{And} calculated $\h^\bullet(B,w \cdot 0)$, where $w \in W$.  Recently, Andersen and Rian \cite{AR} proved some general results on the behavior of $\h^\bullet(B,\la)$ and developed new techniques to enable the calculation of all $B$-cohomology for degree at most 3 when $p>h$.  They also calculated $\h^2(B,\la)$ and $\h^3(B,\la)$ explicitly for $\la \in X(T)$ and $p > h$.  For higher cohomology groups, they proved the following theorem \cite[3.1,6.1]{AR}:
\begin{thm} Suppose $p > h$.  Let $w \in W, \v \in X(T)$.  Then we have for all $i$
\begin{itemize}
\item[(a)] $\h^i(B,w\cdot 0 + p \v) \cong \h^{i-l(w)}(B,\v).$
\item[(b)] $\h^i(B,p\la)=0$ for $i > -2\cdot\mbox{ht}(\la)$.
\end{itemize}
\end{thm} 
Bendel, Nakano, and Pillen \cite{BNP2} calculated $\h^2(B, \la)$ for $\la \in X(T)$ and $p \geq 3$.  This paper computes $\h^2(B, \la)$ for $p=2$.  In both of these papers, $\h^2(B, \la)$ was calculated by previous calculations of (1.1)--(1.4). \\
\indent The focus in the past 15 years of the calculations for the above calculations changed from large primes to small primes.  In 1991, Jantzen \cite{Jan2} calculated (1.1)--(1.3), (1.6) for all primes.  Jantzen used basic facts about the structure of the root systems and isomorphisms relating the different cohomology groups.  Bendel, Nakano, and Pillen used Jantzen's results to get (1.4) and (1.7) for $n=1$ and all $p$ in \cite{BNP1}.  In 2004, Bendel, Nakano, and Pillen \cite{BNP2} worked out (1.1.1)--(1.1.7) for $n=2$ and $p \geq 3$.\\

\subsection{Notation} Throughout the paper, the standard conventions provided in \cite{Jan1} is followed.  Let $G$ be a simply connected semisimple algebraic group over an algebraically closed field, $k$, of prime characteristic, $p > 0$.  Let $\Gg=\mbox{Lie}(G)$ be the Lie algebra of $G$.  For $r \geq 1$, let $G_r$ be the $r$th Frobenius kernel of $G$.  Let $T$ be a maximal split torus in $G$ and $\Phi$ be the root system associated to $(G,T)$.  The positive (respectively negative) roots are $\Phi^+$ (respectively $\Phi^-$), and $\Delta$ is the set of simple roots.  Let $B \supset T$ be the Borel subgroup of $G$ corresponding to the negative roots and let $U$ be the unipotent radical of $B$.  Let $\Uu = \operatorname{Lie}(U)$ be the lie algebra of the unipotent radical.  For a given root system of rank $n$ denote the simple roots $\a_1, \a_2, \ldots, \a_n$, adhering to the ordering used in \cite{Jan2} (following Bourbaki).  In particular, for type $B_n$, $\a_n$ denotes the unique short simple root; for type $C_n$, $\a_n$ denotes the unique long simple root; for type $F_4$ $\a_1$ and $\a_2$ are the short simple roots; for type $G_2$, $\a_1$ is the unique short simple root.  If $\a \in \Phi$, and $\a = \sum_{i=1}^n m_i\a_i$ then the height of $\a$ is defined by ht($\a) \:= \sum_{i=1}^n m_i$.\\
\indent Let $\mathbb{E}$ be the Euclidean space associated with $\Phi$, and the inner product on $\mathbb{E}$ will be denoted by $\left< , \right>$.  For any root $\a$ denote the dual root by $\a^\vee = \frac{2 \a}{\left< \a, \a \right>}$.  Let $\omega_1, \omega_2, \ldots, \omega_n$ be the fundamental weights and $X(T)$ be the integral weight lattice spanned by these fundamental weights.  The set of dominant integral weights is denoted by $X(T)_+$ and the set of $p^r$-restricted weights is $X_r(T)$.  The simple modules for $G$ are indexed by the set $X(T)_+$ and denoted by $L(\la)$, $\la \in X(T)_+$ with $L(\la)$ = soc$_G H^0(\la)$, where $\textnormal{soc}_G H^0(\la)$ is the socle of the $G$-module $H^0(\la)$. where $H^0(\la) =$ ind$_B^G \la$.  Here $\la$ denotes the one-dimensional $B$-module obtained by extending the character $\la \in X(T)_+$ to $U$ trivially.\\
\indent  Given a $G$-module, $M$, then composing a representation of $M$ with $G$ results in a new representation where $G_r$ acts trivially, where $M^{(r)}$ denotes the new module.  For any $\la$ in $X(T)$, the $\la$ weight space of $M$ is the $p^r\la$ weight space of $M^{(r)}$.  One the other hand if $V$ is a $G$-module on which $G_r$ acts trivially, then there is a unique $G$-module $M$, with $V=M^{(r)}$.  We denote $M=V^{(-r)}$.\\

\subsection{Outline of Computations}

In recent work, Bendel, Nakano, and Pillen \cite{BNP2} calculated $\operatorname{H}^2(G_r,H^0(\lambda))$ for $p\geq3$ by reducing the calculations down to $\h^2(\Uu,k)$.  Similar strategies are used to calculate $\h^2(G_r,H^0(\la))$ when $p=2$.  The first step uses the following isomorphism to reduce the calculation to $\h^2(B_r,\la)$:
\begin{eqnarray} \label{BrGrrelation}\operatorname{H}^2(G_r,H^0(\lambda))^{(-r)} \cong \mbox{ind}_B^G(\operatorname{H}^2(B_r,\lambda)^{(-r)}).
\end{eqnarray}
The Lyndon-Hochschild-Serre spectral sequence reduces the problem to $\h^2(B_1,\la)$.  The problem is further reduced to the computation of $\operatorname{H}^2(U_1,k)$ via the isomorphism
\begin{eqnarray} \label{U1B1relation}
\h^2(B_1,\lambda) \cong (\h^2(U_1,k) \t \lambda)^{T_1}.
\end{eqnarray}
This isomorphism tells us that the $B_1$-cohomology can easily be determined by looking at particular weight spaces of $\h^2(U_1,k)$.  That is $\h^2(B_1,\la)\cong \h^2(U_1,k)_{-\la}$.\\
\indent The $B$-cohomology completes the calculations for the second cohomology groups as shown in \cite{BNP2} and  Section \ref{Bcohom}.   The following theorem from Bendel, Nakano, and Pillen \cite{BNP2} states the results for $p \geq 3$: 
\begin{thm}
Let $p \geq 3$ and $\la \in X(T)$.
\begin{itemize}
\item[(a)] Suppose $p > 3$ or $\Phi$ is not of type $G_2$.  Then
\[ \h^2(B,\la) \cong \left\{ \begin{array}{lll}
k && \mbox{if } \la=p^lw \cdot 0, \mbox{ with } 0 \leq l, \mbox{ for } w \in W \mbox{ and } l(w)=2,\\
k && \mbox{if } \la=-p^l\a, \mbox{ with } 0 < l \mbox{ and } \a \in \Delta,\\
k && \mbox{if } \la=-p^k\b-p^l\a, \mbox{ with } 0 \leq l < k \mbox{ and } \a,\b \in \Delta,\\
0 && \mbox{else.}\\
\end{array} \right. \]
\item[(b)] Suppose $p=3$ and $\Phi$ is of type $G_2$.  Then
\[ \h^2(B,\la) \cong \left\{ \begin{array}{lll}
k && \mbox{if } \la=p^lw \cdot 0, \mbox{ with } 0 \leq l, l(w)=2,\\
k && \mbox{if } \la=-p^l\a, \mbox{ with } 0 < l \mbox{ and } \a \in \Delta,\\
k && \mbox{if } \la=-p^k\b-p^l\a, \mbox{ with } 0 \leq l < k \mbox{ and } \a,\b \in \Delta,\\
&& \quad \mbox{where } k \not= l+1 \mbox{ if } \b=\a_1 \mbox{ and } \a=\a_2,\\ 
0 && \mbox{else.}\\
\end{array} \right. \]
\end{itemize}
\end{thm}
$\h^2(B,\la)$ is at most one-dimensional, as shown in  \cite{BNP2} and Section \ref{Bcohom}.

\noindent \section{Restricted Lie algebra cohomology}
\subsection{Observations on $U_1$-cohomology} \label{observations} 
\indent Recall that $\h^\bullet(\Uu,k)$ for $p \geq 3$ may be computed by the cohomology of the following complex, using the exterior algebra:
$$k \stackrel{d_0}{\rightarrow}\mathfrak{u}^* \stackrel{d_1}{\rightarrow} \Lambda^2(\mathfrak{u})^* \stackrel{d_2}{\rightarrow} \Lambda^3(\mathfrak{u})^* \rightarrow \ldots .$$  
$\h^\bullet(U_1,k)$ can be computed from $\h^\bullet(\Uu,k)$ by using the Friedlander-Parshall spectral sequence, which only holds for $p \geq 3$.\\  
\indent To calculate $\operatorname{H}^\bullet(U_1,k)$, for $p=2$, we must take a different approach, using the restricted Lie algebra cohomology.  Tthe restricted Lie algebra cohomology may be computed by the cohomology of the following complex \cite[9.15]{Jan1}.
$$k \stackrel{d_0}{\rightarrow}\mathfrak{u}^* \stackrel{d_1}{\rightarrow} S^2(\mathfrak{u})^* \stackrel{d_2}{\rightarrow} S^3(\mathfrak{u})^* \rightarrow \ldots .$$
The differential $d_1$ is a derivation on $S^\bullet(\mathfrak{u}^*)$ and is thus determined by its restriction to $\mathfrak{u}^*$.  Consider the following composition of maps
$$\mathfrak{u}^* \stackrel{d_1}{\hookrightarrow} S^2(\mathfrak{u})^* \stackrel{\pi}{\rightarrow}{{\Lambda} ^2(\mathfrak{u})^*}$$
where $\pi$ is a surjection with kernel $\{f^2:f \in \mathfrak{u}^*\}$ and $\partial = \pi \circ f$ being the coboundary operator for the ordinary Lie algebra cohomology, i.e. the dual of $\Lambda^2(\mathfrak{u}) \rightarrow \mathfrak{u}$ with $a \wedge b \mapsto [a,b]$.  The differentials are given as follows: $d_0 = 0$ and $d_1:\mathfrak{u}^* \rightarrow S^2(\mathfrak{u})^*$ with
$$(d_1\phi)(x_1 \t x_2) = -\phi([x_1,x_2])$$
where $\phi \in \mathfrak{u}^*$ and $x_1,x_2 \in \mathfrak{u}$.  For higher differentials, we identify $S^n(\mathfrak{u})^* \cong S^n(\mathfrak{u}^*)$ and the differentials are determined by the following product rule:
$$d_{i+j}(\phi \t \psi)=d_i(\phi) \t \psi + (-1)^i\phi \t d_j(\psi).$$

\subsection{Basic Results}
Recall the following theorem from Jantzen, \cite{Jan2}

\begin{thm} \label{thm:Jan} $\h^1(U_1,k) \cong \h^1(\Uu,k)$
\end{thm}

\indent The following results, similar to those found in \cite[2.4]{BNP2}, help identify some limitations on which linear combinations of tensor products, $\phi_\a \t \phi_\b$, can represent cohomology classes when char $k=2$.  Using the additive property of differentials and the fact that differentials preserve the $T$ action, then we are interested in linear combinations of tensor products that have the same weight.  Recall the following definition from \cite{BNP2}.

\begin{defn}
An expression $\sum c_{\a,\b} \p_\a \t \p_\b \in S^2(\Uu^*)$ is in \it{reduced} form if $c_{\a,\b} \not= 0$ and for each pair $(\a,\b)$ $c_{\a,\b}$ appears at most once.
\end{defn}

\begin{prop} \label{prop:d2} Let $p=2$ and $x = \sum c_{\a,\b}\;\p_\a \t \p_\b$ ($\a, \b \in \Phi^+$) be an element in $S^2({\mathfrak u^* })$ in reduced form of weight $\g$ (for some $\g \in X(T)$ and $\g \not\in pX(T)$).  If $d_2(x) = 0$, then $d_1(\p_\a)=0$ for at least one $\a$ appearing in the sum.
\end{prop}
\begin{proof}
For any $\a \in \Phi^+$, if $d_1(\p_\a) = \sum c_{\d,\eta} \p_\d \t \p_\eta$, then ht($\d$) $<$ ht($\a)$ and ht($\eta$) $<$ ht ($\a)$ for all $\d, \eta$.  For all $\a$ and $\b$ appearing in the sum for $x$, choose a root $\sigma$ with ht($\sigma)$ being maximal.  Without loss of generality, we may assume $\p_\sigma$ appears in the second factor of the tensor product.  Consider the corresponding term $c_{\a,\sigma} \p_\a \t \p_\sigma$.  Computing $d_2(x)$, one of the components will be $c_{\a,\sigma} d_1(\p_\a) \t \p_\sigma$.  By height considerations, $\p_\sigma$ appears in no other terms, thus it is not a linear combination of the other terms.  Therefore, $d_1(\p_\a)=0$.
\end{proof}

\begin{cor}\label{cor:d2}Let $p=2$.\\
\begin{itemize}
\item[(a)] Let $x \in \h^2(U_1,k)$ be a representative cohomology class in reduced form having weight $\g$ for some $\g \in X(T)$, $\g \not\in pX(T)$.  Then one of the components of $x$ is of the form $\p_\a \t \p_\b$ for some simple root $\a \in \Delta$ and positive root $\b \in \Phi^+$ (with $\a + \b = \g$).
\item[(b)] Suppose $\p_\a \t \p_\b$ represents a cohomology class in $\h^2(U_1,k)$.  Then one of three things must happen either
\begin{itemize}
\item[(i)]$\a,\b \in \Delta$, 
\item[(ii)]$\a \in \Delta$, then $d_1(\phi_\b)= \sum_{\sigma_1+\sigma_2=\b}c_{\sigma_1,\sigma_2} \; \phi_{\sigma_1} \t \phi_{\sigma_2}$, then $c_{\sigma_1,\sigma_2}=\pm 2$ for all decompositions of $\b$ (that is the structure constant is even), or 
\item[(iii)]$\a=\b$ and $\a \in \Phi^+$.
\end{itemize}
\end{itemize}
\end{cor}

\begin{proof} Part (a) follows immediately from the previous proposition and Theorem \ref{thm:Jan}, since $\h^1(\Uu,k)$ is generated by the simple roots.  For part (b) let's first assume that $\a=\b$, then $$d_2(\p_\a \t \p_\a) = d_1(\p_\a) \t \p_\a + \p_\a \t d_1(\p_\a) = 2d_1(\p_\a) \t \p_\a = 0.$$
Now, assume $\a \not= \b$ and $\a$ is simple. Then, 
$$d_2(\p_\a \t \p_\b) = d_1(\p_\a) \t \p_\b + \p_\a \t d_1(\p_\b) = \p_\a \t d_1(\p_\b)=0.$$
Hence, $d_1(\p_\b) = 0$.  Therefore, $\b \in \Delta$ or if $\b = \sigma_1 + \sigma_2$ then $c_{\sigma_1,\sigma_2}=\pm 2$, for all decompositions of $\b$.
\end{proof}

\subsection{Root Sums}
As previously mentioned, the computation of $\h^2(U_1,k)$ involves information about $B_1$- and $B$-cohomology.  In this process, certain sums involving positive roots arise.  Suppose $x \in \h^2(U_1,k)$ has weight $\g \in X(T)$.  Then by Corollary \ref{cor:d2}, $\g = \a + \b$ for $\a \in \Delta$ and $\b \in \Phi^+$ and $\a \not= \b$.  Given such roots $\a$ and $\b$, we want to know whether there exists a weight $\sigma \in X(T)$, $\b_1, \b_2 \in \Delta$, and integers $0 \leq i \leq p-1$ and $m \geq 0$ such that any of the following hold:
\begin{equation}\label{An:1}
\a + \b = 2 \sigma,
\end{equation}
\begin{equation}\label{An:2}
\a + \b = \b_1 + 2\sigma,
\end{equation}
\begin{equation}\label{An:3}
\a + \b = i\b_1 + 2^m\b_2 + 2 \sigma.
\end{equation}
Given $\g$ a weight of $\h^2(U_1,k)$, then there is a weight $\v \in X(T)$ such that $\h^2(B,-\g+p\v) \not = 0$.  Using results on $B$-cohomology due to Andersen, \cite[2.9]{And}, then $\g$ must satisfy \eqref{An:1}.  Note that \eqref{An:1} and \eqref{An:2} are special cases of (\ref{An:3}) (i.e. when $i=0$ and $m=0$).  Equation (\ref{An:1}) arises from the reduction $\h^2(B_1,k)=\h^2(U_1,k)^{T_1}$.  For more details on how these equations arise see \cite{BNP2}. 

\begin{rem}These sums noted above are only valid when $2$ does not divide the index of connection.\end{rem}

\subsection{$U_1$-cohomology}
We state the theorem for $\operatorname{H}^2(U_1,k)$ when $p=2$, which is a summand of $T$-weights, except for $\Uu^*$.  In the next chapter, we explain the proof for each type.  
\begin{thm} \label{thm:U1} As a $T$-module,\\ 
\begin{itemize}
\item[(a)] If $\Phi=A_n$, then
$$\operatorname{H}^2(U_1,k) \cong (\mathfrak{u}^*)^{(1)}\; \oplus \bigoplus_{\stackrel{\a,\b \in \Delta}{\a+\b\not\in \Phi^+}} -(s_\a s_\b) \cdot 0 \; \oplus \bigoplus_{\stackrel{\stackrel{\a,\b, \g \in \Delta}{\a+\b\not\in \Phi^+}}{{\a+\b+\g \in \Phi^+}}} -(s_\a s_\b) \cdot 0 + 2\g.$$
\item[(b)] If $\Phi=B_n$, then
\begin{align*}
\operatorname{H}^2(U_1,k) &\cong (\mathfrak{u}^*)^{(1)} \; \oplus \bigoplus_{\stackrel{\a,\b \in \Delta}{\a+\b\not\in \Phi^+}} -(s_\a s_\b) \cdot 0 \; \oplus \bigoplus_{\stackrel{\stackrel{\a,\b, \g \in \Delta}{\a+\b\not\in \Phi^+}}{{\a+\b+\g \in \Phi^+}}} -(s_\a s_\b) \cdot 0 + 2\g \\ 
& \quad \oplus \bigoplus_{1\leq i \leq n-3} -(s_{\a_i}s_{\a_{n-1}}) \cdot 0 + 2\a_n \; \oplus \bigoplus_{1\leq i \leq n-1} 2(\a_i + \a_{i+1} + \ldots + \a_n).\\
\end{align*}
\item[(c)] If $\Phi=C_n$, then
\begin{align*}
\operatorname{H}^2(U_1,k) &\cong (\mathfrak{u}^*)^{(1)}\; \oplus \bigoplus_{\stackrel{\a,\b \in \Delta}{\a+\b\not\in \Phi^+}} -(s_\a s_\b) \cdot 0 \; \oplus \bigoplus_{\stackrel{\stackrel{\a,\b, \g \in \Delta}{\a+\b\not\in \Phi^+}}{{\a+\b+\g \in \Phi^+}}} -(s_\a s_\b)\; \cdot 0 + 2\g \\
& \quad \oplus \bigoplus_{1\leq i \leq n-3} -(s_{\a_i}s_{\a_{n}}) \cdot 0 + 2\a_{n-1} \; \oplus \bigoplus_{1\leq i \leq n-1} 2(\a_i + \a_{i+1} + \ldots + \a_n) \\
& \quad \oplus -(s_{\a_n-1}s_{\a_n}) \cdot 0.
\end{align*}
\item[(d)] If $\Phi=D_n$, $n \geq 4$, then
\begin{align*}
\operatorname{H}^2(U_1,k) &\cong (\mathfrak{u}^*)^{(1)} \; \oplus \bigoplus_{\stackrel{\a,\b \in \Delta}{\a+\b\not\in \Phi^+}} -(s_\a s_\b) \cdot 0 \; \oplus \bigoplus_{\stackrel{\stackrel{\a,\b, \g \in \Delta}{\a+\b\not\in \Phi^+}}{{\a+\b+\g \in \Phi^+}}} -(s_\a s_\b) \cdot 0 + 2\g \\
& \quad \oplus -(s_{\a_{n-3}}s_{\a_{n-1}}) \cdot 0 + 2(\a_{n-2}+\a_n )\\
& \quad  \oplus -(s_{\a_{n-3}}s_{\a_n}) \cdot 0 + 2(\a_{n-2}+\a_n )\\
& \quad \oplus \bigoplus_{1 \leq i \leq n-3} -(s_{\a_{n-1}}s_{\a_n}) \cdot 0 + 2(\a_i+\ldots+\a_{n-2}).  \\
\end{align*}
\item[(e)] If $\Phi=E_n$, $n=6,7,8$ then
\begin{align*}
\operatorname{H}^2(\mathfrak{u},k) &\cong \bigoplus_{\stackrel{\alpha,\beta \in \Delta}{\alpha+\beta\not\in \Phi^+}} -(s_\alpha s_\beta) \cdot 0 \; \oplus \bigoplus_{\stackrel{\stackrel{\alpha,\beta, \gamma \in \Delta}{\alpha+\beta \not\in \Phi^+}}{\alpha+\beta+\gamma \in \Phi^+}} -(s_\alpha s_\beta) \cdot 0 + 2\gamma  \\
& \oplus \bigoplus_{i=5}^n -(s_{\alpha_2}s_{\alpha_3}) \cdot 0 + 2(\alpha_4+\ldots + \alpha_i) \oplus -(s_{\alpha_2}s_{\alpha_5}) \cdot 0 + 2\alpha_4 \\
& \oplus -(s_{\alpha_3}s_{\alpha_5}) \cdot 0 + 2(\alpha_2+\alpha_4)  \oplus -(s_{\alpha_2}s_{\alpha_5}) \cdot 0 + 2(\alpha_1+\alpha_3+\alpha_4). 
\end{align*}
\item[(f)] If $\Phi = F_4$, then
\begin{align*}
\operatorname{H}^2(U_1,k) &\cong (\mathfrak{u}^*)^{(1)} \; \oplus \bigoplus_{\stackrel{\a,\b \in \Delta}{\a+\b\not\in \Phi^+}} -(s_\a s_\b) \cdot 0 \; \oplus -(s_{\a_1}s_{\a_3}) \cdot 0 + 2\a_2  \\
& \quad \oplus -(s_{\a_1}s_{\a_3}) \cdot 0 + 2(\a_2+\a_3) \; \oplus 2(\a_2+\a_3) \oplus 2(\a_1+\a_2+\a_3)  \\
& \quad \oplus 2(\a_2+\a_3+\a_4) \oplus 2(\a_1+\a_2+\a_3+\a_4) \; \oplus 2(\a_2+2\a_3+\a_4) \\
& \quad \oplus 2(\a_1+\a_2+\a_3+\a_4).
\end{align*}
\item[(g)] If $\Phi = G_2$, then
$$\operatorname{H}^2(U_1,k) \cong (\mathfrak{u}^*)^{(1)} \; \oplus 2(\a_1+\a_2).$$
\end{itemize}
\end{thm}
Refer to Appendix A to see a complete list of the cohomology classes associated with each $w \in W$.

\section{Proof of Theorem \ref{thm:U1}}
\subsection{Type $A_{2n}$} \label{proof:An,1}
  \indent Suppose $x \in \operatorname{H}^2(U_1,k)$ has weight $\g \in X(T)$.  From corollary \ref{cor:d2}, we know that $\g = \a  + \b$ for some roots $\a \in \Delta$ and $\b \in \Phi^+$, with $\a \not= \b$.   Furthermore, since 2 does not divide the index of connection, then we know $\g = \a + \b$ must satisfy one of the equations \eqref{An:1} -- \eqref{An:3}.

\begin{prop} \label{prop:An,1} Let $\Phi = A_n$ and $p=2,\; \a \in \Delta,\; \b \in \Phi^+,$ and $\a \not= \b$.  Then there is no weight $\sigma \in X(T)$ such that $\a + \b = 2 \sigma$.
\end{prop}
\begin{proof} Consider $\a + \b = 2\sum m_i \a_i$, where $m_i = 0,1, \forall i$. Then, $\a = \b$, which contradicts the hypothesis.
\end{proof}

\begin{prop} \label{prop:An,2} Let $p=2, \; \a \in \Delta,\; \b \in \Phi^+$, and $\a \not= \b$.  Then, there is no simple root $\b_1 \in \Delta$ and $\sigma \in X(T)$ such that $\a + \b = \b_1 + 2\sigma$.
\end{prop}

\begin{proof} Consider $\a + \b = \b_1 + 2\sum m_i\alpha_i$.  Since $\a =\a_i$ for some $i$, then $\sigma \in \Delta$.  Thus the only possibility for $\b$ is $\a_{i-1} + \a_i$ or $\a_i + \a_{i+1}$.  Then $x \in \h^2(U_1,k)$ has only one component and by Corollary \ref{cor:d2}(b), $\b \in \Delta$.  Hence, there does not exist $\b_1 \in \Delta$ and $\sigma \in X(T)$ such that $\a + \b = \b_1+2\sigma.$
\end{proof}

\begin{prop}\label{prop:An,3} Let $p=2$, $\a \in \Delta, \b \in \Phi^+$, and $\a \not= \b$.  If $\a + \b$ is a weight of $\operatorname{H}^2(U_1,k)$ and there exists $\b_1, \b_2 \in \Delta$, $\sigma \in X(T), 0 < i < p$,  and $m \geq 0$ such that
$$\a + \b = i \b_1 +2^m \b_2 +2 \sigma,$$
then one of the following holds
\begin{itemize}
\item[(a)] $\a + \b$ is a solution to equation (\ref{An:1}) or (\ref{An:2}),
\item[(b)] If $n \geq 3$, then $\a + \b = \a_{i-1}+\a_{i+1} + 2\a_i$ or $\a+\b = \a_{i-1}+\a_{i-2}+2\a_i$ or $\a+\b = \a_{i+1}+\a_{i+2}+2\a_i$ for $i \leq n-2$.
\end{itemize}
\end{prop}
\begin{proof}
First note that we only have to consider the cases $i=0,1$  and since $p=2$ if $m \geq 2$, then by choosing a different $\sigma \in X(T)$, these equations reduce down to $m \leq 1$ and $i=0,1$.  If $i=0=m$, then the equation reduces to (\ref{An:2}), which is done.  If $i=0, m=1$, then we have that $\a+\b = 2(\b_2+\sigma)$ and so $\sigma=0$, but then we are back into equation (\ref{An:1}). If $i=1,m=1$, then $\sigma=0$, which is a specific case of equation (\ref{An:2}). The only case we have to check is $i=1,m=0$.  Since $\a=\a_i$, then $\sigma \in \Delta$ and $\b=\a_{i-2}+\a_{i-1}+\a_i$, $\b=\a_i+\a_{i+1}+\a_{i+2}$, or $\b=\a_{i-1}+\a_i+\a_{i+1}$, which are the cases above.  So $\b_1$ and $\b_2$ are either the 2 simple roots on either side of $\sigma$ or the two simple roots to the right or left of $\a$.
\end{proof}

\subsection{Type $A_{2n+1}$}
Note that $X(T)/\Z\Phi = \{t\omega_1 + \Z\Phi:t=0,1, \ldots, n\} \cong \Z_{n+1}$ and $n$ is odd.  Since $2$ divides the index of connection, we must change equations \eqref{An:1}--\eqref{An:3}
$$t\omega_1 = \frac{t}{n+1}(n\a_1+(n-1)\a_2 + \ldots + \a_n).$$
By revising (3.1)-(3.3), for $\a \in \Delta, \b \in \Phi^+$ must satisfy one of the following:
\begin{equation} \label{An:4}
\a + \b = 2t\omega_1 + 2\sigma,
\end{equation}
\begin{equation} \label{An:5}
\a + \b = \b_1 + 2t\omega_1 + 2\sigma,
\end{equation}
\begin{equation} \label{An:6}
\a + \b = i\b_1 + 2^m\b_2 + 2t\omega_1+2\sigma,
\end{equation}
where $\sigma \in \Z\Phi$.  Since $2t\omega_1$ must lie in $\Z\Phi$, $\frac{2t}{n+1} \in \Z$ and $2|n+1,$ it follows that $\frac{t}{s}~ \in~ \Z$, where $s:=\frac{n+1}{2}$.  If $2|\frac{t}{s}$ then (\ref{An:4})-(\ref{An:6}) reduces to the original equations, (\ref{An:1})-(\ref{An:3}) with $\sigma$ lying in the root lattice, and the arguments in Section \ref{proof:An,1} apply.  So, we can assume  that $\frac{t}{s} \not\equiv 0 \pmod 2$.  Consider $\a + \b = \sum_{i=1}^n m_i\a_i$, then $m_i \in \{0,1,2\}$ for $i=1,2, \ldots, n$ and $m_i=2$ for at most one $i$.  To examine the possibilities, reduce $t\w_1$ mod 2, so that we are looking at sequences of 0's and 1's.  The sequence looks like $(1,0,1,0,\ldots, 0,1)$, where one of the zeroes is a 2, and at most two other zeroes can be made into a one by adding a simple root (as in \eqref{An:5} or \eqref{An:6}).  So, we have the following:
$$\sum m_i \a_i = i\b_1+2^m\b_2+(1,0,1,0,\ldots, 0,1) \pmod 2.$$
Since the roots of $A_n$ have consecutive 1's, then $n \geq 9$ has a trivial solution.
\indent Looking at $A_3, A_5,$ and $A_7$ separately it is easy to check that no additional cohomology classes occur, and the only classes that occur are weights of the form $\a + \b = s_{\a}s_{\b} \cdot 0$ and $\a + 2\g + \b = s_{\a}s_{\b}+2\g$, where $\a+\b$ is not a root and $\a+\g+\b$ is a root.

\subsection{Type $B_n$}
For type $B_n$, $X(T)/\Z\Phi \cong \Z_2$ where $\omega_n = \frac{1}{2}(\a_1+2\a_2 + \ldots + n\a_n)$ is a generator, which forces us to revise (\ref{An:4})-(\ref{An:6}) in the following way.  We are looking for $\a,\b \in \Phi^+$ satisfying 
\begin{equation} \label{Bn:1}
\a + \b = 2t\omega_n + 2\sigma,
\end{equation}
\begin{equation} \label{Bn:2}
\a + \b = \b_1 + 2t\omega_n + 2\sigma,
\end{equation}
\begin{equation} \label{Bn:3}
\a + \b = i\b_1 + 2^m\b_2 + 2t\omega_n+2\sigma,
\end{equation}
where $\b_1, \b_2, \sigma \in \Z\Phi$ and $t=0,1$.  Furthermore, when $\Phi$ is of type $B_n$, we have even structure constants.  The structure constants are even when the resulting root has a coefficient of 2 in the $\a_n$ spot and $\a_n$ is broken up between the two roots, i.e. $[\a_n,~ \a_i~+~\a_{i+1}~+~\ldots~+~\a_n]~ =~ 2(\a_i~+~\a_{i+1}~+~\ldots~+~2\a_n)$.\\
\\
{\bf Case 1:}  Suppose $t=0$ and $m_i\not=2$, for any $i$: Arguments in Section \ref{proof:An,1} apply.  Thus, we have the same classes that appear when $\Phi$ is of type $A_n$.\\
{\bf Case 2:}  Suppose $m_i=2$ for some $i$.   For simplicity, the weight will be written as: $(i_1,i_2, \ldots, i_n)$, where $i_j$ is the coefficient of the $\a_j$ term.\\
{\bf Case 2.1:} If $i_j \equiv 0 \pmod 2$ for all $j$, then our weight is twice a root.  When, this occurs, it is easy to check that this results in the case when $\a = \b$ and $\a=\a_i, \b = \a_i + 2\a_{i+1} + 2\a_{i+2} + \dots + 2\a_n$.  It's easy to check that these are always in the kernel. \\
{\bf Case 2.2:} If there is only one odd number, then there must be some 2's.  If all of the two's occur after the odd number, then this is a root and thus in the previous image.  If all of the 2's occur before the odd number, then there can only have one 2 before it because of the pattern of the roots.  In which case, there is only one term in the class, $\phi_{\a_i}\t\p_{\a_i+\a_{i+1}}$.\\
{\bf Case 2.3:} Suppose there are 2 odd numbers.  Note, there can't be more than one 2 between the two odd numbers.  Otherwise, there would be no way to break up the weight into a sum of a simple and a positive root.
Consider the following weight, $(0, 0, \ldots, 0,1, 0, \ldots, 0,1,2)$ where the first 1 is in the $i$th spot, given by $\phi_{\a_i} \t \p_{\a_{n-1}+2\a_n}$.  This is in the kernel because the structure constant is 2 (as noted above).   Besides this weight, the 1's, 2's, and 3's that appear in the weight must all be consecutive.  In particular, the 2's and 3's must be consecutive, which follows from Proposition \ref{prop:d2} and the the structure constants.  Furthermore, there can only be one 3, which follows from Corollary \ref{cor:d2} and the structure of the roots.  If a 3 occurs, then it cannot be in the $\a_n$ spot (unless $t=0$), 2's must follow and a 1 must occur before the 3 with a 2 between the 1 and the 3.  So that leaves us with a weight looking like $(0,\ldots ,0,1,2,3,2,\ldots , 2)$, with the $3$ in the $i$th spot.    Suppose $i \leq n-2$, then the term in the cohomology class is $\p_{\a_i} \t \p_{x_\a}$.  However, after taking the differential the following term appears $\p_{\a_i}\t\p_{x_a}\t\p_{x_b}$ for some $x_a, x_b \in \Phi^+$, which can't cancel out and thus not in the kernel.  The only possible weight is $(0, \ldots , 0,1,2,3,2)$, which is in the kernel.\\  If $t=0$, then the only weight is $(0, \ldots , 0, 1,2,3)$, which is also in the kernel.\\
{\bf Case 2.4:} Suppose there are more than 2 odd numbers.  Then, the weight would look like: $(0,\ldots,0,1,1,2,...,2,3,2,\ldots,2)$.  Suppose the 3 is in the $i^{\text{th}}$ spot and the first 1 is in the $j^{\text{th}}$ spot, then breaking this weight into 2 roots, you must have $\a$ contain $\a_i$ as part of it's sum.  When taking the differential, the term $\a_i \t \a_j \t \p_{x_a}$ appears, which can't cancel with anything else.  Therefore, there can only be at most 2 odd numbers.\\
For a complete list of the cohomology classes that appear refer to Appendix B, found on the ArXiv and http://www.math.uga.edu/~cwright/appendices.pdf.
  
\subsection{Type $C_n$}
For type $C_n$,  $X(T)/\Z\Phi \cong \Z_2$, where $\omega_1 = \a_1+\a_2 + \ldots + \a_{n-1}+\frac{1}{2}\a_n$ is a generator.  This forces us to revise (\ref{An:4})-(\ref{An:6}).  We want $\a,\b \in \Phi^+$ satisfying 
\begin{equation} \label{Cn:1}
\a + \b = 2t\omega_1 + 2\sigma,
\end{equation}
\begin{equation} \label{Cn:2}
\a + \b = \b_1 + 2t\omega_1 + 2\sigma,
\end{equation}
\begin{equation} \label{Cn:3}
\a + \b = i\b_1 + 2^m\b_2 + 2t\omega_1+2\sigma,
\end{equation}
where $\sigma \in \Z\Phi$ and $t=0,1$.\\
{\bf Case 1:} Suppose $t=0$ and $\a = \sum m_i \a_i$, $m_i \not= 2$ for all $i$.  Then arguments in Section \ref{proof:An,1} apply.  Thus, we have the same classes that appear when $\Phi = A_n$ also appear when $\Phi = C_n$.  Furthermore, if $\a + \b = 2\sigma$, then the weight is twice a root.  This results in the cases $\phi_\a \t \p_\a$ and $\a = \a_n, \b = 2\a_i + 2\a_{i+1} + \ldots + 2\a_{n-1} + \a_n$.  It's easy to check that both of these classes are in the kernel.\\
{\bf Case 2:} Suppose there is only one odd number in the weight, then there must be a 1 and 2's must appear.  If the 2's appear only before or after the 1, then there can only be 1 two because of the structure of the roots.  Furthermore, there is only one term in the class with $\b$ not simple and $[\a,\b] \not \equiv 0 \pmod 2$.  However, if the 2 are before and after the 1, then we must have $(0, \ldots, 0, 2,1,2)$.  However, this weight doesn't satisfy Corollary \ref{cor:d2}.  Therefore, there are no weights with only one odd number.\\
{\bf Case 3:} Note that outside of the weight $(1, 0, \ldots, 0, 2, 1)$, (which is in the kernel then all 1's, 2's, and 3's), must be consecutive.  Suppose there are at least 2 odd numbers.  First note that $(0, \ldots, 0, 1, 2, \ldots, 2, 1, 0, \ldots, 0)$ can only be written as a sum of a simple root and a positive root when the last 1 is in the $\a_n$ spot, which then makes this weight a root, thus in the previous image.  Suppose 3's appear in the weight, then there can only be one 3.  Suppose 3 is in the $i$th spot, where $i < n-1$.  Also, there can't be a 1 appearing in the weight before the 3 because of the root structure.  If there are 2's appearing in the weight, $(0, \ldots, 0, 2, \ldots, 2, 3, 2, \ldots, 2, 1)$, then after taking the differential, the term $\p_{\a_i} \t \p_{\a_j} \t \p_{x_a}$ can't be cancelled.  Then, we see that the 3 must be the $n-1$ spot, since it can't be in the spot.  Therefore, we have the following possibilities:\\
(i) $(0, \ldots, 0, 3, 1)$ which is in the kernel, since the structure constant is 2.\\
(ii) $(0, \ldots, 0, 1, 3, 1)$, which doesn't satisfy Corollary \ref{cor:d2}.\\
(iii) $(0, \ldots, 0, 1, 2, 3, 1)$ which isn't in the kernel.\\
(iv) $(0, \ldots, 0, 1, 2, \ldots, 2, 3, 1)$, which isn't in the kernel.  (only two ways to write the weight, but only one way to get $\p_{\a_{i-1}} \t \p_{\a_{n-1}} \t \p_\b$.\\
For a complete list of the cohomology classes that appear refer to Appendix B, found on the ArXiv and http://www.math.uga.edu/~cwright/appendices.pdf.\\

\subsection{Type $D_n$}
Suppose $t=0$ and $\a = \sum m_i \a_i$, $m_i \not= 2$ for all $i$.  Then arguments in Section \ref{proof:An,1} apply.  Thus, we have the same classes that appear when $\Phi = A_n$ also appear when $\Phi = C_n$.  Furthermore, if $\a + \b = 2\sigma$, then the weight is twice a root.  This results in the cases $\phi_\a \t \p_\a$ and $\a = \a_n, \b = 2\a_i + 2\a_{i+1} + \ldots + 2\a_{n-1} + \a_n$.  It's easy to check that both of these classes are in the kernel.\\

\subsubsection{$n$ odd}

For type $D_n$ and $n$ odd, $X(T)/\Z\Phi \cong \Z_4$, where $\omega_n = \frac{1}{2}(\a_1+2\a_2 + \ldots + (n-2)\a_2+\frac{n-2}{2}\a_{n-1}+\frac{n}{2}\a_n)$.  We want $\a, \b \in \Phi^+$ satisfying
\begin{equation} \label{Dn:1}
\a + \b = 2t\omega_n + 2\sigma,
\end{equation}
\begin{equation} \label{Dn:2}
\a + \b = \b_1 + 2t\omega_n + 2\sigma,
\end{equation}
\begin{equation} \label{Dn:3}
\a + \b = i\b_1 + 2^m\b_2 + 2t\omega_n+2\sigma,
\end{equation}
where $\sigma \in \Z\Phi$.  Since we need $\frac{t(n-2)}{2} \in \Z$, then $t \equiv 0 \pmod 2$, then we'll consider the cases $t=0,2$.\\
Note that there are no weights with only one odd number.\\
{\bf Case 1:} Suppose $t=0$, and there are no 3's in the weight.  We claim that if there are no 3's involved in the weight, then there can be at most two 2's involved.  There can't be any 2's prior to the first 1 involved in the weight because there is no way to write the weight as a sum of a positive and a simple root.  When 2's appear after the first 1 then the other 1 must occur in either the $\a_{n-1}$ or the $\a_{n}$ spot.  So, there are two possibilities for weights (i) $(0, \ldots, 0, 1, 2, \ldots, 2, 1)$, where the first 2 is in the $i^{\text{th}}$ spot.  Then after taking the differential, we see that the term $\p_{\a_{n-1}} \t \p_{\a_i} \t \p_{\a_\b}$ can only be cancelled when $i=n-2$.\\
(ii) $(0, \ldots, 0, 1, 2, \ldots, 2, 1,2)$ where the first 2 is in the $i^{\text{th}}$ spot.  Then the same argument holds as above.\\
The only extra classes that occur are $(0, \ldots, 0, 1, 2,2,1)$ and $(0, \ldots, 0, 1,2,1,2)$.\\
{\bf Case 2:} Suppose $t=2$ and there are no 3's in the weight.  The 2's must only come after the first 1 in the weight, and we must have 1's in both the $\a_{n-1}$ and $\a_n$ spots.  Then, we have the following cases:\\
 (i)$(0,\ldots ,0,2, \ldots ,2,1,1)$, which is easy to check is in the kernel.\\
 (ii)$(0,\ldots ,0, 1,2,\ldots ,2,1,1)$, which is in the previous image. \\
 (iii)$(0, \ldots ,0,1,1,2,1,1, \ldots ,2,1,1)$, which isn't in the kernel because the term $\p_{\a_i} \t \p_{\a_j} \t \p_{x_a}$ can't cancel with anything, when $i$ is the place of the first 1 and $j$ is the place of the first 2, which are not connected.
 
 \subsubsection{$n$ even} 

\indent Suppose $n$ is even, then $X(T)/\Z\Phi \cong \Z_2 \times \Z_2$, where $\omega_1=\a_1+\ldots+\a_{n-2}+\frac{1}{2}\a_{n-1}+\frac{1}{2}\a_n$ and $\omega_n = \frac{1}{2}(\a_1+2\a_2+\ldots+(n-2)\a_{n-2})+\frac{n-2}{4}\a_{n-1}+\frac{n}{4}\a_n$ are the two generators.  Then for $\a, \b \in \Phi^+$ must satisfy:  
\begin{equation} \label{Dn:4}
\a + \b = 2t\omega_1 + 2s\omega_n + 2\sigma,
\end{equation}
\begin{equation} \label{Dn:5}
\a + \b = \b_1 + 2t\omega_1 + 2s\omega_n + 2\sigma,
\end{equation}
\begin{equation} \label{Dn:6}
\a + \b = i\b_1 + 2^m\b_2 + 2t\omega_1+2s\omega_n+2\sigma,
\end{equation}
where $i=0,1, m \geq 0, s,t = 0,1, \sigma \in \Z\Phi$, and $\a, \b \in \Phi^+$.\\
{\bf Case 1:} When $s=t=0$, then we can use the same proof as the case when $t=0$ when $n$ is odd.\\
{\bf Case 2:} When $s=1,t=0$, the equations reduce down to the case when $n$ is odd and $t=2$.\\
{\bf Case 3:}  If $t=1,s=0$, then $n \leq 12$ because otherwise there are too many odd numbers in the weight.  However the root structure tells us that if there are any 3's in the weight, then there must be 1's in the $\a_n$ and the $\a_{n-1}$ spots, then we have that $n \leq 6$, but since $n$ is even, then $n=4,6$.\\
{\bf Case 3.1}  When $n=4$, then this is reduced down to the case when $s=t=0$.\\
{\bf Case 3.2}  When $n=6$, then we are reduced down to the case when $s=1,t=0$.\\
{\bf Case 4:}  Now if $s=t=1$, then the only thing that changes are the $\a_{n-1}$ and $\a_{n}$ spots, which we have already considered.\\
For a complete list of the cohomology classes that appear refer to Appendix B.\\

\subsection{The Exceptional cases}
If $\Phi$ is one of the exceptional root systems (i.e., $E_6,E_7,E_8$,$F_4$, or  $G_2$), then determining the $U_1$ cohomology reduces to looking at finitely many cases.  To do this, a program in GAP was written to calculate all different weights that satisfy equations (\ref{An:1})-(\ref{An:3}).  Then calculating all the differentials of the cohomology classes that satisfy the possible weights to see if they were in the kernel.  Note that if $\Phi=E_7$, this is the only case where the index of connection is even.  The program must now be run twice once when $t=0$ and the other time when $t=1$ with $\omega_2=\frac{1}{2}(4\a_1+7\a_2+8\a_3+12\a_4+9\a_5+6\a_6+3\a_7)$ as a generator.  \begin{rem}
Note, that the long simple roots in type $F_4$ are $\a_3, \a_4$ and the long simple root in type $G_2$ is $\a_2$.
\end{rem}
\begin{rem}
A complete list of the possible weights and classes for the exceptional cases is seen in Appendix B.
\end{rem}

\section{$B$--Cohomology}\label{Bcohom}
\subsection{}
In order to calculate $\h^2(B, \la)$, we first had to calculate $\h^2(B_r, \la)$.  For full details on the steps taken to calculate $\h^2(B_r, \la)$  are in Appendix C.\\
\\
Cline, Parshall, and Scott \cite{CPS} gives a relationship between the $B_r$-cohomology and the $B$-cohomology: $$\displaystyle{\h^2(B,\la) \cong \lim_{\longleftarrow } \h^2(B_r,\la)}.$$

\begin{thm}\label{thm:B}
Let $p=2$ and $\la \in X(T)$.
\begin{itemize}
\item[(a)] If $\Phi$ is simply laced, then
\[ \h^2(B,\la) \cong \left\{ \begin{array}{lll}
k && \mbox{if } \la=p^lw \cdot 0, \mbox{ with } 0 \leq l,\; l(w)=2,\\
k && \mbox{if } \la=-p^l\a, \mbox{ with } 0 < l \mbox{ and } \a \in \Delta,\\
k && \mbox{if } \la=-p^t\b-p^l\a, \mbox{ with } 0 \leq l < t \mbox{ and } \a,\b \in \Delta,\\
0 && \mbox{else.}\\
\end{array} \right. \]
\item[(b)] If $\Phi$ is of type $B_n$, then
\[ \h^2(B,\la) \cong \left\{ \begin{array}{lll}
k && \mbox{if } \la=p^lw \cdot 0, \mbox{ with } 0 \leq l,\; l(w)=2,\\
k && \mbox{if } \la=-p^l\a, \mbox{ with } 0 < l \mbox{ and } \a \in \Delta,\\
k && \mbox{if } \la=-p^t\b-p^l\a, \mbox{ with } 0 \leq l < t \mbox{ and } \a,\b \in \Delta,\\
k && \mbox{if } \la=-p^{l+1}(\a_{n-1}+\a_n) - p^l\a_{n-1} \mbox{ with } 0 \leq l,\\
0 && \mbox{else.}\\
\end{array}\right. \]
\item[(c)] If $\Phi$ is of type $C_n$, then
\[ \h^2(B,\la) \cong \left\{ \begin{array}{lll}
k && \mbox{if } \la=p^lw \cdot 0, \mbox{ with } 0 \leq l, \;l(w)=2,\\
k && \mbox{if } \la=-p^l\a, \mbox{ with } 0 < l \mbox{ and } \a \in \Delta,\\
k && \mbox{if } \la=-p^t\b-p^l\a, \mbox{ with } 0 \leq l < t \mbox{ and } \a,\b \in \Delta,\\
k && \mbox{if } \la=-p^{l+1}(\a_{n-1}+\a_n) \mbox{ with } 0 \leq l,\\
0 && \mbox{else.}\\
\end{array} \right. \]
\item[(d)] If $\Phi$ is of type $F_4$, then
\[ \h^2(B,\la) \cong \left\{ \begin{array}{lll}
k && \mbox{if } \la=p^lw \cdot 0, \mbox{ with } 0 \leq l,\; l(w)=2,\\
k && \mbox{if } \la=-p^l\a, \mbox{ with } 0 < l \mbox{ and } \a \in \Delta,\\
k && \mbox{if } \la=-p^t\b-p^l\a, \mbox{ with } 0 \leq l < t \mbox{ and } \a,\b \in \Delta,\\
k && \mbox{if } \la=-p^{l+1}(\a_3+\b) - p^l\a_{2} \mbox{ with } 0 \leq l \mbox{ and } \b \in \{\a_2,\a_4\},\\
0 && \mbox{else.}\\
\end{array} \right. \]
\item[(e)] If $\Phi$ is of type $G_2$, then
\[ \h^2(B,\la) \cong \left\{ \begin{array}{lll}
k && \mbox{if } \la=p^lw \cdot 0, \mbox{ with } 0 \leq l, \;l(w)=2,\\
k && \mbox{if } \la=-p^l\a, \mbox{ with } 0 < l \mbox{ and } \a \in \Delta,\\
k && \mbox{if } \la=-p^t\b-p^l\a, \mbox{ with } 0 \leq l < t \mbox{ and } \a,\b \in \Delta,\\
k && \mbox{if } \la=-p^{l+1}(\a_1+\a_2) - p^l\a_2 \mbox{ with } 0 \leq l,\\
0 && \mbox{else.}\\
\end{array} \right. \]
\end{itemize}
\end{thm}
\begin{proof}
Assume that $\la \in X(T)$ with $\h^2(B,\la) \not= 0$, then $\la \not= 0$.  Choose $s > 0$ such that
\begin{itemize}
\item[(i)] The natural map $\h^n(B,\la) \rightarrow \h^n(B_r,\la)$ is nonzero for all $r \geq s$.

\item[(ii)] By choosing a possibly larger $s$, we can further assume that $|\left<\la,\a^\vee\right>|<p^{s-1}$ for all $\a \in \Delta$.
\end{itemize}

The $\h^2(B_r, \la)$ theorem and condition (ii) gives us $\h^2(B_r,\la)$ is one-dimensional for all $r \geq s$.  Furthermore, since $\h^2(B,\la)$ has trivial $B$-action, then $\h^2(B_r, \la) \cong k$ for all $r \geq s$, by condition (i).\\
On the other hand, if there exists an integer $s$ such that $\h^2(B_r,\la) \cong k$ for all $r \geq s$, then $\displaystyle{\h^2(B,\la) \cong \lim_{\longleftarrow } \h^2(B_r,\la) \cong k}$.  
\end{proof}

\section{$G_r$-Cohomology}
\subsection{}
 The $B_r$-cohomology can aslo be used to calculate $\h^2(G_r, H^0(\la))$ for $\la \in X(T)_+$.  Recall the following theorem from \cite[6.1]{BNP2}.
\begin{thm}\label{thm:ind} Let $\la \in X(T)_+$ and $p$ be an arbitrary prime.  Then
$$\h^2(G_r,H^0(\la))^{(-r)} \cong \operatorname{ind}_B^G\h^2(B_r,\la)^{(-r)}.$$
\end{thm}

Using the isomorphism in Theorem \ref{thm:ind} and the $B_r$ cohomology results (Theorem 5.11) then we get the following theorem.

\begin{thm} Let $p=2$ and $\la \in X(T)_+$.  Then $\h^2(G_r, H^0(\la))$ is given by the following table
\end{thm}
\input{G_rtable}

  \begin{rem}
The module, $M$, is taken to be the appropriate module:
\begin{itemize}
\item $M_{B_n}$ is the module with factors $\a_n, k$.\\
\item $M_{C_n}$ is the module with factors $\a_{n-1},k$.\\
\item $M_{F_4}$ is the module with factors $\a_3, k$.\\
\item $M_{G_2}$ is the module with factors $\a_1,k$.
\end{itemize}
 \end{rem}
 
 If the $B_r$-cohomology involves an indecomposable $B$-module, then it is necessary to determine the structure of the modules for the $G_r$-cohomology, separately.  To see the structure of the modules, refer to Appendix C, found at the ArXiv.

\section*{Acknowledgements}

The author wishes to acknowledge the help and support of her Ph.D. advisor, Daniel K. Nakano.  This work was part of the author's doctoral thesis at The University of Georgia.  This work partially supported by the VIGRE fellowship at University of Georgia.\\

\appendix

\input{appendix1}

\input{appendix3}

\input{Brcohomology}

\input{appendix2}

\end{document}

%% file: G_rtable.tex
\nopagebreak
\begin{longtable}{|l|l|l|}
\hline 
\large{\bf{$\h^2(G_r, H^0(\la)), \; p=2$}} & Root System & \large{\bf{$\la \in X(T)_+$}} \\
\hline\hline 
$\operatorname{ind}_B^G\; \h^2(B_1, w \cdot 0 + p \v)^{(r-1)}$ & $A_n (n\not=3), B_n,$  & $p^{r-1}(w \cdot 0 + p \v), l(w)=0,2$\\
& $E_n, F_4, G_2$  & \\
\hline
$\operatorname{ind}_B^G \h^2(B_1, w \cdot 0 + p \v)^{(r-1)}$ & $A_3$ & $p^r\v$ \\
\hline
$\operatorname{ind}_B^G \; \h^2(B_1, w \cdot 0 + p \v)^{(r)}$ & $C_n, D_n$ & $p^{r-1}(w \cdot 0 + p \v), l(w)=0,2$ \\
\hline
$H^0(\v)^{(r)}$ & All & $p^r\v + p^l w \cdot 0$ with $l(w)=2$ \\
& & $\quad$ and $0\leq l < r-1$ \\
\hline
$H^0(\v)^{(r)}$ & $A_n, D_n, E_6$ & $p^r\v-p^l \a$ with $0 < l < r,$ \\
 & \quad $E_7, E_8$ & \quad $\a \in \Delta$\\
\hline
$H^0(\v)^{(r)}$ & $B_3, B_4, F_4, G_2$ & $p^r\v-p^l \a$ with $0 \leq l < r,$ \\
& & \quad $\a \in \Delta$ \\
\hline
$H^0(\v)^{(r)}$ & $B_n, n \geq 5$ & $p^r\v-p^l \a$ with $0 \leq l < r,$ \\
& & $\a \in \Delta, l \not= r-1$ if $\a=\a_{n-1}$ \\
\hline
$H^0(\v)^{(r)}$ & $C_n$ & $p^r\v-p^l \a$ with $0 \leq l < r,$  \\
 & & $\quad \a \in \Delta; l \not= r-1$ if $\a=\a_n$ \\
\hline
$H^0(\v)^{(r)}$ & All & $p^r\v-p^t\b-p^l \a$ with $\a, \b \in \Delta$ \\
& & $\quad 0 \leq l < t < r,$ \\
\hline
$H^0(\v)^{(r)}$ & $A_3$ & $ p^r \v+p^{r-1} \a_2 - p^l \a$ with\\
& & $\quad 0 \leq l < r-1, \; \a \in \Delta$ \\
\hline
$H^0(\v + \w_1)^{(r)} \oplus H^0(\v + \w_3)^{(r)}$ & $A_3$ & $p^r\v+p^{r-1}\w_2-p^l\a$, with  \\ 
& & $\quad 0 \leq l < r-1, \; \a \in \Delta$ \\
\hline
$H^0(\v + \w_1)^{(r)} \oplus H^0(\v + \w_3)^{(r)}$ & $B_3$ & $p^r\v+p^{r-1}\a_2-p^l\a$, with  \\
& & $\quad 0 \leq l < r-1, \; \a \in \Delta$ \\
\hline
$H^0(\v)^{(r)}$ & $B_n$ & $p^r\v -p^{l+1}(\a_{n-1}+\a_n)-p^l\a_2$  \\
& & \quad with $ 0 \leq l < r-1$ \\
\hline
$\operatorname{ind}_B^G(M_{B_n} \t \v)^{(r)}$ & $B_n, n \not= 4$ & $p^r\v-p^{r-1}\a_{n-1}-p^l\a$ with  \\
& & $\quad 0 \leq l < r-1, \; \a \in \Delta$ \\
\hline
$\operatorname{ind}_B^G((\w_1+M_{B_4}) \t \v)^{(r)}$ & $B_4$ & $p^r\v-p^{r-1}\a_i-p^l\a$ with $\a \in \Delta,$ \\
& & $\quad 0 \leq l < r-1, \; i \in \{1,3\}$ \\
\hline
$\operatorname{ind}_B^G (M_{B_n} \t \v)^{(r)}$ & $B_n$ & $p^r\v-p^{r-1}\a_{n-1} - p^l\a_n$  \\
& & \quad with $0 \leq l < r-1$ \\
\hline
$H^0(\v)^{(r)}$ & $C_n$ & $p^r\v-p^l(\a_{n-1}+\a_n)$ with \\
& & $\quad 0 \leq l < r-1$\\
\hline
$\operatorname{ind}_B^G(M \t \v)^{(r)}$ & $C_n$ & $p^r\v - p^{r-1}\a_n - p^l\a$ with  \\
& & $\quad 0 \leq l < r-1, \; \a \in \Delta$ \\
\hline
$\operatorname{ind}_B^G(M \t \v)^{(r)}$ & $C_n$ & $p^r\v-p^{r-1}\a$ where\\
& &  $\quad \a \in \{\a_{n-1}, \a_n\}$ \\
\hline
$\operatorname{ind}_B^G \; \h^1(B_{r-1},M^{(-1)} \t \la_1)$ & $C_n$ & $p\la_1$ \\
 $\quad  \oplus \operatorname{ind}_B^G \; \h^2(B_{r-1}, \la_1)$ & &\\ 
\hline
$H^0(\v)^{(r)}$ & $D_4$ & $p^r\v+p^{r-1}\a_2-p^l\a$ with\\
& &$\quad  0 \leq l < r-1, \; \a \in \Delta$ \\
\hline
$H^0(\v+\w_1)^{(r)} \oplus H^0(\v + \w_3)^{(r)}$ & $D_4$ & $p^r\v+p^{r-1}\w_2-p^l\a$ with \\
$\quad  \oplus H^0(\v+\w_4)^{(r)}$ & & $\quad  0 \leq l < r-1, \; \a \in \Delta$ \\
\hline
$H^0(\v)^{(r)}$ & $D_n,  (n \geq 5)$ & $p^r\v+p^{r-1}\a_i-p^l\a$ with\\
& & $\quad 0 \leq l < r-1,\;  \a \in \Delta, \; $ \\
& & \quad $i \not= n-1,n$ \\
\hline
$H^0(\v+\w_{n-1})^{(r)} $ & $D_n, (n \geq 5)$ & $p^r\v+p^{r-1}\w_{n-2}-p^l\a$ with  \\
$\quad \oplus H^0(\v+\w_{n})^{(r)}$ & & $\quad  0 \leq l < r-1, \; \a \in \Delta$ \\
\hline
$H^0(\v)^{(r)}$ & $F_4$ & $p^r\v-p^{l+1}(\a_3+\b)-p^l\a_2$ with \\
& & $\quad  0 \leq l < r-1, \; \b \in \{\a_2, \a_4\}$ \\
\hline
$\operatorname{ind}_B^G (M_{F_4} \t \v)^{(r)}$ & $F_4$ & $p^r\v-p^{r-1}\a_2-p^l\a$ with  \\
\nopagebreak
& & $\quad 0 \leq l < r-1, \; \a \in \Delta$ \\
\hline
$\operatorname{ind}_B^G(M \t \v)^{(r)}$ & $F_4$ & $p^r\v - p^{r-1}\a_4 - p^l\a_2$ with \\
\nopagebreak
& & $\quad 0 \leq l < r-1$ \\
\hline
$H^0(\v)^{(r)}$ & $G_2$ & $p^r\v-p^{l+1}(\a_1+\a_2) - p^l\a_2$  \\
& & with $\quad 0 \leq l < r-1,\; \a \in \Delta$ \\
\hline
$\operatorname{ind}_B^G(M_{G_2} \t \v)^{(r)}$ & $G_2$ & $p^r\v-p^{r-1}\a_2 - p^l\a$ with  \\
& & $\quad 0 \leq l < r-1, \; \a \in \Delta$ \\
\hline

\end{longtable}

%% file: appendix1.tex
\section{Cohomology Classes of $\h^2(U_1,k)$}

\begin{longtable}{|l|l|l|}
\caption{Cohomology Classes}
\endfirsthead

\hline
\endhead

\hline
\endfoot

\hline
$\Phi$ & $w$ & cohomology class\\
\hline
$A_n$ & $s_{\a_i}s_{\a_{i+2}}$ & $\p_{\a_{i+1}}\t\p_{\a_i+\a_{i+1}+\a_{i+2}} + \p_{\a_i+\a_{i+1}}\t\p_{\a_{i+1}+\a_{i+2}}$\\
\nopagebreak
\cline{3-3}
\nopagebreak
& \scriptsize{$1 \leq i \leq n-2$} & $\p_{\a_i}\t\p_{\a_{i+2}}$\\
\hline
\pagebreak[3]
\hline
& $s_{\a_i}s_{\a_{i+2}}$ & $\p_{\a_{i+1}}\t\p_{\a_i+\a_{i+1}+\a_{i+2}} + \p_{\a_i+\a_{i+1}}\t\p_{\a_{i+1}+\a_{i+2}}$\\
\nopagebreak
\cline{3-3}
\nopagebreak
& \scriptsize{$1 \leq i \leq n-2$} & $\p_{\a_i}\t\p_{\a_{i+2}}$\\
\cline{2-3}
$B_n$ & $s_{\a_i}s_{\a_{n-1}}$ & $\p_{\a_i}\t\p_{\a_{n-1}+2\a_n}$\\
\nopagebreak
\cline{3-3}
\nopagebreak
& \scriptsize{$1 \leq i \leq n-3$} & $\p_{\a_i}\t\p_{\a_{n-1}}$\\
\pagebreak[3]
\cline{2-3}
\pagebreak[3]
& $e$ & $\p_{\a_i}\t\p_{\a_i+2\a_{i+1}+2\a_{i+2} + \ldots + 2\a_n} + \p_{\a_i + \a_{i+1}}\t\p_{\a_i+\a_{i+1}+2\a_{i+2} + \ldots + 2\a_n}$\\
& & \quad $+ \ldots + \p_{\a_i + \a_{i+1} \ldots + \a_{n-1}}\t\p_{\a_i + \a_{i+1} \ldots + \a_{n-1} + 2\a_n}$\\
& & \quad $+ \p_{\a_i + \a_{i+1} \ldots + \a_{n-1}+\a_n}\t\p_{\a_i + \a_{i+1} \ldots + \a_{n-1}+\a_n}$\\
\hline
\pagebreak[3]
\hline
 & $s_{\a_i}s_{\a_{i+2}}$ & $\p_{\a_{i+1}}\t\p_{\a_i+\a_{i+1}+\a_{i+2}} + \p_{\a_i+\a_{i+1}}\t\p_{\a_{i+1}+\a_{i+2}}$\\
\nopagebreak
\cline{3-3}
\nopagebreak
& \scriptsize{$1 \leq i \leq n-3$} & $\p_{\a_i}\t\p_{\a_{i+2}}$\\
\pagebreak[3]
\cline{2-3}
\pagebreak[3]
& $s_{\a_{n-2}}s_{\a_{n}}$ & $\p_{\a_{n-2}}\t\p_{\a_n}$\\
\pagebreak[3]
\cline{2-3}
\pagebreak[3]
$C_n$ & $s_{\a_{n-1}}s_{\a_n}$ & $\p_{\a_{n-1}}\t\p_{2\a_{n-1}+\a_n}$\\
\pagebreak[3]
\cline{2-3}
\pagebreak[3]
& $s_{\a_i}s_{\a_n}$ & $\p_{\a_i}\t\p_{2\a_{n-1}+\a_n}$\\
\nopagebreak
\cline{3-3}
\nopagebreak
& \scriptsize{$1 \leq i \leq n-3$} & $\p_{\a_i}\t\p_{\a_n}$\\ 
\pagebreak[3]
\cline{2-3}
\pagebreak[3]
& $e$ & $\p_{\a_n}\t\p_{2\a_i + 2\a_{i+1} + \ldots + 2\a_{n-1} + \a_n} + \p_{\a_i + \ldots + \a_n}\t\p_{\a_i + \ldots + \a_n}$\\
\hline
\pagebreak[3]
\hline
& $s_{\a_i}s_{\a_{i+2}}$ & $\p_{\a_{i+1}}\t\p_{\a_i+\a_{i+1}+\a_{i+2}} + \p_{\a_i+\a_{i+1}}\t\p_{\a_{i+1}+\a_{i+2}}$\\
\nopagebreak
\cline{3-3}
\nopagebreak
 & \scriptsize{$1 \leq i \leq n-4$} & $\p_{\a_i}\t\p_{\a_{i+2}}$\\
\pagebreak[3]
\cline{2-3}
\pagebreak[3]
& & $\p_{\a_{n}}\t\p_{\a_{n-3}+2\a_{n-2}+\a_{n-1}+\a_n} + \p_{\a_{n-2}+\a_{n}}\t\p_{\a_{n-3}+\a_{n-2}+\a_{n-1}+\a_n}$\\
& & \quad $ + \p_{\a_{n-3}+\a_{n-2}+\a_{n}}\t\p_{\a_{n-2}+\a_{n-1}+\a_n}$\\
\nopagebreak
\cline{3-3}
\nopagebreak
 & $s_{\a_{n-3}}s_{\a_{n-1}}$ & $\p_{\a_{n-2}}\t\p_{\a_{n-3}+\a_{n-2}+\a_{n-1}} + \p_{\a_{n-3}+\a_{n-2}}\t\p_{\a_{n-2}+\a_{n-1}}$\\
\nopagebreak
\cline{3-3}
\nopagebreak
& & $\p_{\a_{n-3}}\t\p_{\a_{n-1}}$\\
\pagebreak[3]
\cline{2-3}
\pagebreak[3]
& & $\p_{\a_{n-1}}\t\p_{\a_{n-3}+2\a_{n-2}+\a_{n-1}+\a_n} + \p_{\a_{n-2}+\a_{n-1}}\t\p_{\a_{n-3}+\a_{n-2}+\a_{n-1}+\a_n}$\\
& & \quad $ + \p_{\a_{n-3}+\a_{n-2}+\a_{n-1}}\t\p_{\a_{n-2}+\a_{n-1}+\a_n}$\\
\nopagebreak
\cline{3-3}
\nopagebreak
$D_n$ & $s_{\a_{n-3}}s_{\a_{n}}$ & $\p_{\a_{n-2}}\t\p_{\a_{n-3}+\a_{n-2}+\a_n} + \p_{\a_{n-3}+\a_{n-2}}\t\p_{\a_{n-2}+\a_n}$\\
\nopagebreak
\cline{3-3}
\nopagebreak
& & $\p_{\a_{n-3}}\t\p_{\a_n}$\\
\pagebreak[3]
\cline{2-3}
\pagebreak[3]
& & $\p_{\a_i}\t\p_{\a_i+2\a_{i+1}+\ldots + 2\a_{n-2} + \a_{n-1}+\a_n} +$\\
& & \quad $\p_{\a_i+\a_{i+1}}\t\p_{\a_i+\a_{i+1}+2\a_{i+2} + \ldots + 2\a_{n-2} + \a_{n-1} + \a_n}$\\
& $s_{\a_{n-1}}s_{\a_n}$ & \quad $+ \ldots + \p_{\a_i+\ldots+\a_{n-1}}\t\p_{\a_i+\ldots+\a_{n-2}+\a_n}, \quad$ \scriptsize{$1 \leq i \leq n-3$}\\
\nopagebreak
\cline{3-3}
\nopagebreak
&  & $\p_{\a_{n-1}}\t\p_{\a_n}$\\
\pagebreak[3]
\cline{2-3}
\pagebreak[3]
\hline
\pagebreak[3]
\hline
& $s_{\a_1}s_{\a_4}$ & $\p_{\a_3}\t\p_{\a_1+\a_3+\a_4} + \p_{\a_1+\a_3}\t\p_{\a_3+\a_4}$\\
\nopagebreak
\cline{3-3}
\nopagebreak
& & $\p_{\a_1}\t\p_{\a_4}$\\
\pagebreak[3]
\cline{2-3}
\pagebreak[3]
& $s_{\a_4}s_{\a_6}$ & $\p_{\a_5}\t\p_{\a_4+\a_5+\a_6} + \p_{\a_4+\a_5}\t\p_{\a_5+\a_6}$\\
\nopagebreak
\cline{3-3}
\nopagebreak
& & $\p_{\a_4}\t\p_{\a_6}$\\
\pagebreak[3]
\cline{2-3}
\pagebreak[3]
& & $\p_{\a_6}\t\p_{\a_2+\a_3+2\a_4+2\a_5+\a_6} + \p_{\a_5+\a_6}\t\p_{\a_2+\a_3+2\a_4+\a_5+\a_6}$\\
& & \quad $\p_{\a_4+\a_5+\a_6}\t\p_{\a_2+\a_3+\a_4+\a_5+\a_6} + \p_{\a_2+\a_4+\a_5+\a_6}\t\p_{\a_3+\a_4+\a_5+\a_6}$\\
\nopagebreak
\cline{3-3}
\nopagebreak
$E_6$  & $s_{\a_2}s_{\a_3}$ & $\p_{\a_5}\t\p_{\a_2+\a_3+2\a_4+\a_5} + \p_{\a_4+\a_5}\t\p_{\a_2+\a_3+\a_4+\a_5}$\\
& & \quad $+ \p_{\a_2+\a_4+\a_5}\t\p_{\a_3+\a_4+\a_5}$\\
\nopagebreak
\cline{3-3}
\nopagebreak
 & & $\p_{\a_4}\t\p_{\a_2+\a_3+\a_4} + \p_{\a_2+\a_4}\t\p_{\a_3+\a_4}$\\
\nopagebreak
\cline{3-3}
\nopagebreak
& & $\p_{\a_2}\t\p_{\a_3}$\\
\pagebreak[3]
\cline{2-3}

&  & $\p_{\a_1}\t\p_{\a_1+\a_2+2\a_3+2\a_4+\a_5} + \p_{\a_1+\a_3}\t\p_{\a_1+\a_2+\a_3+2\a_4+\a_5}$\\
\nopagebreak
& & \quad $\p_{\a_1+\a_3+\a_4}\t\p_{\a_1+\a_2+\a_3+\a_4+\a_5} + \p_{\a_1+\a_2+\a_3+\a_4}\t\p_{\a_1+\a_3+\a_4+\a_5}$\\
\nopagebreak[3]
\cline{3-3}
\nopagebreak[3]
& $s_{\a_2}s_{\a_5}$ & $\p_{\a_3}\t\p_{\a_2+\a_3+2\a_4+\a_5} + \p_{\a_3+\a_4}\t\p_{\a_2+\a_3+\a_4+\a_5}$\\
\nopagebreak[3]
& & \quad $\p_{\a_2+\a_3+\a_4}\t\p_{\a_3+\a_4+\a_5}$\\
\nopagebreak[3]
\cline{3-3}
\nopagebreak[3]
$E_6$ & & $\p_{\a_4}\t\p_{\a_2+\a_4+\a_5} + \p_{\a_2+\a_4}\t\p_{\a_4+\a_5}$\\
\nopagebreak
\cline{3-3}
\nopagebreak
& & $\p_{\a_2}\t\p_{\a_5}$\\
\pagebreak[3]
\cline{2-3}
\pagebreak[3]
& & $\p_{\a_2}\t\p_{\a_2+\a_3+2\a_4+\a_5} + \p_{\a_2+\a_4}\t\p_{\a_2+\a_3+\a_4+\a_5} $\\
& & $ \quad + \p_{\a_2+\a_3+\a_4}\t\p_{\a_2+\a_4+\a_5}$\\
\nopagebreak
\cline{3-3}
\nopagebreak
& $s_{\a_3}s_{\a_5}$ & $\p_{\a_4}\t\p_{\a_3+\a_4+\a_5} + \p_{\a_3+\a_4}\t\p_{\a_4+\a_5}$\\
\nopagebreak
\cline{3-3}
\nopagebreak
& & $\p_{\a_3}\t\p_{\a_5}$\\
\hline
\pagebreak[3]
\hline
& $s_{\a_1}s_{\a_4}$ & $\p_{\a_3}\t\p_{\a_1+\a_3+\a_4}$\\
\nopagebreak
\cline{3-3}
\nopagebreak
& & $\p_{\a_1}\t\p_{\a_4}$\\
\pagebreak[3]
\cline{2-3}
\pagebreak[3]
& $s_{\a_4}s_{\a_6}$ & $\p_{\a_5}\t\p_{\a_4+\a_5+\a_6} + \p_{\a_4+\a_5}\t\p_{\a_5+\a_6}$\\
\nopagebreak
\cline{3-3}
\nopagebreak
& & $\p_{\a_4}\t\p_{\a_6}$\\
\pagebreak[3]
\cline{2-3}
\pagebreak[3]
 & $s_{\a_5}s_{\a_7}$ & $\p_{\a_6}\t\p_{\a_5+\a_6+\a_7} + \p_{\a_5+\a_6}\t\p_{\a_6+\a_7}$\\
\nopagebreak
\cline{3-3}
\nopagebreak
& & $\p_{\a_5}\t\p_{\a_7}$\\
\pagebreak[3]
\cline{2-3}
\pagebreak[3]
& & $\p_{\a_1}\t\p_{\a_1+\a_2+2\a_3+2\a_4+\a_5} + \p_{\a_1+\a_3}\t\p_{\a_1+\a_2+\a_3+2\a_4+\a_5}$\\
\nopagebreak
& & \quad $\p_{\a_1+\a_3+\a_4}\t\p_{\a_1+\a_2+\a_3+\a_4+\a_5} + \p_{\a_1+\a_2+\a_3+\a_4}\t\p_{\a_1+\a_3+\a_4+\a_5}$\\
\nopagebreak
\cline{3-3}
$E_7$ & $s_{\a_2}s_{\a_5}$ & $\p_{\a_3}\t\p_{\a_2+\a_3+2\a_4+\a_5} + \p_{\a_3+\a_4}\t\p_{\a_2+\a_3+\a_4+\a_5}$\\
& & \quad $+ \p_{\a_2+\a_3+\a_4}\t\p_{\a_3+\a_4+\a_5}$\\
\nopagebreak
\cline{3-3}
\nopagebreak
& & $\p_{\a_4}\t\p_{\a_2+\a_4+\a_5} + \p_{\a_2+\a_4}\t\p_{\a_4+\a_5}$\\
\nopagebreak
\cline{3-3}
\nopagebreak
& & $\p_{\a_2}\t\p_{\a_5}$\\

\pagebreak[2]
\cline{2-3}
\pagebreak[2]

& & $\p_{\a_2}\t\p_{\a_2+\a_3+2\a_4+\a_5} + \p_{\a_2+\a_4}\t\p_{\a_2+\a_3+\a_4+\a_5}$\\
& & \quad $+ \p_{\a_2+\a_3+\a_4}\t\p_{\a_2+\a_4+\a_5}$\\
\pagebreak[3]
\cline{3-3}
\pagebreak[3]
& $s_{\a_3}s_{\a_5}$ & $\p_{\a_4}\t\p_{\a_3+\a_4+\a_5} + \p_{\a_3+\a_4}\t\p_{\a_4+\a_5}$\\
\nopagebreak
\cline{3-3}
\nopagebreak
&  & $\p_{\a_3}\t\p_{\a_5}$\\
\nopagebreak
\cline{2-3}
\nopagebreak
& & $\p_{\a_7}\t\p_{\a_2+\a_3+2\a_4+2\a_5+2\a_6+\a_7} + \p_{\a_6+\a_7}\t\p_{\a_2+\a_3+2\a_4+2\a_5+\a_6+\a_7}$\\
\nopagebreak
& & \quad $+ \p_{\a_5+\a_6+\a_7}\t\p_{\a_2+\a_3+2\a_4+\a_5+\a_6+\a_7}$\\
\nopagebreak
& & \quad $+ \p_{\a_4+\a_5+\a_6+\a_7}\t\p_{\a_2+\a_3+\a_4+\a_5+\a_6+\a_7}$\\
\nopagebreak
& & \quad $+ \p_{\a_2+\a_4+\a_5+\a_6}\t\p_{\a_3+\a_4+\a_5+\a_6+\a_7}$\\
\nopagebreak
\cline{3-3}
\nopagebreak
 & $s_{\a_2}s_{\a_3}$ & $\p_{\a_6}\t\p_{\a_2+\a_3+2\a_4+2\a_5+\a_6} + \p_{\a_5+\a_6}\t\p_{\a_2+\a_3+2\a_4+\a_5+\a_6}$\\
\nopagebreak
& & \quad $\p_{\a_4+\a_5+\a_6}\t\p_{\a_2+\a_3+\a_4+\a_5+\a_6} + \p_{\a_2+\a_4+\a_5+\a_6}\t\p_{\a_3+\a_4+\a_5+\a_6}$\\
\nopagebreak
\cline{3-3}
\nopagebreak
& & $\p_{\a_5}\t\p_{\a_2+\a_3+2\a_4+\a_5} + \p_{\a_4+\a_5}\t\p_{\a_2+\a_3+\a_4+\a_5}$\\
& & \quad $+ \p_{\a_2+\a_4+\a_5}\t\p_{\a_3+\a_4+\a_5}$\\
\nopagebreak
\cline{3-3}
\nopagebreak
& & $\p_{\a_4}\t\p_{\a_2+\a_3+\a_4} + \p_{\a_2+\a_4}\t\p_{\a_3+\a_4}$\\
\nopagebreak
\cline{3-3}
\nopagebreak
& & $\p_{\a_2}\t\p_{\a_3}$\\
\hline
\pagebreak[3]
\hline
& $s_{\a_1}s_{\a_4}$& $\p_{\a_3}\t\p_{\a_1+\a_3+\a_4} + \p_{\a_1+\a_3}\t\p_{\a_3+\a_4}$\\
\nopagebreak
\cline{3-3}
\nopagebreak
& & $\p_{\a_1}\t\p_{\a_4}$\\
\pagebreak[3]
\cline{2-3}
\pagebreak[3]
& $s_{\a_4}s_{\a_6}$ & $\p_{\a_5}\t\p_{\a_4+\a_5+\a_6} + \p_{\a_4+\a_5}\t\p_{\a_5+\a_6}$\\
\nopagebreak
\cline{3-3}
\nopagebreak
& & $\p_{\a_4}\t\p_{\a_6}$\\
\pagebreak[3]
\cline{2-3}
\pagebreak[3]
& $s_{\a_5}s_{\a_7}$ & $\p_{\a_6}\t\p_{\a_5+\a_6+\a_7} + \p_{\a_5+\a_6}\t\p_{\a_6+\a_7}$\\
\nopagebreak
\cline{3-3}
\nopagebreak
& & $\p_{\a_5}\t\p_{\a_7}$\\
\pagebreak[3]
\cline{2-3}
\pagebreak[3]
$E_8$ & $s_{\a_6}s_{\a_8}$ & $\p_{\a_7}\t\p_{\a_6+\a_7+\a_8} + \p_{\a_6+\a_7 }\t\p_{\a_7+\a_8}$\\
\nopagebreak
\cline{3-3}
\nopagebreak
& & $\p_{\a_6}\t\p_{\a_8}$\\
\pagebreak[2]
\cline{2-3}
\pagebreak[2]
& & $\p_{\a_1}\t\p_{\a_1+\a_2+2\a_3+2\a_4+\a_5} + \p_{\a_1+\a_3}\t\p_{\a_1+\a_2+\a_3+2\a_4+\a_5}$\\
\nopagebreak
& & \quad $+ \p_{\a_1+\a_3+\a_4}\t\p_{\a_1+\a_2+\a_3+\a_4+\a_5} + \p_{\a_1+\a_2+\a_3+\a_4}\t\p_{\a_1+\a_3+\a_4+\a_5}$\\
\nopagebreak
\cline{3-3}
& $s_{\a_2}s_{\a_5}$ & $\p_{\a_3}\t\p_{\a_2+\a_3+2\a_4+\a_5} + \p_{\a_3+\a_4}\t\p_{\a_2+\a_3+\a_4+\a_5} + \p_{\a_2+\a_3+\a_4}\t\p_{\a_3+\a_4+\a_5}$\\
\nopagebreak
\cline{3-3}
\nopagebreak
& & $\p_{\a_4}\t\p_{\a_2+\a_4+\a_5} + \p_{\a_2+\a_4}\t\p_{\a_4+\a_5}$\\
\nopagebreak
\cline{3-3}
\nopagebreak
& & $\p_{\a_2}\t\p_{\a_5}$\\
\pagebreak[3]
\cline{2-3}
\pagebreak[3]
& & $\p_{\a_2}\t\p_{\a_2+\a_3+2\a_4+\a_5} + \p_{\a_2+\a_4}\t\p_{\a_2+\a_3+\a_4+\a_5} + \p_{\a_2+\a_3+\a_4}\t\p_{\a_2+\a_4+\a_5}$\\
\nopagebreak
\cline{3-3}
\nopagebreak
& $s_{\a_3}s_{\a_5}$ & $\p_{\a_4}\t\p_{\a_3+\a_4+\a_5} + \p_{\a_3+\a_4}\t\p_{\a_4+\a_5}$\\
\nopagebreak
\cline{3-3}
\nopagebreak
&  & $\p_{\a_3}\t\p_{\a_5}$\\
\cline{2-3}
& & $\p_{\a_8}\t\p_{\a_2+\a_3+2\a_4+2\a_5+2\a_6+2\a_7+\a_8} + \p_{\a_7+\a_8}\t\p_{\a_2+\a_3+2\a_4+2\a_5+2\a_6+\a_7+\a_8}$\\
\nopagebreak
& & \quad $+\p_{\a_6+\a_7+\a_8}\t\p_{\a_2+\a_3+2\a_4+2\a_5+\a_6+\a_7+\a_8}$\\
\nopagebreak
& & \quad $+ \p_{\a_5+\a_6+\a_7+\a_8}\t\p_{\a_2+\a_3+2\a_4+\a_5+\a_6+\a_7+\a_8}$\\
\nopagebreak
& & \quad $+\p_{\a_4+\a_5+\a_6+\a_7+\a_8}\t\p_{\a_2+\a_3+\a_4+\a_5+\a_6+\a_7+\a_8}$\\
\nopagebreak
& & \quad $+ \p_{\a_2+\a_4+\a_5+\a_6+\a_7+\a_8}\t\p_{\a_3+\a_4+\a_5+\a_6+\a_7+\a_8}$\\
\nopagebreak
\cline{3-3}
\nopagebreak
$E_8$ & $s_{\a_2}s_{\a_3}$ & $\p_{\a_7}\t\p_{\a_2+\a_3+2\a_4+2\a_5+2\a_6+\a_7} + \p_{\a_6+\a_7}\t\p_{\a_2+\a_3+2\a_4+2\a_5+\a_6+\a_7}$\\
\nopagebreak
& & \quad $+ \p_{\a_5+\a_6+\a_7}\t\p_{\a_2+\a_3+2\a_4+\a_5+\a_6+\a_7}$\\
\nopagebreak
& & \quad $+ \p_{\a_4+\a_5+\a_6+\a_7}\t\p_{\a_2+\a_3+\a_4+\a_5+\a_6+\a_7}$\\
\nopagebreak
& & \quad $+ \p_{\a_2+\a_4+\a_5+\a_6+\a_7}\t\p_{\a_3+\a_4+\a_5+\a_6+\a_7}$\\
\nopagebreak
\cline{3-3}
\nopagebreak
& & $\p_{\a_6}\t\p_{\a_2+\a_3+2\a_4+2\a_5+\a_6} + \p_{\a_5+\a_6}\t\p_{\a_2+\a_3+2\a_4+\a_5+\a_6}$\\
\nopagebreak
& & \quad $+ \p_{\a_4+\a_5+\a_6}\t\p_{\a_2+\a_3+\a_4+\a_5+\a_6} + \p_{\a_2+\a_4+\a_5+\a_6}\t\p_{\a_3+\a_4+\a_5+\a_6}$\\
\nopagebreak
\cline{3-3}
\nopagebreak
& & $\p_{\a_5}\t\p_{\a_2+\a_3+2\a_4+\a_5} + \p_{\a_4+\a_5}\t\p_{\a_2+\a_3+\a_4+\a_5} + \p_{\a_2+\a_4+\a_5}\t\p_{\a_3+\a_4+\a_5}$\\
\nopagebreak
\cline{3-3}
\nopagebreak
& & $\p_{\a_4}\t\p_{\a_2+\a_3+\a_4} + \p_{\a_2+\a_4}\t\p_{\a_3+\a_4}$\\
\nopagebreak
\cline{3-3}
\nopagebreak
& & $\p_{\a_2}\t\p_{\a_3}$\\

\hline
\pagebreak[3]
\hline
& & $\p_{\a_3}\t\p_{\a_1+2\a_2+2\a_3} + \p_{\a_2+\a_3}\t\p_{\a_1+\a_2+2\a_3} + \p_{\a_1+\a_2+\a_3}\t\p_{\a_2+2\a_3}$\\
\nopagebreak
\cline{3-3}
\nopagebreak
& $s_{\a_1}s_{\a_3}$ & $\p_{\a_2}\t\p_{\a_1+\a_2+\a_3} + \p_{\a_1+\a_2}\t\p_{\a_2+\a_3}$\\
\nopagebreak
\cline{3-3}
\nopagebreak
& & $\p_{\a_1}\t\p_{\a_3}$\\
\pagebreak[3]
\cline{2-3}
\pagebreak[3]
& & $\p_{\a_2}\t\p_{2\a_1+3\a_2+4\a_3+2\a_4} + \p_{\a_1+\a_2}\t\p_{\a_1+3\a_2+4\a_3+2\a_4}$\\
\nopagebreak
& & \quad $+ \p_{\a_1+2\a_2+2\a_3}\t\p_{\a_1+2\a_2+2\a_3+2\a_4} + \p_{\a_1+2\a_2+2\a_3+\a_4}\t\p_{\a_1+2\a_2+2\a_3+\a_4}$\\
\nopagebreak
\cline{3-3} 
\nopagebreak
$F_4$ & & $\p_{\a_1}\t\p_{\a_1+2\a_2+4\a_3+2\a_4} + \p_{\a_1+\a_2+2\a_3}\t\p_{\a_1+\a_2+2\a_3+2\a_4}$\\
\nopagebreak
& & \quad $+ \p_{\a_1+\a_2+2\a_3+\a_4}\t\p_{\a_1+\a_2+2\a_3+\a_4}$\\
\nopagebreak
\cline{3-3}
\nopagebreak
& $e$ & $\p_{\a_1}\t\p_{\a_1+2\a_2+2\a_3+2\a_4} + \p_{\a_1+\a_2}\t\p_{\a_1+\a_2+2\a_3+2\a_4} + \p_{\a_1+\a_2+\a_3+\a_4}$\\
\nopagebreak
\cline{3-3}
\nopagebreak
& & $\p_{\a_1}\t\p_{\a_1+2\a_2+2\a_3} + \p_{\a_1+\a_2}\t\p_{\a_1+\a_2+2\a_3} + \p_{\a_1+\a_2+\a_3}\t\p_{\a_1+\a_2+\a_3}$\\
\nopagebreak
\cline{3-3}
\nopagebreak
& & $\p_{\a_2}\t\p_{\a_2+2\a_3+2\a_4} + \p_{\a_2+\a_3+\a_4}\t\p_{\a_2+\a_3+\a_4}$ \\
\nopagebreak
\cline{3-3}
\nopagebreak
& & $\p_{\a_2}\t\p_{\a_2+2\a_3} + \p_{\a_2+\a_3}\t\p_{\a_2+\a_3}$\\
\hline
\pagebreak[3]
\hline
$G_2$ & $e$ & $\p_{\a_1}\t\p_{\a_1+2\a_2}$\\
\hline
\end{longtable}

\begin{rem} In the above table, $e$ denotes the identity element.  When $w=s_{\a_i}s_{\a_j}$ when $i$ and $j$ are not connected and aren't separated by a single vertex, then in all of the classical types, there is a single cohomology class: $\p_{\a_i}\t\p_{\a_j}$.  In the exceptional cases, there are some other classes that occur, and are demonstrated in Appendix B.
\end{rem}

%% file: appendix3.tex
\section{Weights and Cohomology Classes in the Exception Cases}

\begin{longtable}{|l|l|l|}
\caption{weights and cohomology classes in exceptional cases}
\endfirsthead

\hline
\endhead

\hline
\endfoot
\hline
$\Phi$ & cohomology class & weight\\
\hline
&$\phi_{\a_5}\t\p_{\a_2+\a_3+2\a_4+\a_5}+\phi_{\a_4+\a_5}\t\p_{\a_2+\a_3+\a_4+\a_5}$ & $(0,1,1,2,2,0)$\\
\nopagebreak
&\quad $+\phi_{\a_2+\a_4+\a_5}\t\p_{\a_3+\a_4+\a_5}$&\\
\cline{2-3}
&$\p_{\a_6}\t\p_{\a_2+\a_3+2\a_4+2\a_5+\a_6}+\p_{\a_5+\a_6}\t\p_{\a_2+\a_3+2\a_4+\a_5+\a_6}$ & $(0,1,1,2,2,2)$\\
\nopagebreak
& \quad $+\p_{\a_4+\a_5+\a_6}\t\p_{\a_2+\a_3+\a_4+\a_5+\a_6}$&\\
\nopagebreak
& \quad $+\p_{\a_2+\a_4+\a_5+\a_6}\t\p_{\a_3+\a_4+\a_5+\a_6}$&\\
\cline{2-3}
$E_6$& $\p_{\a_3}\t\p_{\a_2+\a_3+2\a_4+\a_5}+\p_{\a_3+\a_4}\t\p_{\a_2+\a_3+\a_4+\a_5}$&$(0,1,2,2,1,0)$\\
\nopagebreak
&\quad $+\p_{\a_2+\a_3+\a_4}\t\p_{\a_3+\a_4+\a_5}$&\\
\cline{2-3}
&$\p_{\a_1}\t\p_{\a_1+\a_2+2\a_3+2\a_4+\a_5} + \p_{\a_1+\a_3}\t\p_{\a_1+\a_2+\a_3+2\a_4+\a_5}$&$(2,1,2,2,1,0)$\\
\nopagebreak
&\quad $+\p_{\a_1+\a_3+\a_4}\t\p_{\a_1+\a_2+\a_3+\a_4+\a_5} + \p_{\a_1+\a_2+\a_3+\a_4}\t\p_{\a_1+\a_3+\a_4+\a_5}$&\\
\cline{2-3}
&$\p_{\a_2}\t\p_{\a_2+\a_3+2\a_4+\a_5}+\p_{\a_2+\a_4}\t\p_{\a_2+\a_3+\a_4+\a_5}$&$(0,2,1,2,1,0)$\\
\nopagebreak
&\quad $+\p_{\a_2+\a_3+\a_4}\t\p_{\a_2+\a_4+\a_5}$&\\
\hline
\pagebreak[3]
\hline
& $\p_{\a_5}\t\p_{\a_2+\a_3+2\a_4+\a_5} + \p_{\a_4+\a_5}\t\p_{\a_2+\a_3+\a_4+\a_5}$&$(0,1,1,2,2,0,0)$\\
\nopagebreak
$E_7$&\quad $+ \p_{\a_2+\a_4+\a_5}\t\p_{\a_3+\a_4+\a_5}$&\\
\cline{2-3}
& $\p_{\a_6}\t\p_{\a_2+\a_3+2\a_4+2\a_5+\a_6} + \p_{\a_5+\a_6}\t\p_{\a_2+\a_3+2\a_4+\a_5+\a_6}$&$(0,1,1,2,2,2,0)$\\
\nopagebreak
&\quad $+\p_{\a_4+\a_5+\a_6}\t\p_{\a_2+\a_3+\a_4+\a_5+\a_6}$&\\
\nopagebreak
&\quad $+ \p_{\a_2+\a_4+\a_5+\a_6}\t\p_{\a_3+\a_4+\a_5+\a_6}$&\\
\cline{2-3}
 & $\p_{\a_7}\t\p_{\a_2+\a_3+2\a_4+2\a_5+2\a_6+\a_7}$&$(0,1,1,2,2,2,2)$\\
 \nopagebreak
&\quad $+ \p_{\a_6+\a_7}\t\p_{\a_2+\a_3+2\a_4+2\a_5+\a_6+\a_7}$&\\
\nopagebreak
&\quad $+\p_{\a_5+\a_6+\a_7}\t\p_{\a_2+\a_3+2\a_4+\a_5+\a_6+\a_7}$&\\
\nopagebreak
&\quad $+ \p_{\a_4+\a_5+\a_6+\a_7}\t\p_{\a_2+\a_3+\a_4+\a_5+\a_6+\a_7}$&\\
\nopagebreak
&\quad $+\p_{\a_2+\a_4+\a_5+\a_6+\a_7}\t\p_{\a_3+\a_4+\a_5+\a_6+\a_7}$&\\
\cline{2-3}
$E_7$ & $\p_{\a_3}\t\p_{\a_2+\a_3+2\a_4+\a_5} + \p_{\a_3+\a_4}\t\p_{\a_2+\a_3+\a_4+\a_5}$&$(0,1,2,2,1,0,0)$\\
\nopagebreak
&\quad $+ \p_{\a_2+\a_3+\a_4}\t\p_{\a_3+\a_4+\a_5}$ &\\
\cline{2-3}
&$\p_{\a_1}\t\p_{\a_1+\a_2+2\a_3+2\a_4+\a_5}+\p_{\a_1+\a_3}\t\p_{\a_1+\a_2+\a_3+2\a_4+\a_5}$&$(2,1,2,2,1,0,0)$\\
\nopagebreak
&\quad $+\p_{\a_1+\a_3+\a_4}\t\p_{\a_1+\a_2+\a_3+\a_4+\a_5}$&\\
\nopagebreak
&\quad $+\p_{\a_1+\a_2+\a_3+\a_4}\t\p_{\a_1+\a_3+\a_4+\a_5}$&\\
\cline{2-3}
& $\p_{\a_2}\t\p_{\a_2+\a_3+2\a_4+\a_5} + \p_{\a_2+\a_4}\t\p_{\a_2+\a_3+\a_4+\a_5}$&$(0,2,1,2,1,0,0)$\\
\nopagebreak
&\quad $ + \p_{\a_2+\a_3+\a_4}\t\p_{\a_2+\a_4+\a_5}$&\\
\hline
\pagebreak[3]
\hline
& $\p_{\a_5}\t\p_{\a_2+\a_3+2\a_4+\a_5} + \p_{\a_4+\a_5}\t\p_{\a_2+\a_3+\a_4+\a_5}$&$(0,1,1,2,2,0,0,0)$\\
\nopagebreak
&\quad $+ \p_{\a_2+\a_4+\a_5}\t\p_{\a_3+\a_4+\a_5}$&\\
\cline{2-3}
& $\p_{\a_6}\t\p_{\a_2+\a_3+2\a_4+2\a_5+\a_6} + \p_{\a_5+\a_6}\t\p_{\a_2+\a_3+2\a_4+\a_5+\a_6}$&$(0,1,1,2,2,2,0,0)$\\
\nopagebreak
&\quad $+\p_{\a_4+\a_5+\a_6}\t\p_{\a_2+\a_3+\a_4+\a_5+\a_6}$&\\
\nopagebreak
&\quad $+ \p_{\a_2+\a_4+\a_5+\a_6}\t\p_{\a_3+\a_4+\a_5+\a_6}$&\\
\cline{2-3}
$E_8$ & $\p_{\a_7}\t\p_{\a_2+\a_3+2\a_4+2\a_5+2\a_6+\a_7}$&$(0,1,1,2,2,2,2,0)$\\
\nopagebreak
&\quad $ + \p_{\a_6+\a_7}\t\p_{\a_2+\a_3+2\a_4+2\a_5+\a_6+\a_7}$&\\
\nopagebreak
&\quad $+\p_{\a_5+\a_6+\a_7}\t\p_{\a_2+\a_3+2\a_4+\a_5+\a_6+\a_7}$&\\
\nopagebreak
&\quad $+ \p_{\a_4+\a_5+\a_6+\a_7}\t\p_{\a_2+\a_3+\a_4+\a_5+\a_6+\a_7}$&\\
\nopagebreak
&\quad $+\p_{\a_2+\a_4+\a_5+\a_6+\a_7}\t\p_{\a_3+\a_4+\a_5+\a_6+\a_7}$&\\
\cline{2-3}
& $\p_{\a_3}\t\p_{\a_2+\a_3+2\a_4+\a_5} + \p_{\a_3+\a_4}\t\p_{\a_2+\a_3+\a_4+\a_5}$&$(0,1,2,2,1,0,0,0)$\\
\nopagebreak
& \quad $+ \p_{\a_2+\a_3+\a_4}\t\p_{\a_3+\a_4+\a_5}$&\\
\cline{2-3}
& $\p_{\a_8}\t\p_{\a_2+\a_3+2\a_4+2\a_5+2\a_6+2\a_7+\a_8}$&$(0,1,1,2,2,2,2,2)$\\
\nopagebreak
&\quad $ + \p_{\a_7+\a_8}\t\p_{\a_2+\a_3+2\a_4+2\a_5+2\a_6+\a_7+\a_8}$&\\
\nopagebreak
&\quad $+\p_{\a_6+\a_7+\a_8}\t\p_{\a_2+\a_3+2\a_4+2\a_5+\a_6+\a_7+\a_8}$&\\  
\nopagebreak
&\quad $+\p_{\a_5+\a_6+\a_7+\a_8}\t\p_{\a_2+\a_3+2\a_4+\a_5+\a_6+\a_7+\a_8}$&\\
\nopagebreak
&\quad $+\p_{\a_4+\a_5+\a_6+\a_7+\a_8}\t\p_{\a_2+\a_3+\a_4+\a_5+\a_6+\a_7+\a_8}$&\\
\nopagebreak
&\quad $+ \p_{\a_2+\a_4+\a_5+\a_6+\a_7+\a_8}\t\p_{\a_3+\a_4+\a_5+\a_6+\a_7+\a_8}$&\\
\cline{2-3}
$E_8$ & $\p_{\a_1}\t\p_{\a_1+\a_2+2\a_3+2\a_4+\a_5} + \p_{\a_1+\a_3}\t\p_{\a_1+\a_2+\a_3+2\a_4+\a_5}$&$(2,1,2,2,1,0,0,0)$\\
\nopagebreak
&\quad $+ \p_{\a_1+\a_3+\a_4}\t\p_{\a_1+\a_2+\a_3+\a_4+\a_5}$&\\
\nopagebreak
&\quad $+\p_{\a_1+\a_2+\a_3+\a_4}\t\p_{\a_1+\a_3+\a_4+\a_5}$&\\
\cline{2-3}
& $\p_{\a_2}\t\p_{\a_2+\a_3+2\a_4+\a_5} + \p_{\a_2+\a_4}\t\p_{\a_2+\a_3+\a_4+\a_5}$&$(0,2,1,2,1,0,0,0)$\\
\nopagebreak
&\quad $+ \p_{\a_2+\a_3+\a_4}\t\p_{\a_2+\a_4+\a_5}$&\\
\hline
\pagebreak[3]
\hline
& $\p_{\a_2}\t\p_{\a_1+\a_2+\a_3} + \p_{\a_1+\a_2}\t\p_{\a_2+\a_3}$ & $(1,2,1,0)$\\
\cline{2-3}
& $\p_{\a_3}\t\p_{\a_1+2\a_2+2\a_3} + \p_{\a_2+\a_3}\t\p_{\a_1+\a_2+2\a_3} + \p_{\a_1+\a_2+\a_3}\t\p_{\a_2+2\a_3}$ & $(1,2,3,0)$\\
\cline{2-3}
$F_4$ & $\p_{\a_2} \t \p_{\a_2+2\a_3}$ & $(0,2,2,0)$\\
\cline{2-3}
& $\p_{\a_2} \t \p_{\a_2+2\a_3+2\a_4}$ & $(0,2,2,2)$\\
\cline{2-3}
& $\p_{\a_1}\t\p_{\a_1+2\a_2+2\a_3} + \p_{\a_1+\a_2}\t\p_{\a_1+\a_2+2\a_3}$ & $(2,2,2,0)$\\
\cline{2-3}
& $\p_{\a_1}\t\p_{\a_1+2\a_2+2\a_3+2\a_4} + \p_{\a_1+\a_2}\t\p_{\a_1+\a_2+2\a_3+2\a_4}$ & $(2,2,2,2)$\\
\cline{2-3}
& $\p_{\a_1}\t\p_{\a_1 + 2\a_2 + 4\a_3 +2\a_4} + \p_{\a_1+\a_2+2\a_3}\t\p_{\a_1+\a_2 + 2\a_3 + 2\a_4}$ & $(2,2,4,2)$\\
\cline{2-3}
& $\p_{\a_2}\t\p_{2\a_1 + 3\a_2 + 4\a_3+2\a_4} + \p_{\a_1+\a_2}\t\p_{2\a_1+2\a_2+4\a_3+2\a_4}$ & $(2,4,4,2)$\\
\nopagebreak
&\quad $ + \p_{\a_1+2\a_2+2\a_3}\t\p_{\a_1+2\a_2+2\a_3+2\a_4}$&\\
\hline
\pagebreak[3]
\hline
$G_2$ & $\p_{\a_1} \t \p_{\a_1+2\a_2}$& $(2,4)$\\
\cline{2-3}
& $\p_{\a_1}\t\p_{\a_1+3\a_2} + \p_{\a_1+\a_2}\t\p_{\a_1+2\a_2}$ & $(2,3)$\\

\end{longtable}

%% file: Brcohomology.tex
\section{Calculating $B_r$-Cohomology}

To calculate $\h^2(B_,\la)$ and $\h^2(G_r, H^0(\la))$, we must calculate $\h^2(B_1, \la)$ (as both a $T$-module and a $B$-module), and $\h^2(B_r,\la)$ first.  Recall, for $U_1 \unlhd B_1$ that $\h^2(B_1, \la)$ as a $T$-module can be calculated by the Lyndon-Hochschild-Serre (LHS) spectral sequence.  THe LHS spectral sequence, along with the fact $T_1 \cong B_1/U_1$ gives us the following characterization:
\begin{equation} \label{U,B_cohom} \h^2(B_1, \la) \cong \h^2(U_1,k)^{T_1} \cong (\h^2(U_1,k) \t \la)^{T_1}. \end{equation}
Thus it suffices to determine the $-\la$ weight space of $\h^2(U_1,k)$ relative to $T_1$:
$$\h^2(B_1,\la) \cong \h^2(U_1,k)_{-\overline{\la}}, \text{where} \overline{\la} = \{\la + p \nu \text{is a weight of} \h^2(U_1,k) | \nu \in X(T)\}.$$

\noindent \subsection{$B_1$-cohomology}
\subsubsection{}

\begin{thm}\label{thm:B1,T0}
Let $p=2$.  Then as a $T$-module, $\operatorname{H}^2(B_1,k) \cong (\mathfrak{u}^*)^{(1)}$, except in the following cases:
\begin{itemize}
\item[(a)] If $\Phi=A_3$, then $$\displaystyle{\operatorname{H}^2(B_1,k) \cong (\mathfrak{u}^*)^{(1)} \oplus (\omega_1-\omega_2+\omega_3)^{(1)} \oplus \omega_2^{(1)}}.$$
\item[(b)] If $\Phi=B_3$, then $$\displaystyle{\operatorname{H}^2(B_1,k) \cong (\mathfrak{u}^*)^{(1)} \oplus \omega_1^{(1)} \oplus (-\omega_1+\omega_2)^{(1)}}.$$
\item[(c)] If $\Phi=B_4$, then 
\indent \begin{align*}
\operatorname{H}^2(B_1,k) &\cong (\mathfrak{u}^*)^{(1)} \oplus (\omega_1-\omega_2+\omega_3-\omega_4)^{(1)} \oplus (\omega_2-\omega_4)^{(1)} \oplus (-\omega_2+\omega_3)^{(1)}\\
& \quad \oplus (-\omega_1 + \omega_2)^{(1)} \oplus \omega_1^{(1)}.
\end{align*}
\item[(d)] If $\Phi=B_n$ for $n \not= 3,4$, then
$$\displaystyle{\operatorname{H}^2(B_1,k) \cong (\mathfrak{u}^*)^{(1)} \oplus \bigoplus_{2 \leq i \leq n-1} (-\omega_{i-1} + \omega_i)^{(1)} \oplus \omega_1^{(1)}}.$$
\item[(e)] If $\Phi=C_n$, then
\indent \begin{align*}
\operatorname{H}^2(B_1,k) &\cong (\mathfrak{u}^*)^{(1)} \oplus \bigoplus_{2 \leq i \leq n-2} (-\omega_{i-1}+\omega_i - \omega_{n-1} + \omega_n)^{(1)}\\
& \quad \oplus (-\omega_{n-2}+\omega_n)^{(1)} \oplus (\omega_1-\omega_{n-1}+\omega_n)^{(1)}.
\end{align*}
\item[(f)] If $\Phi=D_4$, then 
\indent \begin{align*}
\operatorname{H}^2(B_1,k) &\cong (\mathfrak{u}^*)^{(1)} \oplus (\omega_1-\omega_2+\omega_3)^{(1)} \oplus (\omega_2 - \omega_4)^{(1)} \oplus (\omega_2 - \omega_3)^{(1)}\\
& \quad  \oplus (\omega_1-\omega_2 + \omega_4)^{(1)} \oplus (-\omega_2+\omega_3+\omega_4)^{(1)} \oplus (-\omega_1 + \omega_2)^{(1)}.
\end{align*}
\item[(g)]If $\Phi=D_n$, for $n \geq 5$, then
$$\operatorname{H}^2(B_1,k) \cong (\mathfrak{u}^*)^{(1)} \oplus \bigoplus_{i=1}^{n-3} (-\omega_i + \omega_{i+1})^{(1)} \oplus (\omega_1)^{(1)} \oplus (-\omega_{n-2}+\omega_{n-1}+\omega_n)^{(1)}.$$
\item[(h)] If $\Phi=F_4$, then
\indent \begin{align*}
\operatorname{H}^2(B_1,k) &\cong (\mathfrak{u}^*)^{(1)} \oplus (-\omega_1+\omega_2-\omega_4)^{(1)} \oplus (-\omega_1+\omega_2-\omega_3+\omega_4)^{(1)}\\
& \quad \oplus (\omega_1-\omega_3+\omega_4)^{(1)} \oplus (\omega_1-\omega_2+\omega_3)^{(1)} \oplus (\omega_2 - \omega_3)^{(1)} \oplus (\omega_1 - \omega_4)^{(1)}.
\end{align*}
\item[(i)] If $\Phi=G_2$, then $$\displaystyle{\operatorname{H}^2(B_1,k) \cong (\mathfrak{u}^*) \oplus (-\omega_1+\omega_2)^{(1)}}.$$
\end{itemize}
\end{thm}

\begin{proof}
Express the weight in $\h^2(U_1,k)$ (Theorem 1.4.1) in terms of the fundamental weight basis and determine which ones are $T_1$-invariant (i.e. multiples of 2).
\end{proof}

\subsubsection{}

Recall, it is enough to compute $\h^2(B_1,\la)$ for $\la \in X_1(T)$.  The following lemma gives the unique weight $\nu$ such that $\la = w \cdot 0 + 2\nu \in X_1(T)$.

\begin{lem}\label{lem:B1,T}
Let $p=2$.  For $w = s_{\a_i}s_{\a_j} \in W$ with $i < j$ and $\a_i$ and $\a_j$ not connected, we define
\[ \v_w = 
      \left\{ \begin{array}{l} 
   \displaystyle{\omega_i-\omega_k+\omega_j,\;\;\; \mbox{ if  } \a_i, \a_j \mbox{ separated by a single vertex, } \a_k,} \\ \;\;\;\;\;\displaystyle{\omega_i+\omega_j,\;\;\;\;\;\;\; \mbox{ otherwise.}} \\
 \end{array} \right. \]
 except in the following cases:\\
 \[ \v_w = 
      \left\{ \begin{array}{l} 
   \displaystyle{\omega_{n-3}-\omega_{n-2}+\omega_{n-1}-\omega_n,\;\;\; \mbox{ if  } \Phi=B_n, w=s_{\a_{n-3}}s_{\a_{n-1}}}, \\ \displaystyle{\omega_i+\omega_{n-1}-\omega_n,\;\;\;\;\;\;\;\;\;\;\;\;\;\;\;\;\;\;\; \mbox{ if } \Phi=B_n, w=s_{\a_i}s_{\a_{n-1}},\;i \not= n-3},\\
   \displaystyle{\omega_i - \omega_{n-1}+\omega_n, \;\;\;\;\;\;\;\;\;\;\;\;\;\;\;\;\;\;\; \mbox{ if } \Phi=C_n, w=s_{\a_i}s_{\a_n}, \;i\not=n-2},\\
   \displaystyle{-\omega_{n-2}+2\omega_{n-1}, \;\;\;\;\;\;\;\;\;\;\;\;\;\;\;\;\;\; \mbox{ if } \Phi=C_n, w=s_{\a_{n-1}}s_{\a_n}}.\\
 \end{array} \right. \]
 Then $s_{\a_i}s_{\a_j} \cdot 0 + 2\v_w \in X_1(T)$.
\end{lem}
 
\begin{proof}
For $\la$ a weight in $\h^2(U_1,k)$, then write $\la=s_{\a_i}s_{\a_j} \cdot 0 + 2\v$.  Note that $\a_i$ and $\a_j$ are not connected except when $\Phi$ is of type $C_n$ and $w=s_{\a_{n-1}}s_{\a_n}$.  Without loss of generality we can assume $i<j$.  The $\v_w$ is determined after writing $w \cdot 0$ in the fundamental weight basis.
\end{proof}

In the following theorem, we are only concerned with $\la$  a weight of $\h^2(U_1,k)$.

\begin{thm}\label{thm:B1,T}
Let $p=2$.  Then, as a $T$-module, 
 If $\lambda = s_{\a_i}s_{\a_j} \cdot 0 + 2 \v_w \in X_1(T)$ then
 
\[ \h^2(B_1,\la) \cong \left\{ \begin{array}{l} 
   \displaystyle{\quad \v_w^{(1)}, \quad \quad \quad \mbox{ if  } j \not = n-1, n-2}  \\ \displaystyle{\v_w^{(1)} \oplus \omega_{i+1}^{(1)},\quad \mbox{ if } j=n-2} \\
 \end{array} \right. \]

except in the following cases
\begin{itemize}
\item[(i)] If $\Phi=B_n$ and $j=n$, then
$$\operatorname{H}^2(B_1,\lambda) \cong \v_w^{(1)} \oplus (\omega_{n-1}-\omega_n)^{(1)}.$$
\item[(ii)] If $\Phi=B_n$ and $w=s_{\a_{n-3}}s_{\a_{n-1}}$, then
$$\operatorname{H}^2(B_1,\lambda) \cong \v_w^{(1)} \oplus (\omega_{n-2}-\omega_{n})^{(1)}.$$
\item[(iii)] If $\Phi=B_n$ and $w=s_{\a_i}s_{\a_{n-1}}$ for $i\not=n-3$, then
$$\operatorname{H}^2(B_1,\lambda) \cong \v_w^{(1)} \oplus(\omega_i + \omega_n)^{(1)}.$$
\item[(iv)] If $\Phi=C_n$ and $w=s_{\a_i}s_{\a_n}$ for $i \not= n-1$, then
$$\operatorname{H}^2(B_1,\lambda) \cong \v_w^{(1)} \oplus (\omega_i-\omega_{n-2}+\omega_{n-1})^{(1)}.$$
\item[(v)] If $\Phi=C_n$ and $w=s_{\a_{n-1}}s_{\a_n}$, then 
$$\operatorname{H}^2(B_1,\lambda) \cong \v_w^{(1)}.$$
\item[(vi)] If $\Phi=D_n$ and $w=s_{\a_{n-3}}s_{\a_{n-1}}$, then
$$\operatorname{H}^2(B_1,\lambda) \cong \v_w^{(1)} \oplus (\omega_{n-2}-\omega_n)^{(1)} \oplus (\omega_n)^{(1)}.$$
\item[(vii)] If $\Phi=D_n$ and $w=s_{\a_{n-3}}s_{\a_n}$, then
$$\operatorname{H}^2(B_1,\lambda) \cong \v_w^{(1)} \oplus (\omega_{n-2}-\omega_{n-1})^{(1)} \oplus (\omega_{n-1})^{(1)}.$$
\item[(viii)] If $\Phi=E_6$ and $w=s_{\a_2}s_{\a_3}$, then
$$\operatorname{H}^2(B_1,\lambda) \cong \v_w^{(1)} \oplus (\omega_4-\omega_5)^{(1)} \oplus (\omega_5-\omega_6)^{(1)} \oplus \omega_6^{(1)}.$$
\item[(ix)] If $\Phi=E_6,E_7,E_8$ and $w=s_{\a_2}s_{\a_5}$, then
$$\operatorname{H}^2(B_1,\lambda) \cong \v_w^{(1)} \oplus (\omega_4-\omega_3)^{(1)} \oplus (\omega_3-\omega_1)^{(1)} \oplus \omega_1^{(1)}.$$
\item[(x)] If $\Phi=E_6,E_7,E_8$ and $w=s_{\a_3}s_{\a_5}$, then
$$\operatorname{H}^2(B_1,\lambda) \cong \v_w^{(1)} \oplus (\omega_4-\omega_2)^{(1)} \oplus \omega_2^{(1)}.$$
\item[(xi)] If $\Phi=E_7$ and $w=s_{\a_2}s_{\a_3}$, then
$$\operatorname{H}^2(B_1,\lambda) \cong \v_w^{(1)} \oplus (\omega_4-\omega_5)^{(1)} \oplus (\omega_5-\omega_6)^{(1)} \oplus (\omega_6 - \omega_7)^{(1)} \oplus \omega_7^{(1)}.$$
\item[(xii)] If $\Phi=E_8$ and $w=s_{\a_2}s_{\a_3}$, then
\begin{align*}
\operatorname{H}^2(B_1,\lambda) &\cong \v_w^{(1)} \oplus (\omega_4-\omega_5)^{(1)} \oplus (\omega_5-\omega_6)^{(1)} \oplus (\omega_6 -\omega_7)^{(1)} \\
& \quad \oplus (\omega_7-\omega_8)^{(1)} \oplus \omega_8^{(1)}.
\end{align*}
\item[(xiii)] If $\Phi=F_4$ and $w=s_{\a_1}s_{\a_3}$, then
$$\operatorname{H}^2(B_1,\lambda) \cong \v_w^{(1)} \oplus (\omega_2-\omega_3)^{(1)} \oplus (\omega_3 - \omega_4)^{(1)}.$$
\end{itemize}
\end{thm}

\begin{proof}
The relationship between the $U_1-$ and the $B_1$-cohomology (as a $T$-module), as given in equation \eqref{U,B_cohom}, gives us,
$$\operatorname{H}^2(B_1,\lambda) \cong \operatorname{H}^2(U_1,k)_{-\overline{\lambda}} \cong \bigoplus_\nu \operatorname{H}^2(U_1,k)_{-\lambda + p\nu},$$
for $\la \in X_1(T)$.  From Lemma \ref{lem:B1,T}, the unique weight $\v_w$ such that $\lambda = s_{\a_i}s_{\a_j} + p\v_w \in X_1(T)$.  So we now have that 
$$\operatorname{H}^2(B_1, \lambda) \cong \bigoplus_\sigma \operatorname{H}^2(U_1,k)_{\la + 2 \sigma}.$$ 
We want to find $\sigma$ such that $\la=-(s_{\a_i}s_{\a_j}) + 2 \sigma$ is a weight of $\h^2(U_1,k)$.\\
\end{proof}

\subsubsection{$B$-module structure}

The $T$-module structure of $\h^2(B_1,\la)$ gives rise the $B$-module structure of $\h^2(B_1,\la)$.  Recall that $\operatorname{H}^2(B_1,\lambda)$ is a subquotient of  $S^2(\mathfrak{u}^*)_{-\lambda}$.  Since $B$ acts on $\mathfrak{u}$ by the adjoint action, then for $\a \in \Phi$, $\mbox{Dist}(B) = \left< {H_i \choose m}, \frac{X_{\alpha}^n}{n!}\right>$ acts on $\mathfrak{u}^*$, where $H_i = (d\phi_i)(1)$ and the $\phi_i$'s are a basis for $\mbox{Hom}(G_m,T) \cong \Z^r$ \cite[II.1.11]{Jan1}.  In particular since $\mathfrak{u}$ corresponds to the negative roots, then it's only necessary to look at $\frac{X_{\alpha}^n}{n!}$ for $\alpha \in \Phi^-$.  The action from \cite[26.3]{Hum} is defined by
\begin{equation}\label{modstructure}
\frac{X_\alpha^n}{n!}(u \t v) = \sum_k (\frac{X_\alpha^k}{k!}u \t \frac{X_\alpha^{n-k}}{(n-k)!}v).
\end{equation}
Using the results in the previous section.  If $\h^2(B_1,\la)$ has an answer consisting of only one factor, then this forms a 1-dimensional module.  However, it remains to be determined if the $m$-dimensional submodule of the cohomology is an indecomposable module.  As before, we will first look at this with the trivial module, then move on to arbitrary weights.

\begin{thm}\label{B1,B0}
Let $p=2$ then as a $B$-module, $\operatorname{H}^2(B_1,k) \cong (\mathfrak{u}^*)^{(1)}$.  Except in the following cases:\\
\input{B1k-mod}

\end{thm}
\begin{proof} Using \eqref{modstructure} and Theorem \ref{thm:B1,T0}, then we can determine which weights form an indecomposable module.  If the action takes one cohomology class to another, then the factors form an indecomposable module.
\end{proof}

\begin{thm}\label{thm:B1,B}
Let $p=2$ and $\la=w \cdot 0 + p\v_w$, then $\h^2(B_1,\la)$ is given by the following table for specified $w$.
\input{B1la-mod}

\end{thm}

\begin{proof}
Applying the action, \eqref{modstructure} to the various cohomology classes, which are explicitly listed in Appendix A.
\end{proof}  

\begin{rem} For the indecomposable $B$-modules in the preceding theorems, the factors are listed in order: where the first factor is the head and the last factor is the socle.
\end{rem}
\begin{cor}\label{cor:B1}
Let $p=2$ and $\lambda,\gamma \in X(T)$.
\begin{itemize}
\item[(a)]If $\lambda \not \in pX(T)$ and $\lambda \not= w \cdot 0 + p\sigma$ for some $w \in W$ with $l(w)=2$ and $\sigma \in X(T)$, then $\operatorname{H}^2(B_1,\lambda)=0.$
\item[(b)] If $\a \in \Delta$, then $\operatorname{H}^2(B_1,p\gamma - \alpha) = 0$.
\end{itemize}
\end{cor}

\subsection{$B_r$-Cohomology}
To calculate $\h^2(B_r,\la)$, we first make a few observations.

\subsubsection{} The following observations on $\la \in pX(T)$ will be used in $\h^2(B_r,\la)$ for $\la \in X(T)$.
\begin{lem}\label{lem:Br,0} Let $p=2$ and $\lambda \in X(T)$.  Then $\displaystyle{\operatorname{Hom}_{B_r/B_l}(k,\h^2(B_1,k)^{(l-1)}\t p^l\la)}$ is given by the following table, with $\g \in X(T)$
\input{HOMtable}
\end{lem}

\begin{proof}
 We have
\begin{align*} 
\operatorname{Hom}_{B_r/B_l}(k,\operatorname{H}^2(B_1,k)^{(l-1)} \t p^l\lambda) &\cong \operatorname{Hom}_{B_{r-l}}(k,\operatorname{H}^2(B_1,k)^{(-1)} \t \lambda)^{(l)} \\
& \cong \operatorname{Hom}_{B_{r-l}}(-\lambda,\operatorname{H}^2(B_1,k)^{(-1)})^{(l)}.
\end{align*}
It's necessary to consider the $B$-socle of $\operatorname{H}^2(B_1,k)$.  In general, this is the $B$-socle of $\mathfrak{u}^*$, which is $\sum_{\b \in \Delta} k_\b$ by \cite{Jan2}.  When the cohomology is not $\Uu^*$, this accounts for the additional weights.  (The extra cases come from the cases when the socle is not $\Uu^*$.
\end{proof}

The following proposition gives an easy formula to compute $\operatorname{H}^2(B_r,k)$ from $\operatorname{H}^2(B_1,k).$  Use the Lyndon-Hochschild-Serre (LHS) spectral sequence for $B_1 \underline{\triangleleft } B_r$:\begin{equation}\label{LHS}
\operatorname{E}_2^{i,j} = \operatorname{H}^i(B_r/B_1, \operatorname{H}^j(B_1,\lambda)) \Longrightarrow \operatorname{H}^{i+j}(B_r,\lambda).
\end{equation}

\begin{prop}\label{prop:Br,0}
Let $p=2$.  Then $\operatorname{H}^2(B_r,k) \cong \operatorname{H}^2(B_1,k)^{(r-1)}.$
\end{prop}

\begin{proof}  Use induction on $r$.  Assume that $r>1$.  Consider \eqref{LHS} for $j=1$.  $\operatorname{H}^1(B_1,k) = 0$ and thus $\operatorname{E}_2^{i,1} = 0$ by \cite{Jan2}.  Now applying the previous lemma with $l=1$ and $\lambda=0$, we have $\operatorname{E}_2^{0,2} = \operatorname{Hom}_{B_r/B_1}(k,\operatorname{H}^2(B_r/B_1,k)) =0.$  Since all differentials going into and out of $E_2^{2,0}$ are zero, and by the induction hypothesis we have that $$\operatorname{H}^2(B_r,k) \cong \operatorname{E}_2^{2,0} = \operatorname{H}^2(B_{r-1},k)^{(1)} \cong \operatorname{H}^2(B_1,k)^{(r-1)}.$$
\end{proof}

\subsubsection{}  

\begin{lem}\label{Br,p}
Let $0 \leq l < r$, and $\a \in \Delta$.\\
\begin{itemize}
\item[(a)]Then
\[ \operatorname{H}^2(B_r,-p^l\a) \cong 
      \left\{ \begin{array}{ll} 
   \displaystyle{k}& \mbox{ if  } l > 0,  \\
0 & \mbox{ otherwise.} \\
 \end{array} \right. \]
 
\item[(b)] Suppose $\Phi$ is of type $A_3$.  Then
\[ \operatorname{H}^2(B_r,-p^l\omega_2 ) \cong 
      \left\{ \begin{array}{ll} 
   \displaystyle{k}& \mbox{ if  } l > 0,  \\
0 & \mbox{ otherwise.} \\
 \end{array} \right. \]
 
\item[(c)] Suppose $\Phi$ is of type $B_3$.  Then
\[ \operatorname{H}^2(B_r,-p^l(\omega_2 - \omega_1)) \cong 
      \left\{ \begin{array}{ll} 
   \displaystyle{k}& \mbox{ if  } l > 0,  \\
0 & \mbox{ otherwise.} \\
 \end{array} \right. \]
 
\item[(d)] Suppose $\Phi$ is of type $B_4$ and $\lambda = -p^l(-\omega_2+\omega_3)$ or $\lambda=-p^l(\omega_1-\omega_2+\omega_3-\omega_4)$.  Then
\[ \operatorname{H}^2(B_r,\lambda) \cong 
      \left\{ \begin{array}{ll} 
   \displaystyle{k}& \mbox{ if  } l > 0,  \\
0 & \mbox{ otherwise.} \\
 \end{array} \right. \]
 
 \item[(e)] Suppose $\Phi$ is of type $B_n$.  Then
\[ \operatorname{H}^2(B_r,-p^l(-\omega_{n-2} + \omega_{n-1})) \cong 
      \left\{ \begin{array}{ll} 
   \displaystyle{k}& \mbox{ if  } l > 0,  \\
0 & \mbox{ otherwise.} \\
 \end{array} \right. \]
 
\item[(f)] Suppose $\Phi$ is of type $C_n$.  Then
\[ \operatorname{H}^2(B_r,-p^l(-\omega_{n-2} + \omega_n)) \cong 
      \left\{ \begin{array}{ll} 
   \displaystyle{k}& \mbox{ if  } l > 0,  \\
0 & \mbox{ otherwise.} \\
 \end{array} \right. \]
 
 \item[(g)] Suppose $\Phi$ is of type $D_4$ and $\lambda = -p^l(\omega_1-\omega_2+\omega_i)$ for $i \in \{3,4\}$\\ or $\lambda=-p^l(-\omega_2+\omega_3+\omega_4)$.  Then
\[ \operatorname{H}^2(B_r,\lambda) \cong 
      \left\{ \begin{array}{ll} 
   \displaystyle{k}& \mbox{ if  } l > 0,  \\
0 & \mbox{ otherwise.} \\
 \end{array} \right. \]
 
\item[(h)] Suppose $\Phi$ is of type $D_n$.  Then
\[ \operatorname{H}^2(B_r,-p^l(-\omega_{n-2} + \omega_{n-1} + \omega_n)) \cong 
      \left\{ \begin{array}{ll} 
   \displaystyle{k}& \mbox{ if  } l > 0,  \\
0 & \mbox{ otherwise.} \\
 \end{array} \right. \]
 
 \item[(i)] Suppose $\Phi$ is of type $F_4$.  Then
\[ \operatorname{H}^2(B_r,-p^l(-\omega_1+\omega_2 - \omega_4)) \cong 
      \left\{ \begin{array}{ll} 
   \displaystyle{k}& \mbox{ if  } l > 0,  \\
0 & \mbox{ otherwise.} \\
 \end{array} \right. \]
 
\item[(j)] Suppose $\Phi$ is of type $G_2$.  Then
\[ \operatorname{H}^2(B_r,-p^l(\omega_2 - \omega_1)) \cong 
      \left\{ \begin{array}{ll} 
   \displaystyle{k}&\mbox{ if  } l > 0,  \\
0 & \mbox{ otherwise.} \\
 \end{array} \right. \]
 \end{itemize}
 \end{lem}
 
\begin{proof}
(a) When $r=1$, $\operatorname{H}^2(B_1,-\a)=0$ by Corollary \ref{cor:B1}.  Suppose $r>1$ and $l=0$.  By Equation \eqref{LHS}, then $\displaystyle{E_2^{i,0}= \operatorname{H}^i(B_r/B_1,\operatorname{Hom}_B(k,-\a)) = 0}$, and $E_2^{i,2}= \operatorname{H}^i(B_r/B_1,\operatorname{H}^2(B_1,-\a)) = 0$, by Theorem \ref{thm:B1,B}.  Furthermore, we have $E_2^{1,1} = \operatorname{H}^1(B_r/B_1,\operatorname{H}^1(B_1,-\a))$.   Generally $\operatorname{H}^1(B_1,-\a) = k$ except when $\Phi$ is of type $A_3,D_4,D_n$.\\
(i) If $\Phi = A_3$, then $$\operatorname{H}^1(B_1,-\a_1) = k \oplus (-\omega_1+\omega_3)^{(1)}$$ and $$\operatorname{H}^1(B_1,-\a_3) = k \oplus (\omega_1-\omega_3)^{(1)}$$\\
(ii) If $\Phi=D_4$, then $$\operatorname{H}^1(B_1,-\a_1) = k \oplus (-\omega_1+\omega_3)^{(1)} \oplus (-\omega_1 + \omega_4)^{(1)},$$ $$\operatorname{H}^1(B_1,-\a_3) = k \oplus (\omega_1-\omega_3)^{(1)} \oplus (-\omega_3 + \omega_4)^{(1)},$$ and $$\operatorname{H}^1(B_1,-\a_4) = k \oplus (\omega_1-\omega_4)^{(1)} \oplus (\omega_3 - \omega_4)^{(1)}.$$\\
(iii) If $\Phi=D_n$, then $$\operatorname{H}^1(B_1,-\a_{n-1}) = k \oplus (-\omega_{n-1}+\omega_n)^{(1)}$$ and $$\operatorname{H}^1(B_1,-\a_n) = k \oplus (\omega_{n-1}-\omega_n)^{(1)}.$$  So, now we can apply \cite[Thm 2.8(C)]{BNP1} and we have that $E_2^{1,1}$ vanishes.  Hence\\ $E_2^{2,0}=E_2^{1,1}=E_2^{0,2}=0$ in all cases and so $\operatorname{H}^2(B_r,-\a)=0$.\\
Now assume $l>0$.  We use \eqref{LHS} again. 
$E_2^{i,1} = \operatorname{H}^i(B_r/B_l,\operatorname{H}^1(B_l,k) \t -p^l\a)$.  By \cite[Thm 2.8(C)]{BNP1}.  From the case $l=0$, we have
$$E_2^{2,0} = \operatorname{H}^2(B_r/B_l,-p^l\a) \cong \operatorname{H}^2(B_{r-l}, -\a)^{(l)}=0.$$
So, $\h^2(B_r,-p^l\a)=E_2^{0,2}$.  From Lemma \ref{lem:Br,0}, 
$$E_2^{0,2} = \operatorname{Hom}_{B_r/B_l}(k,\operatorname{H}^2(B_l,k)\t\-p^l\a) \cong \operatorname{Hom}_{B_r/B_l}(k,\operatorname{H}^2(B_1,k)^{(l-1)} \t -p^l\a) \cong k.$$
Hence, the result follows. 
 
\indent For (b)-(j) and $\la$ defined in the theorem.  As before, the case $r=1$ follows from Corollary \ref{cor:B1}.  For $r>1$ and $l=0$, we use \eqref{LHS}.
Note for all cases $\operatorname{H}^1(B_1,-\omega_i)=0$ and so $E_2^{i,1}=0$.
\end{proof}
 
\subsubsection{} \label{sec:Br}
Now consider the case when $\la \not\in p^rX(T)$.  Define the following indecomposable $B$-modules, where all of the factors are listed from the head to the socle:\\
\\
$\indent \bullet$ $N_{B_n}$ is the two-dimensional indecomposable $B$-module with factors $\a_3$ and $k$.\\ \indent  Furthermore, $N_{B_n} \cong \h^1(B_1,-\a_{n-1})^{(-1)}$ \cite[3.6]{Jan2}.\\
$\indent \bullet$ $N_{C_n}$ is the $n$-dimensional indecomposable $B$-module with factors\\ $\indent \omega_1, \omega_2-\omega_1, \omega_3-\omega_2, \ldots, \omega_n - \omega_{n-1}$.\\ \indent Furthermore, $N_{C_n} \cong \h^1(B_1,k)^{(-1)}$. \cite[3.6]{Jan2}\\
$\indent \bullet$ $N_{F_4}$ is the three-dimensional indecomposable $B$-module with factors $\a_3+\a_4$, $\a_3$, $k$.\\ \indent  Furthermore $N_{F_4} \cong \h^1(B_1,-\a_2)^{(-1)}$ \cite[3.6]{Jan2}.\\
$\indent \bullet$ $N_{G_2}$ is the two-dimensional indecomposable $B$-module with factors $\a_2$ and $k$.\\ \indent  Furthermore $N_{G_2} \cong \h^1(B_1,-\a_2)^{(-1)}$ \cite[3.7]{Jan2}.\\
\\
Notice that if $N$ is one of the above modules, then $N \t \la$ remains indecomposable for any weight $\la$.
\begin{lem}\label{lem:U1,M}
\begin{itemize}
\item[(a)] If $\Phi$ is of type $B_n,F_4,G_2$, then $\h^1(U_1,N) \cong \h^1(\Uu,N)$, where $N$ is one of the modules defined above.
\item[(b)] If $\Phi$ is of type $B_n,F_4,G_2$, then $\h^1(U_1,N)$ is defined as follows:
\begin{itemize}
\item[(i)] $\h^1(U_1,N_{B_n})$ has a $T$-basis $\{\a_1,\a_2,\ldots,\a_{n-1}, \a_{n-1}+\a_n, 2\a_n, \a_{n-1}+2\a_n\}$.
\item[(ii)] $\h^1(U_1,N_{C_n})$ has a $T$-basis\\  $\{\a_1,\a_2,\ldots,\a_{n-2},\a_n,\a_{n-1}+\a_n,2\a_{n-1},2\a_{n-1}+\a_n,\a_{n-2}+2\a_{n-1}\}$
\item[(iii)] $\h^1(U_1,N_{F_4})$ has a $T$-basis $\{\a_1,\a_2,\a_4,\a_2+\a_3,2\a_3,\a_2+2\a_3,2\a_3+\a_4\}$.
\item[(iv)] $\h^1(U_1,N_{G_2})$ has a $T$-basis $\{\a_1, \a_1+\a_2, 2\a_2, 2\a_1+\a_2\}$.\\
\end{itemize}
\end{itemize}
\end{lem}

\begin{proof}
(a) Consider the following exact sequence:
$$0 \rightarrow \h^1(U_1,N) \hookrightarrow \h^1(\Uu,N) \rightarrow \operatorname{Hom}^s(\Uu,N^\Uu),$$
where $\operatorname{Hom}^s(\mathfrak{u},N^\mathfrak{u})$ is the set of all maps that are additive and satisfy the property: $\phi(ax)=a^p \phi(x)$.\\  
\indent In our case we have that $\operatorname{Hom}^s(\Uu,N^\Uu) \cong (\Uu^*)^{(1)}$  So, we now have the following exact sequence:
$0 \rightarrow \h^1(U_1,N) \hookrightarrow \Lambda^1(\Uu^*)\t N \rightarrow (\Uu^*)^{(1)}$
Since $\Uu^*$ is spanned by the negative roots, then this last map is the $0$-map.  Therefore, $\h^1(U_1,N) \cong \h^1(\Uu,N)$.\\
(b) To calculate $\h^1(\Uu,N)$, we use \cite[I.9.15]{Jan1}:
$k \t N \rightarrow \Uu^* \t N \rightarrow \Lambda^2(\Uu^*)\t N.$\\
So, it is first necessary to calculate $\ker(\Uu^* \t N \rightarrow \Lambda^2(\Uu^*) \t N)$ by \cite[I.9.17]{Jan2} $$d_i'(m \t \psi)= \sum_j m_j\t(\phi_j \wedge \psi) + m \t d_i(\psi),$$
where $d_i$ is the differential defined in Section 2.1 and $\psi \in \Phi$.
  For example, when $\Phi$ is of type $B_3$, then $N$ is an indecomposable $B$-module with factors $\a_3=m_1$ and $k=m_0$.   Then $x_\gamma \cdot m_0 = 0$ and 
  \[ x_\gamma \cdot m_1 = \left\{ \begin{array}{ll}
  m_0 & \gamma=-\a_3,\\
  0 & \mbox{else.}\\
  \end{array} \right. \]
  Therefore, $d_0'(m_0)~=~0$ and $d_0'(m_1)~=~x_{-\a_3}^* ~\t~ m_0$.  Furthermore, we get that
  $d_1'(m_0 ~\t~ \psi)~=~m_0 \t d_1(\psi)$ and $d_1'(m_1 \t \psi)~ =~ m_0 \t (x_{-\a_3} \wedge \psi)~ +~ m_1 \t d_1(\psi)$.  It is necessary to determine when any linear combination of these maps (both with the same map) returns 0.   $$\{\a_1,\a_2,\a_3,\a_2+\a_3,\a_2+2\a_3,2\a_3\},$$ 
and check which weights are in the previous image $\displaystyle{k \t N \rightarrow \Uu^* \t N}$. The weight $\a_3$ is in the previous image because $\a_3$ is one of the factors of our module and $1 \t \a_3 = \a_3 \t 1$.  The other calculations are similar. 
\end{proof}
\begin{thm}\label{thm:Br,M}
\begin{itemize}
\item[(a)] If $\Phi$ is of type $B_n$, then
 \[ \operatorname{H}^2(B_r,N_{B_n}\t\la)^{(-r)}  \cong 
      \left\{ \begin{array}{ll}
   \displaystyle{\v} & \quad \mbox{ if  } \la=p^r\v-p^l\a, \a\in \Delta, 0 \leq l < r \mbox{ and } \a \in \Delta\\
& \qquad \mbox{for } l\not=r-1 \mbox{ if } \a=\a_{n-1}, \mbox{ and } l \not= 0 \mbox{ if } \a=\a_n,  \\
   \displaystyle{\v} & \quad \mbox { if } \la=p^\g - (\a_{n-1}+\a_n), \\
   \displaystyle{M\t\v} & \quad \mbox { if } \la=p^r\v - p^{r-1}\a_{n-1}, \\
   \displaystyle{M\t\v} & \quad \mbox { if } \la=p^r\v, \\
0 & \quad \mbox{ otherwise.} \\
 \end{array} \right. \]
where $M$ is an indecomposable module with factors $\a_n$ and $k$.
\item[(b)] If $\Phi$ is of type $C_n$
\[ \h^2(B_r,N_{C_n}\t\la)^{(-r)}  \cong
	\left\{ \begin{array}{lll}
	\v && \mbox{if } \la=p^r\v-p^l\a, \mbox{ with } \a \in \Delta, \; 0 \leq l < r \mbox{ and } \\
	&& \quad \mbox{where } l \not= r-1 \mbox{ if } \a=\a_n, \a_{n-2} \mbox{ and } l \not=0 \mbox{ if }\\
&& \quad \a=\a_{n-1},\\
	 \v && \mbox{if } \la=p^r\v-(\a_{n-1}+\a_n), \\
	 M\t\v && \mbox{if } \la=p^r\v-p^{r-1}\a \mbox{ where } \a\in \{\a_{n-2},\a_n\},\\
	 M\t\v && \mbox{if } \la=p^r\v, \\
	 0 && \mbox{else,} \\
	 \end{array} \right. \]
	 where $M$ is an indecomposable module with factors $\a_{n-1}$ and $k$.
\item[(c)] If $\Phi$ is of type $F_4$ 
\[ \operatorname{H}^2(B_r,N_{F_4}\t\la)^{(-r)}  \cong 
      \left\{ \begin{array}{ll} 
   \displaystyle{\v} & \quad \mbox{ if  } \la=p^r\v-p^l\a, \mbox{ with }  \leq l < r \mbox{ and } \a \in \Delta\\
& \qquad l\not=r-1 \mbox{ when } \a\in \{\a_2,\a_4\},  \\
   \displaystyle{\v} & \quad \mbox { if } \la=p^\g - (\a_2+\a_3), \\
   \displaystyle{\v} & \quad \mbox { if } \la=p^\g - (\a_3+\a_4), \\
   \displaystyle{M\t\v} & \quad \mbox { if } \la=p^r\v - p^{r-1}\a \mbox{ for } \a \in \{a_2,\a_4\}, \\
   \displaystyle{\v} & \quad \mbox {if } \la=p^r\v, \\
0 & \quad \mbox{ otherwise,} \\
	\end{array} \right. \]
where $M$ is the two-dimensional indecomposable module with factors $\a_3$ and $k$.
\item[(d)] If $\Phi$ is of type $G_2$
\[ \operatorname{H}^2(B_r,N_{G_2}\t\la)^{(-r)} \cong 
      \left\{ \begin{array}{ll} 
   \displaystyle{\v} & \quad \mbox{ if  } \la=p^r\v-p^l\a, \a\in \Delta 0 \leq l < r-1,\\
& \qquad l\not=0 \mbox{ when } \a=\a_2,  \\
   \displaystyle{\v} & \quad \mbox { if } \la=p^\g - (\a_1+\a_2), \\
   \displaystyle{M\t\v} & \quad \mbox { if } \la=p^r\v - p^{r-1}\a_2, \\
   \displaystyle{\v} & \quad \mbox {if } \la=p^r\v, \\
0 & \quad \mbox{ otherwise,} \\
	\end{array} \right. \]
	where $M_1$ is the two-dimensional indecomposable module with factors $\a_1$ and $k$.
\end{itemize}
\end{thm}
 
\subsubsection{Statement of Theorem}
The previous calculations helps us compute $\operatorname{H}^2(B_r,\lambda)$ for any $r$ and $\lambda \in X_r(T)$.

\begin{thm}\label{thm:Br}
Let $p=2$ and $\la \in X_r(T)$.  Then $\operatorname{H}^{2}(B_{r},\lambda)$ is given 
by the following table.

\input{B_rtable}

\end{thm}
 
 \begin{rem}
If $\la$ is not listed above, then $\h^2(B_r,\la)=0$.  $M$ is the module with factors listed in Theorem \ref{thm:Br,M}.  \end{rem}
 
 \begin{proof}
 We will use induction on $r$.  For $r=1$, the claim reduces to Theorem \ref{thm:B1,B}.  Suppose $r>1$.  Set $\lambda=\lambda_0+p\lambda_1$ where $\lambda_0 \in X_1(T)$ and $\lambda_1 \in X(T)$.  Use \eqref{LHS}.\\
{\bf Case 1:} $\lambda_0 \not\equiv 0$ and $\lambda_0 \not\equiv -\a \mod pX(T)$, with $\a \in \Delta$.\\
In this case we have $E_2^{i,0}=0$ and $E_2^{i,1}=0$ \cite[3.2]{Jan2}.  Then,
$$\operatorname{H}^2(B_r,\lambda) = E^2 \cong E_2^{0,2} = \operatorname{Hom}_{B_r/B_1} (k,\operatorname{H}^2(B_1,\lambda_0) \t p\lambda_1.$$
By Lemma \ref{lem:Br,0} this expression is zero unless $\lambda_0 = w \cdot 0 + p\v_w$ for some $w\in W$ with $l(w)=2$ and $\v_w$ as given in Lemma \ref{lem:B1,T}.  Assume $\lambda_0$ is of this form.  Then by Theorem \ref{thm:B1,B}, the $B$-module $\operatorname{H}^2(B_1,\lambda_0)$ has simple socle of weight $p\v_w$ and $E_2^{0,2}$ vanishes unless $p(\v_w + \lambda_1) \in p^rX(T)$.  This implies $\lambda_0 = w\cdot0 + p^r\v$ with $l(w)=2$ and $\v \in X(T)$.  Moreover, $\operatorname{H}^2(B_r,\lambda) \cong \v^{(r)}$ for such weights.  To summarize:
if $\la_0 \not \equiv 0$ and $\la_0 \not \equiv -\a \mod pX(T)$, with $\a \in \Delta$, and $r > 1$ then
\[ \h^2(B_r,\la) \cong \left\{
\begin{array}{lll}
\v^{(r)} && \mbox{ if } \la=p^r\v + w \cdot 0, \mbox{ with } l(w)=2, \\
0 && \mbox{ else.}
\end{array} \right. \]
{\bf Case 2:} Suppose $\lambda_0 \equiv -\a  \mod pX(T)$, with $\a \in \Delta$.\\
Then $\lambda_0=p\omega_\a - \a$ by \cite[3.3]{Jan2}, except in the following cases:\\
\indent $\bullet$ $\Phi$ is of type $B_n$ with $\a=\a_{n-1}$, then $\la_0 = 2(\omega_{n-1}-\omega_n)-\a_{n-1}$.\\
\indent $\bullet$ $\Phi$ is of type $C_n$ with $\a=\a_n$, then $\la_0=2(-2\omega_{n-1}+\omega_n)-\a_n$.\\
\indent $\bullet$ $\Phi$ is of type $F_4$ with $\a=\a_2$, then $\la_0=2(\omega_2-\omega_3)-\a_2$.\\
\indent $\bullet$ $\Phi$ is of type $G_2$ with $\a=\a_2$, then $\la_0=2(-\omega_1+\omega_2)-\a_2$.\\
\noindent  If $\Phi \not= C_n$ with $\la \equiv -\a_n \mod pX(T)$, then $E_2^{i,0}=0$.  By Corollary \ref{cor:B1}, it follows that $E_2^{i,2}=0$.  Therefore 
$$\h^2(B_r,\la) \cong E_2^{1,1} = \operatorname{H}^1(B_r/B_1,\operatorname{H}^1(B_1,\lambda_0)\t p\lambda_1).$$
If $\Phi$ is of type $C_n$.  Then by \cite[3.5]{Jan2} we have 
$$E_2^{1,1} \cong \operatorname{H}^1(B_r/B_1,p(\omega_\a + \lambda_1)) \cong \operatorname{H}^1(B_{r-1},\omega_\a+\lambda_1)^{(1)}.$$  
\cite[2.8]{BNP2} implies that $E_2^{1,1}=0$ unless $\omega_\a + \lambda_1 = p^{r-1}\v - p^{k-1}\b$ for some $\b \in \Delta$ and some $0 < k < r$.  Moreover, in this case $\operatorname{H}^2(B_r,\lambda) \cong \v^{(r)},$ except in the following cases:\\
\indent $\bullet$ if $\Phi=B_n, n\not=4, \b=\a_{n-1}$ and $k=r-1$; in which case $\h^2(B_r,\la) \cong M_{B_n}^{(r)}\t\v^{(r)}$\\
\indent $\bullet$ if $\Phi=B_4, \b=\a_i, i \in {1,3}$ and $k=r-1$; in which case $\h^2(B_r,\la) \cong M_{B_4}^{(r)}\t\v^{(r)}$\\
\indent $\bullet$ if $\Phi=C_n, \b=\a_n$ and $k=r-1$; in which case $\h^2(B_r,\la) \cong M_{C_n}\t\v^{(r)}$\\
\indent $\bullet$ if $\Phi=F_4, \b=\a_2$ and $k=r-1$; in which case $\h^2(B_r,\la) \cong M_{F_4}\t\v^{(r)}$\\
\indent $\bullet$ if $\Phi=G_2, \b=\a_2$ and $k=r-1$; in which case $\h^2(B_r,\la) \cong M_{G_2}\t\v^{(r)}$\\

\noindent If $\Phi$ is of type $A_3$, then $p\omega_i - \a_i = \omega_2$ for $i \in \{1,3\}$.  Applying \cite[3.5(b)]{Jan2} yields $E_2^{1,1} \cong \operatorname{H}^1(B_r/B_1,p(\omega_1 \oplus \omega_3) \t p\lambda_1) \cong \operatorname{H}^1(B_{r-1},\omega_1+\lambda_1)^{(1)} \oplus \operatorname{H}^1(B_{r-1},\omega_3+\lambda_1)^{(1)}$.  As before, by \cite[2.8]{BNP2} the cohomology vanishes unless $\lambda=p^r\v - p^k\a_i - \a$  where $i \in \{1,3\}$.  Moreover, $\operatorname{H}^2(B_r,\lambda) \cong \gamma^{(r)}$, unless $k=r-1$, in which case 
 $$\operatorname{H}^2(B_r,\lambda)^{(-r)} \cong \gamma \oplus (\gamma + (-1)^j(\omega_1 - \omega_3)).$$  Adding $p^r\gamma_w$ as defined in Lemma \ref{lem:B1,T} to $\lambda$ results in the more symmetric statement $\operatorname{H}^2(B_r,p^r\gamma + p^{r-1}\omega_2 - \a) \cong (\g + \omega_1)^{(r)} \oplus (\g+\omega_2)^{(r)}$.\\
 \indent Similarly if $\Phi=D_4$ then $p\w_i - \a_i = \w_2$ for $i \in \{1,3,4\}$, $\operatorname{H}^2(B_r,p^r\gamma + p^{r-1}\omega_2 - \a) \cong (\g + \omega_1)^{(r)} \oplus (\g+\omega_3)^{(r)} \oplus (\g+\omega_4)^{(r)}$; if $\Phi=D_n$ for $n \geq 5$ then $p\w_i - \a_i = \w_{n-2}$ for $i \in \{n-1, n\}$,  $\operatorname{H}^2(B_r,p^r\gamma + p^{r-1}\omega_{n-2} - \a) \cong (\g + \omega_{n-1})^{(r)} \oplus (\g+\omega_n)^{(r)}$.\\
\indent If $\Phi$ is of type $B_n$ with $\la_0 \equiv \a_{n-1} \pmod {pX(T)}$, or $\Phi$ is of type $F_4$ with $\la_0 \equiv \a_2 \pmod {pX(T)}$, or $\Phi$ is of type $G_2$ with $\la_0 \equiv \a_2 \pmod{pX(T)}$.  Then, define $\v \in X(T)$ via $\la=p\v-\a_{n-1}$.  Then from \cite[3.6]{Jan2} $\h^1(B_1,\la) \cong M^{(1)}\t\v^{(1)}$ and so $$\h^2(B_r,\la) \cong \h^1(B_{r-1},M\t\v)^{(1)}.$$ Now apply Theorem \ref{thm:Br,M}.\\  
\indent Similarly, if $\Phi$ is of type $B_3$, then $\h^2(B_r,\la) \cong (\gamma+\omega_1)^{(r)}\oplus(\gamma+\omega_3)^{(r)}$ when $\la=p^r\omega_1-p^{r-1}\a_1-\a=p^r\omega_3 - p^{r-1}\a_3-\a$ and  $\Phi$ is of type $B_4$, 
$\h^2(B_r,\la) \cong (\omega_1+\g)^{(r)} \oplus (M_{B_4}\t\g)^{(r)}$ when $\la=p^r\omega_1-p^{r-1}\a_1-\a=p^r(\omega_3-\omega_4)-p^{r-1}\a_3-\a$\\  
\indent If $\Phi=C_n$ and $\a=\a_n$, then $\la_0 \equiv 0 \pmod {pX(T)}$, which is excluded.

 {\bf Case 3:} Now assume $\lambda_0=0$.  First assume $\Phi$ is not of type $C_n$. Then $E_2^{i,1}=0$ for all $i$ by \cite[3.3]{Jan2}.  From Lemma \ref{lem:Br,0} one obtains that $E_2^{0,2}=0$ unless $\lambda=p^r\gamma - p\a$, with $\a \in \Delta$ or if $\la$ is one of those listed in (i)-(iv) or (vi)-(viii), then by Lemma \ref{lem:Br,0}(B), then $\operatorname{H}^2(B_r,\lambda) \cong \v^{(r)}$, as claimed. 
 Now if $E_2^{0,2}=0$.
  This implies that
 $$E^2=E_2^{2,0} \cong \operatorname{H}^2(B_r/B_1,p\lambda_1) \cong \operatorname{H}^2(B_{r-1},\lambda_1)^{(1)}.$$
 If $\Phi$ is of type $C_n$, then by Lemma \ref{lem:Br,0} 
$$E_2^{0,2}=\operatorname{Hom}_{B_r/B_1}(k,\h^2(B_1,k)^{(1)}\t p\la_1) \cong \operatorname{Hom}_{B_{r-1}}(k,\h^2(B_1,k) \t \la_1)^{(1)} = 0.$$  Now, consider, $E_2^{2,0}=\h^2(B_r/B_1,\operatorname{Hom}_{B_1}(k,p\la_1)) \cong \h^2(B_{r-1},\la_1)^{(1)}$ and\\ $\displaystyle{E_2^{0,1}=\operatorname{Hom}_{B_r/B_1}(k,\h^1(B_1,k)\t p\la_1) \cong \operatorname{Hom}_{B_{r-1}}(k,M \t \la_1)^{(1)}}$ and the map\\ $d_2:E_2^{0,1} \rightarrow E_2^{2,0}$.  We want to show that $d_2$ is the 0-map.  But, $E_2^{0,1} \not= 0$ if and only if $\la_1=p^{r-1}\v-(\omega_n-\omega_{n-1})$, and $E_2^{2,0} \cong \h^2(B_{r-1},p^{r-1}\v-(\omega_n - \omega_{n-1}))=0$, by our previous calculations with case (i).\\
\indent Now, consider $E_2^{1,1} = \h^1(B_r/B_1,\h^1(B_1,k)\t p\la_1) \cong \h^1(B_{r-1},M \t \la_1)$, where $M$ is defined as in Section \ref{sec:Br} (as in \cite{Jan2}).  Since $E_2^{1,1} \not= 0$, then we must consider $E_2^{3,0}$.  We want to show $d_2:E_2^{1,1} \rightarrow E_2^{3,0}$ is the zero-map.  Note that the factors of $M$ are not in the root lattice, and so for $E_2^{1,1} \not=0$, then $\la_1 \not\in \Z\Phi$.  For $\la_1 \not\in \Z\Phi$, then $E_2^{3,0}=\h^3(B_{r-1},\la_1)=0$ by \cite[II.4.10]{Jan1}.  Therefore, $d_2$ is the zero-map.
 \end{proof}

%% file: B1k-mod.tex
\nopagebreak
\begin{longtable}{|l|l|l|}
\hline 
Root System & \large{\bf{$\h^2(B_1, k), \; p=2$}} & Module factors (listed from top to bottom) \\
\hline\hline 
$A_3$ & $(\Uu^*)^{(1)} \oplus M$ & $(\w_1-\w_2+\w_3)^{(1)}, \w_2^{(1)}$ \\
\hline
$B_3$ & $(\Uu^*)^{(1)} \oplus M$ & $\w_1^{(1)}, (-\w_1+\w_2)^{(1)}$ \\
\hline
$B_4$ & $(\Uu^*)^{(1)} \oplus M_1 \oplus M_2$ & $M_1: (\w_2-\w_4)^{(1)}, (\w_1-\w_2+\w_3-\w_4)^{(1)}$\\
& & $M_2: \w_1^{(1)}, (-\w_1+\w_2)^{(1)}, (-\w_2+\w_3)^{(1)}$ \\
\hline
$B_n, \; n \geq 5$ & $(\Uu^*)^{(1)} \oplus M$ & $\w_1^{(1)}, (-\w_1+\w_2)^{(1)}, (-\w_2+\w_3)^{(1)}, \ldots,$\\
& & $\quad (-\w_{n-2}+\w_{n-1})^{(1)}$ \\
\hline
$C_n$ & $(\Uu^*)^{(1)} \oplus M$ & $(\w_1-\w_{n-1}+\w_n)^{(1)}, (-\w_1+\w_2-\w_{n-1}+\w_n)^{(1)},$\\
& & $\quad (-\w_2+\w_3-\w_{n-1}+\w_n)^{(1)}, \ldots,$\\
& & $\quad (-\w_{n-3}+\w_{n-2}-\w_{n-1}+\w_n)^{(1)}, (-\w_{n-2}+\w_n)^{(1)}$ \\
\hline
$D_4$ & $(\Uu^*)^{(1)} \oplus M_1 \oplus M_2 \oplus M_3$ & $M_1:(\w_2-\w_4)^{(1)}, (\w_1-\w_2+\w_3)^{(1)}$\\
& & $M_2: (\w_2-\w_3)^{(1)}, (\w_1-\w_2+\w_4)^{(1)}$\\
& & $M_3: \w_1^{(1)}, (-\w_1+\w_2)^{(1)}, (-\w_2+\w_3+\w_4)^{(1)}$ \\
\hline
$D_n, \; n \geq 5$ & $(\Uu^*)^{(1)} \oplus M$ & $\w_1^{(1)}, (-\w_1+\w_2)^{(1)}, (-\w_2+\w_3)^{(1)}, \ldots, $\\
& & $\quad (-\w_{n-3}+\w_{n-2})^{(1)}, (-\w_{n-2}+\w_{n-1}+\w_n)^{(1)}$ \\
\hline
$F_4$ & $(\Uu^*)^{(1)} \oplus M$ & $(\w_2-\w_3)^{(1)}, (\w_1-\w_2+\w_3)^{(1)}, (\w_1-\w_3+\w_4)^{(1)},$\\
& & $\quad  (-\w_1+\w_2-\w_3+\w_4)^{(1)}, (\w_1-\w_4)^{(1)},$\\
& & $\quad  (-\w_1+\w_2-\w_4)^{(1)}$ \\
\hline
$G_2$ & $(\Uu^*)^{(1)} \oplus (-\w_1+\w_2)^{(1)}$ & \\
\hline
\end{longtable}

%% file: B1la-mod.tex
\nopagebreak
\begin{longtable}{|l|l|l|}
\hline 
Root System & $w$ in $\la=w \cdot 0 + p\v$ & Module factors (listed from top to bottom) \\
\hline
\hline
All & $w=s_{\a_i}s_{\a_{i+2}}$ & $\w_{i+1}^{(1)}, (\w_i-\w_{i+1}+\w_{i+2})^{(1)}$\\
\hline
$B_n$ & $w=s_{\a_{n-2}}s_{\a_n}$ & $(\w_{n-1}-\w_n)^{(1)}, (\w_{n-2}-\w_{n-1}+\w_n)^{(1)}$\\
\hline
$B_n$ & $w=s_{\a_{n-3}}s_{\a_{n-1}}$ & $(\w_{n-2}-\w_n)^{(1)}, (\w_{n-3}-\w_{n-2}+\w_{n-1}-\w_n)^{(1)}$\\
\hline
$B_n$ & $w=s_{\a_i}s_{\a_{n-1}}, \; i \not= n-3$ & $(\w_i+\w_n)^{(1)}, (\w_i+\w_{n-1}-\w_n)^{(1)}$\\
\hline
$C_n$ & $w=s_{\a_i}s_{\a_n}, \; i \not= n-2$ & $(\w_i-\w_{n-2}+\w_{n-1})^{(1)}, (\w_i-\w_{n-1}+\w_n)^{(1)}$\\
\hline
$D_n$ & $w=s_{\a_{n-3}}s_{\a_{n-1}}$ & $\w_n^{(1)}, (\w_{n-2}-\w_n)^{(1)}, (\w_{n-3}-\w_{n-2}+\w_{n-1})^{(1)}$\\
\hline
$D_n$ & $w=s_{\a_{n-3}}s_{\a_n}$ & $\w_{n-1}^{(1)}, (\w_{n-2}-\w_{n-1})^{(1)}, (\w_{n-3}-\w_{n-2}+\w_n)^{(1)}$\\
\hline
$E_6$ & $w=s_{\a_2}s_{\a_3}$ & $\w_6^{(1)}, (\w_5-\w_6)^{(1)}, (\w_4-\w_5)^{(1)}, (\w_2+\w_3-\w_4)^{(1)}$\\
\hline
$E_6, E_7, E_8$ & $w=s_{\a_2}s_{\a_5}$ & $\w_1^{(1)}, (\w_3-\w_1)^{(1)}, (\w_4-\w_3)^{(1)}, (\w_2-\w_4+\w_5)^{(1)}$\\
\hline
$E_6, E_7, E_8$ & $w=s_{\a_3}s_{\a_5}$ & $\w_2^{(1)}, (\w_4-\w_2)^{(1)}, (\w_3-\w_4+\w_5)^{(1)}$\\
\hline
$E_7$ & $w=s_{\a_2}s_{\a_3}$ & $\w_7^{(1)}, (\w_6-\w_7)^{(1)}, (\w_5-\w_6)^{(1)}, (\w_4-\w_5)^{(1)},$\\
& & $\quad (\w_2+\w_3-\w_4)^{(1)}$\\
\hline
$E_8$ & $w=s_{\a_2}s_{\a_3}$ & $\w_8^{(1)}, (\w_7-\w_8)^{(1)}, (\w_6-\w_7)^{(1)}, (\w_5-\w_6)^{(1)},$\\
& & $\quad (\w_4-\w_5)^{(1)}, (\w_2+\w_3-\w_4)^{(1)}$\\
\hline
$F_4$ & $w-s_{\a_1}s_{\a_3}$ & $(\w_3-\w_4)^{(1)}, (\w_2-\w_3)^{(1)}, (\w_1-\w_2+\w_3)^{(1)}$ \\
\hline

\end{longtable}

%% file: HOMtable.tex
\begin{longtable}{|l|l|l|}
\hline
\large{\bf{$\operatorname{Hom}_{B_r/B_l}(k,\h^2(B_1,k)^{(l-1)} \t p^l \la)$}} & Root System & \large{\bf{$\la \in X(T)$}}  \\
\hline
\hline
$p^r \g$ & All & $\la=p^{r-l}\g - \a, \a \in \Delta$ \\
\hline
$p^r \g$ & $A_3$ & $\la=p^{r-l}\g-\w_2$ \\
\hline
$p^r \g$ & $B_3$ & $\la=p^{r-l}\g-(\w_2-\w_1)$ \\
\hline
$p^r \g$ & $B_4$ & $\la=p^{r-l}\g-(\w_2+\w_3)$ \\
\hline
$p^r \g$ & $B_4$ & $\la=p^{r-l}\g-(\w_1-\w_2+\w_3-\w_4)$ \\
\hline
$p^r \g$ & $B_n$ & $\la=p^{r-l}\g-(\w_{n-1}-\w_{n-2})$ \\
\hline
$p^r \g$ & $C_n$ & $\la=p^{r-l}\g-(\w_{n}-\w_{n-2})$ \\
\hline
$p^r \g$ & $D_4$ & $\la=p^{r-l}\g-(\w_1-\w_2+\w_i), i \in \{3,4\}$ \\
\hline
$p^r \g$ & $D_4$ & $\la=p^{r-l}\g-(-\w_2+\w_3-\w_4)$ \\
\hline
$p^r \g$ & $D_n$ & $\la=p^{r-l}\g-(-\w_{n-2}+\w_{n-1}+\w_n)$ \\
\hline
$p^r \g$ & $F_4$ & $\la=p^{r-l}\g-(-\w_1+\w_2-\w_4)$ \\
\hline
$p^r \g$ & $G_2$ & $\la=p^{r-l}\g-(-\w_1+\w_2)$ \\
\hline
\end{longtable}

%% file: B_rtable.tex
\begin{longtable}{|l|l|l|}
\hline 
\large{\bf{$\h^2(B_r, \la), \; p=2$}} & Root System & \large{\bf{$\la \in X_r(T)$}} \\
\hline\hline 
$\h^2(B_1, w \cdot 0 + p \v)^{(r-1)}$ & $A_n (n\not=3), B_n, E_6$  & $p^{r-1}(w \cdot 0 + p \v), l(w)=0,2$\\
& $E_7, E_8, F_4, G_2$  & \\
\hline
$\h^2(B_1, w \cdot 0 + p \v)^{(r-1)}$ & $A_3$ & $p^r\v$ \\
\hline
$\h^2(B_1, w \cdot 0 + p \v)^{(r)}$ & $C_n, D_n$ & $p^{r-1}(w \cdot 0 + p \v), l(w)=0,2$ \\
\hline
$\v^{(r)}$ & All & $p^r\v + p^l w \cdot 0$ with $l(w)=2$ \\
& & $\quad$ and $0\leq l < r-1$ \\
\hline
$\v^{(r)}$ & $A_n, D_n, E_6$ & $p^r\v-p^l \a$ with $0 < l < r,\; \a \in \Delta$ \\
 & \quad $E_7, E_8$ & \\
\hline
$\v^{(r)}$ & $B_3, B_4, F_4, G_2$ & $p^r\v-p^l \a$ with $0 \leq l < r,\; \a \in \Delta$ \\
\hline
$\v^{(r)}$ & $B_n, n \geq 5$ & $p^r\v-p^l \a$ with $0 \leq l < r,\; \a \in \Delta;$ \\
& & $\quad l \not= r-1$ if $\a=\a_{n-1}$ \\
\hline
$\v^{(r)}$ & $C_n$ & $p^r\v-p^l \a$ with $0 \leq l < r,\; \a \in \Delta;$ \\
 & & $\quad l \not= r-1$ if $\a=\a_n$ \\
\hline
$\v^{(r)}$ & All & $p^r\v-p^t\b-p^l \a$ with $0 \leq l < t < r,$ \\
& & $\quad \a, \b \in \Delta$ \\
\hline
$\v^{(r)}$ & $A_3$ & $ p^r \v+p^{r-1} \a_2 - p^l \a$ with $0 \leq l < r-1,$ \\
& & $\quad \a \in \Delta$ \\
\hline
$(\v + \w_1)^{(r)} \oplus (\v + \w_3)^{(r)}$ & $A_3$ & $p^r\v+p^{r-1}\w_2-p^l\a$, with $0 \leq l < r-1$ \\ 
& & $\quad \a \in \Delta$ \\
\hline
$(\v + \w_1)^{(r)} \oplus (\v + \w_3)^{(r)}$ & $B_3$ & $p^r\v+p^{r-1}\a_2-p^l\a$, with $0 \leq l < r-1$ \\
& & $\quad \a \in \Delta$ \\
\hline
$\v^{(r)}$ & $B_n$ & $p^r\v -p^{l+1}(\a_{n-1}+\a_n)-p^l\a_2$ with \\
& & $\quad 0 \leq l < r-1$ \\
\hline
$M_{B_n}^{(r)} \t \v^{(r)}$ & $B_n, n \not= 4$ & $p^r\v-p^{r-1}\a_{n-1}-p^l\a$ with $0 \leq l < r-1,$ \\
& & $\quad \a \in \Delta$ \\
\hline
$(\w_1+M_{B_4})^{(r)} \t \v^{(r)}$ & $B_4$ & $p^r\v-p^{r-1}\a_i-p^l\a$ with $0 \leq l < r-1, \;$ \\
& & $\quad \a \in \Delta, \; i \in \{1,3\}$ \\
\hline
$M_{B_n}^{(r)} \t \v^{(r)}$ & $B_n$ & $p^r\v-p^{r-1}\a_{n-1} - p^l\a_n$ with $0 \leq l < r-1$ \\
\hline
$\v^{(r)}$ & $C_n$ & $p^r\v-p^l(\a_{n-1}+\a_n)$ with $0 \leq l < r-1$\\
\hline
$M^{(r)} \t \v^{(r)}$ & $C_n$ & $p^r\v - p^{r-1}\a_n - p^l\a$ with $0 \leq l < r-1,$ \\
& & $\quad  \a \in \Delta$ \\
\hline
$M^{(r)} \t \v^{(r)}$ & $C_n$ & $p^r\v-p^{r-1}\a$ where $\a \in \{\a_{n-1}, \a_n\}$ \\
\hline
$\h^1(B_{r-1},M^{(-1)} \t \la_1)$ & $C_n$ & $p\la_1$ \\
 $\quad  \oplus \h^2(B_{r-1}, \la_1)$ & &\\ 
\hline
$\v^{(r)}$ & $D_4$ & $p^r\v+p^{r-1}\a_2-p^l\a$ with $0 \leq l < r-1, $\\
& &$\quad \a \in \Delta$ \\
\hline
$(\v+\w_1)^{(r)} \oplus (\v + \w_3)^{(r)}$ & $D_4$ & $p^r\v+p^{r-1}\w_2-p^l\a$ with $0 \leq l < r-1,$ \\
& & $\quad  \a \in \Delta$ \\
$\quad  \oplus (\v+\w_4)^{(r)}$ & & \\
\hline
$\v^{(r)}$ & $D_n,  (n \geq 5)$ & $p^r\v+p^{r-1}\a_i-p^l\a$ with $0 \leq l < r-1,$ \\
\nopagebreak
& & $\quad \a \in \Delta, \; i \not=n-1,n$ \\
\hline
$(\v+\w_{n-1})^{(r)} \oplus (\v+\w_{n})^{(r)}$ & $D_n, (n \geq 5)$ & $p^r\v+p^{r-1}\w_{n-2}-p^l\a$ with $0 \leq l < r-1, $ \\
& & $\quad  \a \in \Delta$ \\
\hline
$\v^{(r)}$ & $F_4$ & $p^r\v-p^{l+1}(\a_3+\b)-p^l\a_2$ with \\
& & $\quad  0 \leq l < r-1, \; \b \in \{\a_2, \a_4\}$ \\
\hline
$M_{F_4}^{(r)} \t \v^{(r)}$ & $F_4$ & $p^r\v-p^{r-1}\a_2-p^l\a$ with $0 \leq l < r-1,$ \\
& & $\quad \a \in \Delta$ \\
\hline
$M^{(r)} \t \v^{(r)}$ & $F_4$ & $p^r\v - p^{r-1}\a_4 - p^l\a_2$ with \\
& & $\quad 0 \leq l < r-1$ \\
\hline
$\v^{(r)}$ & $G_2$ & $p^r\v-p^{l+1}(\a_1+\a_2) - p^l\a_2$ with \\
& & $\quad 0 \leq l < r-1,\; \a \in \Delta$ \\
\hline
$M_{G_2}^{(r)} \t \v^{(r)}$ & $G_2$ & $p^r\v-p^{r-1}\a_2 - p^l\a$ with $0 \leq l < r-1,$ \\
& & $\quad \a \in \Delta$ \\
\hline

\end{longtable}

%% file: appendix2.tex
\section{$G_r$-cohomology module structure}

If the $B_r$-cohomology involves an indecomposable $B$-module, then it is necessary to determine the $G_r$-cohomology structure separately.  Consider the following modules, with the factors listed from top to bottom.
\begin{itemize}
\item[$\bullet$] $M$ has factors $\a_i$ and $k$ corresponding to $w=s_\a s_\b$ where $\a,\b$ are separated by a single vertex, unless otherwise noted in one of the following modules.
\item[$\bullet$] If $\Phi=B_n$ and $w=s_{\a_{n-2}}s_{\a_n}$, then $M$ has factors $\a_{n-1}$ and $k$. 
\item[$\bullet$] If $\Phi=B_n$ and $w=s_{\a_{n-3}}s_{\a_{n-1}}$ then $M$ has factors $\a_{n-2}$ and $k$. 
\item[$\bullet$] If $\Phi=B_n$ and $w=s_{\a_i}s_{\a_{n-1}}$ with $i\not= n-3$ then $M$ has factors $\a_n$ and $k$.
\item[$\bullet$] If $\Phi=C_n$ and $w=s_{\a_i}s_{\a_n}$ with $i \not= n-2$ then $M$ has factors $\a_{n-1}$ and $k$.
\item[$\bullet$] If $\Phi=D_n$ and $w=s_{\a_{n-3}}s_{\a_{n-1}}$ then $M$ has factors $\a_{n-2}+\a_n, \a_{n-2}$ and $k$.
\item[$\bullet$] If $\Phi=D_n$ and $w=s_{\a_{n-3}}s_{\a_n}$ then $M$ has factors $\a_{n-2}+\a_{n-1}, \a_{n-2}$ and $k$.
\item[$\bullet$] If $\Phi=E_6$ and $w=s_{\a_2}s_{\a_3}$ then $M$ has factors $\a_4+\a_5+\a_6, \a_4+\a_5, \a_4$ and $k$.
\item[$\bullet$] If $\Phi$ is of type $E_6,E_7,E_8$ and $w=s_{\a_2}s_{\a_5}$ then $M$ has factors $\a_1+\a_3+\a_4,\a_3+\a_4, \a_4$ and $k$.
\item[$\bullet$] If $\Phi$ is of type $E_6,E_7,E_8$ and $w=s_{\a_3}s_{\a_5}$ then $M$ has factors $\a_2+\a_4, \a_4$ and $k$.
\item[$\bullet$] If $\Phi=E_7$ and $w=s_{\a_2}s_{\a_3}$ then $M$ has factors $\a_4+\a_5+\a_6+\a_7, \a_4+\a_5+\a_6,\a_4+\a_5, \a_4$ and $k$.
\item[$\bullet$] If $\Phi=E_8$ and $w=s_{\a_2}s_{\a_3}$ then $M$ has factors $\a_4+\a_5+\a_6+\a_7+\a_8, \a_4+\a_5+\a_6+\a_7,\a_4+\a_5+\a_6,\a_4+\a_5, \a_4$ and $k$.
\item[$\bullet$] If $\Phi=F_4$ and $w=s_{\a_1}s_{\a_3}$ then $M$ has factors $\a_2+\a_3, \a_2$ and $k$.
\end{itemize}
If $\la = w \cdot 0 + p\nu \in X(T)_+$ where $w \in W$ is one of the reflections listed above, then $$\h^2(G_r,H^0(\la)) \cong \operatorname{ind}_B^G(M \t \nu)^{(1)}.$$  The structure of these modules are listed below.  Since $\la \in X(T)_+$, we must look at the module structure carefully.  Furthermore, the module decomposition demonstrates that $\h^2(G_r,H^0(\la))$ has a good filtration.  The good filtration satisfies Donkin's Conjecture for $V = H^0 (\la)$, \cite{VdK} factors listed below are listed from top to bottom.\\

\begin{lem}
Let $p=2$ and $M$ be a module as above with corresponding $w \in W$.  Suppose $\v \in X(T)$ with $w \cdot 0 + p\nu \in X(T)_+$.
\begin{itemize}
\item[(a)] Suppose $w= s_\a s_\b$ with $\a + \b \not\in \Phi^+$ and $\a + \b + \gamma \in \Phi^+$.  Then $\left<\nu,\a_i^\vee\right> \geq 0$ for $\a_i \not= \a,\b,\g$, $\left<\nu,\a^\vee\right> \geq 1, \left<\nu,\g^\vee\right> \geq -1, \left<\nu,\b^\vee\right> \geq 1$.  Let $\d \in \Delta$ be such that $\d + \a \in \Phi^+$ and $\d \not= \g$.  Furthermore,\\
\begin{itemize}
\item[(i)] If $\left<\nu,\g^\vee\right>=-1$, then $\operatorname{ind}_B^G(M \t \nu)=0$.
\vspace{.25cm}
\item[(ii)] If $\left<\nu,\g^\vee\right>=0$, then $\operatorname{ind}_B^G(M \t \nu)\cong H^0(\nu)$.
\vspace{.25cm}
\item[(iii)] If $\left<\nu,\g^\vee\right> \geq 1, \left<\nu,\d^\vee\right> =0$, then $\operatorname{ind}_B^G(M \t \nu) \cong H^0(\nu)$.
\vspace{.25cm}
\item[(iv)] If $\left<\nu,\g^\vee\right> \geq 1, \left<\nu, \d^\vee\right> \geq 1$, then $\operatorname{ind}_B^G(M \t \nu)$ has factors $H^0(\a+\nu)$ and $H^0(\nu)$.\\
\\
\end{itemize}
\item[(b)] $\Phi$ is of type $B_n$ with $w=s_{\a_{n-2}}s_{\a_n}$.  Then $\left<\nu,\a_i^\vee\right> \geq 0$ for $1 \leq i \leq n-3$, $\displaystyle{\left<\nu,\a_{n-2}^\vee\right> \geq 1, \left<\nu,\a_{n-1}^\vee\right> \geq -1, \left<\nu,\a_n^\vee\right> \geq 1}$.  Furthermore,\\
\begin{itemize}
\item[(i)] If $\left<\v,\a_{n-1}^\vee\right>=-1$ and $\left<\v,\a_n^\vee\right>=1$, then $\operatorname{ind}_B^G(M \t \nu)=0$.
\vspace{.25cm}
\item[(ii)] If $\left<\v,\a_{n-1}^\vee\right>=-1$ and $\left<\v,\a_n^\vee\right>\geq 2$, then $\operatorname{ind}_B^G(M \t \nu)=H^0(\a_{n-1}+\v)$.
\vspace{.25cm}
\item[(iii)] If $\left<\v,\a_{n-1}^\vee\right>\geq 0$ and $\left<\v,\a_n^\vee\right>=1$, then $\operatorname{ind}_B^G(M \t \nu)=H^0(\v)$.
\vspace{.25cm}
\item[(iv)] If $\left<\v,\a_{n-1}^\vee\right>\geq 0$ and $\left<\v,\a_n^\vee\right> \geq 2$, then $\operatorname{ind}_B^G(M \t \nu)$ has factors $H^0(\a_{n-1}+\v)$ and $H^0(\v)$.\\
\\
\end{itemize}
\item[(c)] $\Phi$ is of type $B_n$ with $w=s_{\a_{n-3}}s_{\a_{n-1}}$.  Then $\left<\nu,\a_i^\vee\right> \geq 0$ for $1 \leq i \leq n-4$, $\displaystyle{\left<\nu,\a_{n-3}^\vee\right> \geq 1, \left<\nu,\a_{n-2}^\vee\right> \geq -1, \left<\nu,\a_{n-1}^\vee\right> \geq 1}$ and $\left<\nu,\a_n^\vee\right> \geq -1$.  Furthermore,\\
\begin{itemize}
\item[(i)] If $\left<\v,\a_{n}^\vee\right>=-1$, then $\operatorname{ind}_B^G(M \t \nu)=0$.
\vspace{.25cm}
\item[(ii)] If $\left<\v,\a_{n}^\vee\right> \geq 0$ and $\left<\v,\a_{n-2}^\vee\right>=-1$, then $\operatorname{ind}_B^G(M \t \nu)=H^0(\a_{n-2}+\v)$.
\vspace{.25cm}
\item[(iii)] If $\left<\v,\a_{n}^\vee\right>\geq 0$ and $\left<\v,\a_{n-2}^\vee\right> \geq 0$, then $\operatorname{ind}_B^G(M \t \nu)$ has factors $H^0(\a_{n-2}+\v)$ and $H^0(\v)$.\\
\\
\end{itemize}
\item[(d)] $\Phi$ is of type $B_n$ with $w=s_{\a_j}s_{\a_{n-1}}$ with $j\not= n-3$.  Then $\left<\nu,\a_i^\vee\right> \geq 0$ for $\displaystyle{1 \leq i \leq n-2, i\not=j}$, $\displaystyle{\left<\nu,\a_{j}^\vee\right> \geq 1, \left<\nu,\a_{n-1}^\vee\right> \geq 1}$ and $\left<\nu,\a_n^\vee\right> \geq -1$.  Furthermore,\\
\begin{itemize}
\item[(i)] If $\left<\v,\a_{n}^\vee\right>=-1$, then $\operatorname{ind}_B^G(M \t \nu)=H^0(\a_n+\v)$.
\item[(ii)] If $\left<\v,\a_{n}^\vee\right> \geq 0$, then $\operatorname{ind}_B^G(M \t \nu)$ has factors $H^0(\a_n+\v)$ and $H^0(\v)$.
\\
\end{itemize}
\item[(e)] $\Phi$ is of type $C_n$ with $w=s_{\a_{j}}s_{\a_{n}}$ with $j \not=n-2$.  Then $\left<\nu,\a_i^\vee\right> \geq 0$ for $1 \leq i \leq n-2$ and $i\not=j$, $\left<\nu,\a_{j}^\vee\right> \geq 1, \left<\nu,\a_{n-1}^\vee\right> \geq -1$, and $\left<\nu,\a_{n}^\vee\right> \geq 1$.  Furthermore,\\
\begin{itemize}
\item[(i)] If $\left<\v,\a_{n-2}^\vee\right>=0$ and $\left<\v,\a_{n-1}^\vee\right>=-1$, then $\operatorname{ind}_B^G(M \t \nu)=0$.
\vspace{.25cm}
\item[(ii)] If $\left<\v,\a_{n-2}^\vee\right> \geq 1$ and $\left<\v,\a_{n-1}^\vee\right>=-1$, then $\operatorname{ind}_B^G(M \t \nu)=H^0(\a_{n-1}+\v)$.
\vspace{.25cm}
\item[(iii)] If $\left<\v,\a_{n-2}^\vee\right>= 0$ and $\left<\v,\a_{n-1}^\vee\right> \geq 0$, then $\operatorname{ind}_B^G(M \t \nu) \cong H^0(\v)$.
\vspace{.25cm}
\item[(iv)] If $\left<\v,\a_{n-2}^\vee\right> \geq 1$ and $\left<\v,\a_{n-1}^\vee\right> \geq 0$ then $\operatorname{ind}_B^G(N_{C_n} \t \v)$ has factors $H^0(\a_{n-1}+\v)$ and $H^0(\v)$.\\
\\
\end{itemize}
\item[(f)] $\Phi$ is of type $C_n$ with $w=s_{\a_{n-2}}s_{\a_{n}}$.  Then $\left<\nu,\a_i^\vee\right> \geq 0$ for $1 \leq i \leq n-3$, $\displaystyle{\left<\nu,\a_{n-2}^\vee\right> \geq 1, \left<\nu,\a_{n-1}^\vee\right> \geq -1}$, and $\left<\nu,\a_n^\vee\right> \geq 1$.  Furthermore,\\
\begin{itemize}
\item[(i)] If $\left<\v,\a_{n-1}^\vee\right>=-1$, then $\operatorname{ind}_B^G(M \t \nu)=0$.
\vspace{.25cm}
\item[(ii)] If $\left<\v,\a_{n-3}^\vee\right> = 0$ or $\left<\v,\a_{n-1}^\vee\right>=0$, then $\operatorname{ind}_B^G(M \t \nu)=H^0(\v)$.
\vspace{.25cm}
\item[(iii)] If $\left<\v,\a_{n-3}^\vee\right>\geq 1$ and $\left<\v,\a_{n-1}^\vee\right> \geq 1$, then $\operatorname{ind}_B^G(M \t \nu)$ has factors $H^0(\a_{n-2}+\v)$ and $H^0(\v)$.\\
\\
\end{itemize}
\item[(g)] $\Phi$ is of type $D_n$ with $w=s_{\a_{n-3}}s_{\a_{n-1}}$.  Then $\left<\v,\a_i^\vee\right> \geq 0$ for $1 \leq i \leq n-4$, $\displaystyle{\left<\v,\a_{n-3}^\vee\right> \geq 1, \left<\v,\a_{n-2}^\vee\right> \geq -1}$, $\left<\v,\a_{n-1}^\vee\right> \geq 1$, and $\left<\v,\a_n^\vee\right> \geq 0$.  Furthermore,\\
\begin{itemize}
\item[(i)]If $\left<\v,\a_{n-2}^\vee\right> =-1$ and $\left<\v,\a_n^\vee\right>=0$, then $\operatorname{ind}_B^G(M \t \nu)=H^0(\a_{n-2}+\a_n+\v)$.
\vspace{.25cm}
\item[(ii)]If $\left<\v,\a_{n-2}^\vee\right> =-1$ and $\left<\v,\a_n^\vee\right> \geq 1$, then $\operatorname{ind}_B^G(M \t \nu)$ has factors $H^0(\a_{n-2}+\a_n+\v)$ and $H^0(\a_{n-2}+\v)$.
\vspace{.25cm}
\item[(iii)] If $\left<\v,\a_{n-2}^\vee\right> \geq 0$ and $\left<\v,\a_n^\vee\right>=0$, then $\operatorname{ind}_B^G(M \t \nu)$ has factors $H^0(\a_{n-2}+\a_n+\v)$ and $H^0(\v)$.
\vspace{.25cm}
\item[(iv)] If $\left<\v,\a_{n-2}^\vee\right> \geq 0$ and $\left<\v,\a_n^\vee\right> \geq 1$, then $\operatorname{ind}_B^G(M \t \nu)$ has factors\\ $H^0(\a_{n-2}+\a_n+\v),H^0(\a_{n-2}+\v)$ and $H^0(\v)$.\\
\\
\end{itemize}
\item[(h)] $\Phi$ is of type $D_n$ with $w=s_{\a_{n-3}}s_{\a_{n}}$.  Then $\left<\v,\a_i^\vee\right> \geq 0$ for $1 \leq i \leq n-4$, $\displaystyle{\left<\v,\a_{n-3}^\vee\right> \geq 1, \left<\v,\a_{n-2}^\vee\right> \geq -1}$, $\left<\v,\a_{n-1}^\vee\right> \geq 0$, and $\left<\v,\a_n^\vee\right> \geq 1$.  Furthermore,\\
\begin{itemize}
\item[(i)]If $\left<\v,\a_{n-1}^\vee\right> =-1$ and $\left<\v,\a_n^\vee\right>=0$, then $\operatorname{ind}_B^G(M \t \nu)=H^0(\a_{n-2}+\a_{n-1}+\v)$.
\vspace{.25cm}
\item[(ii)]If $\left<\v,\a_{n-1}^\vee\right> =-1$ and $\left<\v,\a_n^\vee\right> \geq 1$, then $\operatorname{ind}_B^G(M \t \nu)$ has factors\\ $H^0(\a_{n-2}+\a_{n-1}+\v)$ and $H^0(\a_{n-2}+\v)$.
\vspace{.25cm}
\item[(iii)] If $\left<\v,\a_{n-1}^\vee\right> \geq 0$ and $\left<\v,\a_n^\vee\right>=0$, then $\operatorname{ind}_B^G(M \t \nu)$ has factors $H^0(\a_{n-2}+\a_{n-1}+\v)$ and $H^0(\v)$.
\vspace{.25cm}
\item[(iv)] If $\left<\v,\a_{n-1}^\vee\right> \geq 0$ and $\left<\v,\a_n^\vee\right> \geq 1$, then $\operatorname{ind}_B^G(M \t \nu)$ has factors\\ $H^0(\a_{n-2}+\a_{n-1}+\v),H^0(\a_{n-2}+\v)$ and $H^0(\v)$.\\
\\
\end{itemize}
\item[(i)] $\Phi$ is of type $E_6$ with $w=s_{\a_2}s_{\a_3}$.  Then $\left<\v,\a_1^\vee\right> \geq 0, \left<\v,\a_2^\vee\right> \geq 1, \left<\v,\a_3^\vee\right> \geq 1$, $\displaystyle{\left<\v,\a_4^\vee\right> \geq -1, \left<\v,\a_5^\vee\right> \geq 0}$, and $\left<\v,\a_6^\vee\right> \geq 0$.  Furthermore,\\
\begin{itemize}
\item[(i)] If $\left<\v,\a_4^\vee\right> = -1, \left<\v,\a_5^\vee\right> = 0$ and $\left<\v,\a_6^\vee\right> =0$.  Then\\ $\operatorname{ind}_B^G(M \t \nu)=H^0(\a_{4}+\a_{5}+\a_6+\v)$.
\vspace{.25cm}
\item[(ii)]If $\left<\v,\a_4^\vee\right> = -1, \left<\v,\a_5^\vee\right> = 0$ and $\left<\v,\a_6^\vee\right> \geq 1$.  Then $\operatorname{ind}_B^G(M \t \nu)$ has factors $H^0(\a_4+\a_5+\a_6+\v)$ and $H^0(\a_4+\a_5+\v)$.
\vspace{.25cm}
\item[(iii)] If $\left<\v,\a_4^\vee\right> = -1, \left<\v,\a_5^\vee\right> \geq 1$ and $\left<\v,\a_6^\vee\right> =0$.  Then $\operatorname{ind}_B^G(M \t \nu)$ has factors $H^0(\a_4+\a_5+\a_6+\v)$ and $H^0(\a_4+\v)$.
\vspace{.25cm}
\item[(iv)] If $\left<\v,\a_4^\vee\right> = -1, \left<\v,\a_5^\vee\right> \geq 1$ and $\left<\v,\a_6^\vee\right> \geq 1$.  Then $\operatorname{ind}_B^G(M \t \nu)$ has factors $H^0(\a_{4}+\a_{5}+\a_6+\v)$, $H^0(\a_4+\a_5+\v)$ and $H^0(\a_4+\v)$.
\vspace{.25cm}
\item[(v)] If $\left<\v,\a_4^\vee\right> \geq 0, \left<\v,\a_5^\vee\right> = 0$ and $\left<\v,\a_6^\vee\right> =0$.  Then $\operatorname{ind}_B^G(M \t \nu)$ has factors\\ $H^0(\a_{4}+\a_{5}+\a_6+\v)$ and $H^0(\v)$.
\vspace{.25cm}
\item[(vi)] If $\left<\v,\a_4^\vee\right> \geq 0, \left<\v,\a_5^\vee\right> = 0$ and $\left<\v,\a_6^\vee\right> \geq 1$.  Then $\operatorname{ind}_B^G(M \t \nu)$ has factors\\ $H^0(\a_{4}+\a_{5}+\a_6+\v), H^0(\a_4+\a_5+\v)$ and $H^0(\v)$.
\vspace{.25cm}
\item[(vii)] If $\left<\v,\a_4^\vee\right> \geq 0, \left<\v,\a_5^\vee\right> \geq 1$ and $\left<\v,\a_6^\vee\right> =0$.  Then $\operatorname{ind}_B^G(M \t \nu)$ has factors\\ $H^0(\a_{4}+\a_{5}+\a_6+\v), H^0(\a_4+\v)$ and $H^0(\v)$.
\vspace{.25cm}
\item[(viii)] If $\left<\v,\a_4^\vee\right> \geq 0, \left<\v,\a_5^\vee\right> \geq 1$ and $\left<\v,\a_6^\vee\right> \geq 1$.  Then $\operatorname{ind}_B^G(M \t \nu)$ has factors $H^0(\a_{4}+\a_{5}+\a_6+\v)$, $H^0(\a_4+\a_5+\v), H^0(\a_4+\v)$ and $H^0(\v)$.\\
\\
\end{itemize}
\item[(j)] $\Phi$ is of type $E_6,E_7,E_8$ with $w=s_{\a_2}s_{\a_5}$.  Then $\left<\v,\a_i^\vee\right> \geq 0$, if $i \in \{1,3,6,7,8\}$, $\displaystyle{\left<\v,\a_2^\vee\right> \geq 1,  \left<\v,\a_4^\vee\right> \geq -1}$, and $\left<\v,\a_5^\vee\right> \geq 1$.  Furthermore,\\
\begin{itemize}
\item[(i)] If $\left<\v,\a_1^\vee\right> = 0, \left<\v,\a_3^\vee\right> = 0$ and $\left<\v,\a_4^\vee\right> =-1$.  Then\\ $\operatorname{ind}_B^G(M \t \nu)=H^0(\a_{1}+\a_{3}+\a_4+\v)$.
\vspace{.25cm}
\item[(ii)]If $\left<\v,\a_1^\vee\right> = 0, \left<\v,\a_3^\vee\right> \geq 1$ and $\left<\v,\a_4^\vee\right> = -1$.  Then $\operatorname{ind}_B^G(M \t \nu)$ has factors $H^0(\a_1+\a_3+\a_4+\v)$ and $H^0(\a_4+\v)$.
\vspace{.25cm}
\item[(iii)] If $\left<\v,\a_1^\vee\right> \geq 1, \left<\v,\a_3^\vee\right> =0$ and $\left<\v,\a_4^\vee\right> =-1$.  Then $\operatorname{ind}_B^G(M \t \nu)$ has factors $H^0(\a_1+\a_3+\a_4+\v)$ and $H^0(\a_3+\a_4+\v)$.
\vspace{.25cm}
\item[(iv)] If $\left<\v,\a_1^\vee\right> \geq 1, \left<\v,\a_3^\vee\right> \geq 1$ and $\left<\v,\a_4^\vee\right> =-1$.  Then $\operatorname{ind}_B^G(M \t \nu)$has factors $H^0(\a_{1}+\a_{3}+\a_4+\v)$, $H^0(\a_3+\a_4+\v)$ and $H^0(\a_4+\v)$.
\vspace{.25cm}
\item[(v)] If $\left<\v,\a_1^\vee\right> = 0, \left<\v,\a_3^\vee\right> = 0$ and $\left<\v,\a_4^\vee\right> \geq 0$.  Then $\operatorname{ind}_B^G(M \t \nu)$ has factors\\ $H^0(\a_{1}+\a_{3}+\a_4+\v)$ and $H^0(\v)$.
\vspace{.25cm}
\item[(vi)] If $\left<\v,\a_1^\vee\right> = 0, \left<\v,\a_3^\vee\right> \geq 1$ and $\left<\v,\a_4^\vee\right> \geq 0$.  Then $\operatorname{ind}_B^G(M \t \nu)$ has factors\\ $H^0(\a_{1}+\a_{3}+\a_4+\v), H^0(\a_4+\v)$ and $H^0(\v)$.
\vspace{.25cm}
\item[(vii)] If $\left<\v,\a_1^\vee\right> \geq 1, \left<\v,\a_3^\vee\right> =0$ and $\left<\v,\a_6^\vee\right> \geq 0$.  Then $\operatorname{ind}_B^G(M \t \nu)$ has factors\\ $H^0(\a_{1}+\a_{3}+\a_4+\v), H^0(\a_3+\a_4+\v)$ and $H^0(\v)$.
\vspace{.25cm}
\item[(viii)] If $\left<\v,\a_1^\vee\right> \geq 1, \left<\v,\a_3^\vee\right> \geq 1$ and $\left<\v,\a_4^\vee\right> \geq 0$.  Then $\operatorname{ind}_B^G(M \t \nu)$ has factors $H^0(\a_{1}+\a_{3}+\a_4+\v)$, $H^0(\a_3+\a_4+\v), H^0(\a_4+\v)$ and $H^0(\v)$.\\
\\
\end{itemize}
\item[(k)] $\Phi$ is of type $E_6,E_7,E_8$ with $w=s_{\a_3}s_{\a_5}$.  Then $\left<\v,\a_i^\vee\right> \geq 0$, if $i \in \{1,2,6,7,8\}$, $\displaystyle{\left<\v,\a_3^\vee\right> \geq 1,  \left<\v,\a_4^\vee\right> \geq -1}$, and $\left<\v,\a_5^\vee\right> \geq 1$.  Furthermore,\\
\begin{itemize}
\item[(i)] If $\left<\v,\a_2^\vee\right> = 0$ and $\left<\v,\a_4^\vee\right> =-1$.  Then $\operatorname{ind}_B^G(M \t \nu)=H^0(\a_{2}+\a_4+\v)$.
\vspace{.25cm}
\item[(ii)]If $\left<\v,\a_2^\vee\right> \geq 1$ and $\left<\v,\a_4^\vee\right> = -1$.  Then $\operatorname{ind}_B^G(M \t \nu)$ has factors $H^0(\a_2+\a_4+\v)$ and $H^0(\a_4+\v)$.
\vspace{.25cm}
\item[(iii)] If $\left<\v,\a_2^\vee\right> = 0$ and $\left<\v,\a_4^\vee\right> \geq 0$.  Then $\operatorname{ind}_B^G(M \t \nu)$ has factors $H^0(\a_2+\a_4+\v)$ and $H^0(\v)$.
\vspace{.25cm}
\item[(iv)] If $\left<\v,\a_2^\vee\right> \geq 1$ and $\left<\v,\a_4^\vee\right> \geq 0$.  Then $\operatorname{ind}_B^G(M \t \nu)$has factors $H^0(\a_2+\a_4+\v)$, $H^0(\a_4+\v)$, and $H^0(\v)$.\\
\\
\end{itemize}
\item[(l)] $\Phi$ is of type $E_6,E_7,E_8$ with $w=s_{\a_4}s_{\a_6}$.  Then $\left<\v,\a_i^\vee\right> \geq 0$ for $i \in \{1,2,3,7,8\}$, $\left<\v,\a_4^\vee\right> \geq 1$, $<\left<\v,\a_5^\vee\right> \geq -1$, and $\left<\v,a_6^\vee\right> \geq 1$.  Furthermore,\\
\begin{itemize}
\item[(i)] If $\left<\v,a_5^\vee\right> = -1$, then $\operatorname{ind}_B^G(M \t \v)=0$.
\vspace{.25cm}
\item[(ii)] If $\left<\v,a_5^\vee\right>=0$, then $\operatorname{ind}_B^G(M \t \v) \cong H^0(\v)$.
\vspace{.25cm}
\item[(iii)] If $\left<\v,a_5^\vee\right> \geq 1$ and ($\left<\v,a_2^\vee\right> = 0$ or $\left<\v,a_3^\vee\right> = 0$), then $\operatorname{ind}_B^G(M \t \v) \cong H^0(\v)$.
\vspace{.25cm}
\item[(iv)] If $\left<\v,\a_5^\vee\right> \geq 1, \left<\v,\a_2^\vee\right> \geq 1$ and $\left<\v,\a_3^\vee\right> \geq 1$ then $\operatorname{ind}_B^G(M \t \v)$ has factors $H^0(\a_4+\v)$ and $H^0(\v)$.\\
\\
\end{itemize}
\item[(m)] $\Phi$ is of type $E_7$ with $w=s_{\a_2}s_{\a_3}$.  Then $\left<\v,\a_i^\vee\right> \geq 0$, if $i \in \{1,5,6,7\}$, $\left<\v,\a_2^\vee\right> \geq 1$,\\  $\left<\v,\a_3^\vee\right> \geq 1$, and $\left<\v,\a_4^\vee\right> \geq -1$.  Furthermore,\\
\begin{itemize}
\item[(i)] If $\left<\v,\a_4^\vee\right> = -1, \left<\v,\a_5^\vee\right> = 0, \left<\v,\a_6^\vee\right> = 0$ and $\left<\v,\a_7^\vee\right> = 0$.  Then\\ $\operatorname{ind}_B^G(M \t \nu)=H^0(\a_{4}+\a_{5}+\a_6+\a_7+\v)$.
\vspace{.25cm}
\item[(ii)]If $\left<\v,\a_4^\vee\right> = -1, \left<\v,\a_5^\vee\right> = 0, \left<\v,\a_6^\vee\right> = 0$ and $\left<\v,\a_7^\vee\right> \geq 1$.  Then $\operatorname{ind}_B^G(M \t \nu)$ has factors $H^0(\a_4+\a_5+\a_6+\a_7+\v)$ and $H^0(\a_4+\a_5+\a_6+\v)$.
\vspace{.25cm}
\item[(iii)] If $\left<\v,\a_4^\vee\right> = -1, \left<\v,\a_5^\vee\right> =0, \left<\v,\a_6^\vee\right> \geq 1$ and $\left<\v,\a_7^\vee\right> =0$.  Then $\operatorname{ind}_B^G(M \t \nu)$ has factors $H^0(\a_4+\a_5+\a_6+\a_7+\v)$ and $H^0(\a_4+\a_5+\v)$.
\vspace{.25cm}
\item[(iv)] If $\left<\v,\a_4^\vee\right> = -1, \left<\v,\a_5^\vee\right> \geq 1, \left<\v,\a_6^\vee\right>=0$ and $\left<\v,\a_7^\vee\right> =0$.  Then $\operatorname{ind}_B^G(M \t \nu)$ has factors $H^0(\a_{4}+\a_{5}+\a_6+\a_7+\v)$ and $H^0(\a_4+\v)$.
\vspace{.25cm}
\item[(v)] If $\left<\v,\a_4^\vee\right> = -1, \left<\v,\a_5^\vee\right> = 0, \left<\v,\a_6^\vee\right> \geq 1$ and $\left<\v,\a_7^\vee\right> \geq 1$.  Then $\operatorname{ind}_B^G(M \t \nu)$ has factors $H^0(\a_{4}+\a_{5}+\a_6+\a_7+\v)$, $H^0(\a_4+\a_5+\a_6+\v),$ and $H^0(\a_4+\a_5)$.
\vspace{.25cm}
\item[(vi)] If $\left<\v,\a_4^\vee\right> = -1, \left<\v,\a_5^\vee\right> \geq 1, \left<\v,\a_6^\vee\right> = 0$ and $\left<\v,\a_7^\vee\right> \geq 1$.  Then $\operatorname{ind}_B^G(M \t \nu)$ has factors $H^0(\a_{4}+\a_{5}+\a_6+\a_7+\v), H^0(\a_4+\a_5+\a_6+\v)$ and $H^0(\a_4+\v)$.
\vspace{.25cm}
\item[(vii)] If $\left<\v,\a_4^\vee\right> =-1, \left<\v,\a_5^\vee\right> \geq 1, \left<\v,\a_6^\vee\right> \geq 1$ and $\left<\v,\a_7^\vee\right> = 0$.  Then $\operatorname{ind}_B^G(M \t \nu)$ has factors $H^0(\a_{4}+\a_{5}+\a_6+\a_7+\v), H^0(\a_4+\a_5+\v)$ and $H^0(\a_4+\v)$.
\vspace{.25cm}
\item[(viii)] If $\left<\v,\a_4^\vee\right> = -1, \left<\v,\a_5^\vee\right> \geq 1, \left<\v,\a_6^\vee\right> \geq 1$ and $\left<\v,\a_7^\vee\right> \geq 1$.  Then $\operatorname{ind}_B^G(M \t \nu)$ has factors $H^0(\a_{4}+\a_{5}+\a_6+\a_7+\v)$, $H^0(\a_4+\a_5+\a_6+\v), H^0(\a_4+\a_5+\v)$ and $H^0(\a_4+\v)$.
\vspace{.25cm}
\item[(ix)] If $\left<\v,\a_4^\vee\right> \geq 0, \left<\v,\a_5^\vee\right> = 0, \left<\v,\a_6^\vee\right> = 0$ and $\left<\v,\a_7^\vee\right> = 0$.  Then $\operatorname{ind}_B^G(M \t \nu)$ has factors $H^0(\a_{4}+\a_{5}+\a_6+\a_7+\v)$ and $H^0(\v)$.
\vspace{.25cm}
\item[(x)]If $\left<\v,\a_4^\vee\right> \geq 0, \left<\v,\a_5^\vee\right> = 0, \left<\v,\a_6^\vee\right> = 0$ and $\left<\v,\a_7^\vee\right> \geq 1$.  Then $\operatorname{ind}_B^G(M \t \nu)$ has factors $H^0(\a_4+\a_5+\a_6+\a_7+\v), H^0(\a_4+\a_5+\a_6+\v)$ and $H^0(\v)$.
\vspace{.25cm}
\item[(xi)] If $\left<\v,\a_4^\vee\right> \geq 0, \left<\v,\a_5^\vee\right> =0, \left<\v,\a_6^\vee\right> \geq 1$ and $\left<\v,\a_7^\vee\right> =0$.  Then $\operatorname{ind}_B^G(M \t \nu)$ has factors $H^0(\a_4+\a_5+\a_6+\a_7+\v), H^0(\a_4+\a_5+\v)$ and $H^0(\v)$.
\vspace{.25cm}
\item[(xii)] If $\left<\v,\a_4^\vee\right> \geq 0, \left<\v,\a_5^\vee\right> \geq 1, \left<\v,\a_6^\vee\right>=0$ and $\left<\v,\a_7^\vee\right> =0$.  Then $\operatorname{ind}_B^G(M \t \nu)$has factors $H^0(\a_{4}+\a_{5}+\a_6+\a_7+\v), H^0(\a_4+\v)$ and $H^0(\v)$.
\vspace{.25cm}
\item[(xiii)] If $\left<\v,\a_4^\vee\right> \geq 0, \left<\v,\a_5^\vee\right> = 0, \left<\v,\a_6^\vee\right> \geq 1$ and $\left<\v,\a_7^\vee\right> \geq 1$.  Then $\operatorname{ind}_B^G(M \t \nu)$ has factors $H^0(\a_{4}+\a_{5}+\a_6+\a_7+\v)$, $H^0(\a_4+\a_5+\a_6+\v), H^0(\a_4+\a_5)$ and $H^0(\v)$.
\vspace{.25cm}
\item[(xiv)] If $\left<\v,\a_4^\vee\right> \geq 0, \left<\v,\a_5^\vee\right> \geq 1, \left<\v,\a_6^\vee\right> = 0$ and $\left<\v,\a_7^\vee\right> \geq 1$.  Then $\operatorname{ind}_B^G(M \t \nu)$ has factors $H^0(\a_{4}+\a_{5}+\a_6+\a_7+\v), H^0(\a_4+\a_5+\a_6+\v), H^0(\a_4+\v)$ and $H^0(\v)$.
\vspace{.25cm}
\item[(xv)] If $\left<\v,\a_4^\vee\right> \geq 0, \left<\v,\a_5^\vee\right> \geq 1, \left<\v,\a_6^\vee\right> \geq 1$ and $\left<\v,\a_7^\vee\right> = 0$.  Then $\operatorname{ind}_B^G(M \t \nu)$ has factors $H^0(\a_{4}+\a_{5}+\a_6+\a_7+\v), H^0(\a_4+\a_5+\v), H^0(\a_4+\v)$ and $H^0(\v)$.
\vspace{.25cm}
\item[(xvi)] If $\left<\v,\a_4^\vee\right> = -1, \left<\v,\a_5^\vee\right> \geq 1, \left<\v,\a_6^\vee\right> \geq 1$ and $\left<\v,\a_7^\vee\right> \geq 1$.  Then $\operatorname{ind}_B^G(M \t \nu)$ has factors $H^0(\a_{4}+\a_{5}+\a_6+\a_7+\v)$, $H^0(\a_4+\a_5+\a_6+\v), H^0(\a_4+\a_5+\v), H^0(\a_4+\v),$ and $H^0(\v)$.\\
\\
\end{itemize}
\item[(n)] $\Phi$ is of type $E_8$ with $w=s_{\a_2}s_{\a_3}$.  Then $\left<\v,\a_i^\vee \right> \geq 0$, if $i \in \{1,5,6,7,8\}$, $\left<\v,\a_2^\vee \right> \geq 1$,\\  $\left<\v,\a_3^\vee \right> \geq 1$, and $\left<\v,\a_4^\vee \right> \geq -1$.  Furthermore,\\
\begin{itemize}
\item[(i)] If $\left<\v,\a_4^\vee\right> = -1, \left<\v,\a_5^\vee\right> = 0, \left<\v,\a_6^\vee\right> = 0, \left<\v,\a_7^\vee\right> = 0$ and $\left<\v,\a_8^\vee\right> = 0$.  Then $\operatorname{ind}_B^G(M \t \nu)=H^0(\a_{4}+\a_{5}+\a_6+\a_7+\a_8+\v)$.
\vspace{.25cm}
\item[(ii)]If $\left<\v,\a_4^\vee\right> = -1, \left<\v,\a_5^\vee\right> = 0, \left<\v,\a_6^\vee\right> = 0, \left<\v,\a_7^\vee\right> = 0$ and $\left<\v,\a_8^\vee\right> \geq 1$.  Then $\operatorname{ind}_B^G(M \t \nu)$ has factors $H^0(\a_4+\a_5+\a_6+\a_7+\a_8+\v)$ and $H^0(\a_4+\a_5+\a_6+\a_7+\v)$.
\vspace{.25cm}
\item[(iii)] If $\left<\v,\a_4^\vee\right> = -1, \left<\v,\a_5^\vee\right> =0, \left<\v,\a_5^\vee\right> =0, \left<\v,\a_7^\vee\right> \geq 1$ and $\left<\v,\a_8^\vee\right> =0$.  Then $\operatorname{ind}_B^G(M \t \nu)$ has factors $H^0(\a_4+\a_5+\a_6+\a_7+\a_8+\v)$ and $H^0(\a_4+\a_5+\a_6+\v)$.
\vspace{.25cm}
\item[(iv)] If $\left<\v,\a_4^\vee\right> = -1, \left<\v,\a_5^\vee\right> =0, \left<\v,\a_6^\vee\right> \geq 1, \left<\v,\a_7^\vee\right>=0$ and $\left<\v,\a_8^\vee\right> =0$.  Then $\operatorname{ind}_B^G(M \t \nu)$has factors $H^0(\a_{4}+\a_{5}+\a_6+\a_7+\a_8+\v)$ and $H^0(\a_4+\a_5+\v)$.
\vspace{.25cm}
\item[(v)] If $\left<\v,\a_4^\vee\right> = -1, \left<\v,\a_5^\vee\right> \geq 1,\left<\v,\a_6^\vee\right> = 0,\left<\v,\a_7^\vee\right> = 0,$ and $\left<\v,\a_8^\vee\right> =0$.  Then $\operatorname{ind}_B^G(M \t \nu)$ has factors $H^0(\a_{4}+\a_{5}+\a_6+\a_7+\a_8+\v)$ and $H^0(\a_4+\v)$.
\vspace{.25cm}
\item[(vi)] If $\left<\v,\a_4^\vee\right> = -1, \left<\v,\a_5^\vee\right> = 0,\left<\v,\a_6^\vee\right> = 0,\left<\v,\a_7^\vee\right> \geq 1,$ and $\left<\v,\a_8^\vee\right> \geq 1$.  Then $\operatorname{ind}_B^G(M \t \nu)$ has factors $H^0(\a_{4}+\a_{5}+\a_6+\a_7+\a_8+\v), H^0(\a_4+\a_5+\a_6+\a_7+\v)$ and $H^0(\a_4+\a_5+\a_6+\v)$.
\vspace{.25cm}
\item[(vii)] If $\left<\v,\a_4^\vee\right> =-1, \left<\v,\a_5^\vee\right> = 0,\left<\v,\a_6^\vee\right> \geq 1, \left<\v,\a_7^\vee\right> = 0$ and $\left<\v,\a_8^\vee\right> \geq 1$.  Then $\operatorname{ind}_B^G(M \t \nu)$ has factors $H^0(\a_{4}+\a_{5}+\a_6+\a_7+\a_8+\v), H^0(\a_4+\a_5+\a_6+\a_7+\v)$ and $H^0(\a_4+\a_5+\v)$.
\vspace{.25cm}
\item[(viii)] If $\left<\v,\a_4^\vee\right> = -1, \left<\v,\a_5^\vee\right> \geq 1, \left<\v,\a_6^\vee\right> = 0,\left<\v,\a_7^\vee\right> = 0,$ and $\left<\v,\a_8^\vee\right> \geq 1$.  Then $\operatorname{ind}_B^G(M \t \nu)$ has factors $H^0(\a_{4}+\a_{5}+\a_6+\a_7+\a_8+\v)$, $H^0(\a_4+\a_5+\a_6+\a_7+\v),$ and $H^0(\a_4+\v)$.
\vspace{.25cm}
\item[(ix)] If $\left<\v,\a_4^\vee\right> = -1, \left<\v,\a_5^\vee\right> = 0, \left<\v,\a_6^\vee\right> \geq 1,\left<\v,\a_7^\vee\right> \geq 1,$ and $\left<\v,\a_8^\vee\right> =0$.  Then $\operatorname{ind}_B^G(M \t \nu)$ has factors $H^0(\a_{4}+\a_{5}+\a_6+\a_7+\a_8+\v)$, $H^0(\a_4+\a_5+\a_6+\v),$ and $H^0(\a_4+\a_5+\v)$.
\vspace{.25cm}
\item[(x)] If $\left<\v,\a_4^\vee\right> = -1, \left<\v,\a_5^\vee\right> \geq 1, \left<\v,\a_6^\vee\right> = 0,\left<\v,\a_7^\vee\right> \geq 1,$ and $\left<\v,\a_8^\vee\right> = 0$.  Then $\operatorname{ind}_B^G(M \t \nu)$ has factors $H^0(\a_{4}+\a_{5}+\a_6+\a_7+\a_8+\v)$, $H^0(\a_4+\a_5+\a_6+\v),$ and $H^0(\a_4+\v)$.
\vspace{.25cm}
\item[(xi)] If $\left<\v,\a_4^\vee\right> = -1, \left<\v,\a_5^\vee\right> \geq 1, \left<\v,\a_6^\vee\right> \geq 1,\left<\v,\a_7^\vee\right> = 0,$ and $\left<\v,\a_8^\vee\right> =0$.  Then $\operatorname{ind}_B^G(M \t \nu)$ has factors $H^0(\a_{4}+\a_{5}+\a_6+\a_7+\a_8+\v)$, $H^0(\a_4+\a_5+\v),$ and $H^0(\a_4+\v)$.
\vspace{.2cm}
\item[(xii)] If $\left<\v,\a_4^\vee\right> = -1, \left<\v,\a_5^\vee\right> =0, \left<\v,\a_6^\vee\right> \geq 1,\left<\v,\a_7^\vee\right> \geq 1,$ and $\left<\v,\a_8^\vee\right> \geq 1$.  Then $\operatorname{ind}_B^G(M \t \nu)$ has factors $H^0(\a_{4}+\a_{5}+\a_6+\a_7+\a_8+\v)$, $H^0(\a_4+\a_5+\a_6+\a_7+\v)$, $H^0(\a_4+\a_5+\a_6+\v)$ and $H^0(\a_4+\a_5+\v)$.
\vspace{.25cm}
\item[(xiii)] If $\left<\v,\a_4^\vee\right> = -1, \left<\v,\a_5^\vee\right> \geq 1, \left<\v,\a_6^\vee\right> = 0,\left<\v,\a_7^\vee\right> \geq 1,$ and $\left<\v,\a_8^\vee\right> \geq 1$.  Then $\operatorname{ind}_B^G(M \t \nu)$ has factors $H^0(\a_{4}+\a_{5}+\a_6+\a_7+\a_8+\v)$, $H^0(\a_4+\a_5+\a_6+\a_7+\v)$, $H^0(\a_4+\a_5+\a_6+\v)$ and $H^0(\a_4+\v)$.
\vspace{.25cm}
\item[(xiv)] If $\left<\v,\a_4^\vee\right> = -1, \left<\v,\a_5^\vee\right> \geq 1, \left<\v,\a_6^\vee\right> = \geq 1,\left<\v,\a_7^\vee\right> = 0,$ and $\left<\v,\a_8^\vee\right> \geq 1$.  Then $\operatorname{ind}_B^G(M \t \nu)$ has factors $H^0(\a_{4}+\a_{5}+\a_6+\a_7+\a_8+\v)$, $H^0(\a_4+\a_5+\a_6+\v)$, $H^0(\a_4+\a_5+\v)$ and $H^0(\a_4+\v)$.
\vspace{.25cm}
\item[(xv)] If $\left<\v,\a_4^\vee\right> = -1, \left<\v,\a_5^\vee\right> \geq 1, \left<\v,\a_6^\vee\right> \geq 1,\left<\v,\a_7^\vee\right> \geq 1,$ and $\left<\v,\a_8^\vee\right> = 0$.  Then $\operatorname{ind}_B^G(M \t \nu)$ has factors $H^0(\a_{4}+\a_{5}+\a_6+\a_7+\a_8+\v)$, $H^0(\a_4+\a_5+\a_6+\v),H^0(\a_4+\a_5+\v)$ and $H^0(\a_4+\v)$.
\vspace{.25cm}
\item[(xvi)] If $\left<\v,\a_4^\vee\right> = -1, \left<\v,\a_5^\vee\right> \geq 1, \left<\v,\a_6^\vee\right> \geq 1,\left<\v,\a_7^\vee\right> \geq 1,$ and $\left<\v,\a_8^\vee\right> \geq 1$.  Then $\operatorname{ind}_B^G(M \t \nu)$ has factors $H^0(\a_{4}+\a_{5}+\a_6+\a_7+\a_8+\v)$,\\ $H^0(\a_4+\a_5+\a_6+\a_7+\v), H^0(\a_4+\a_5+\a_6+\v), H^0(\a_4+\a_5+\v),$ and $H^0(\a_4+\v)$.
\vspace{.25cm}
\item[(xvii)] If $\left<\v,\a_4^\vee\right> \geq 0, \left<\v,\a_5^\vee\right> = 0, \left<\v,\a_6^\vee\right> = 0, \left<\v,\a_7^\vee\right> = 0$ and $\left<\v,\a_8^\vee\right> = 0$.  Then $\operatorname{ind}_B^G(M \t \nu)$ has factors $H^0(\a_{4}+\a_{5}+\a_6+\a_7+\a_8+\v)$ and $H^0(\v)$.
\vspace{.25cm}
\item[(xviii)]If $\left<\v,\a_4^\vee\right> \geq 0, \left<\v,\a_5^\vee\right> = 0, \left<\v,\a_6^\vee\right> = 0, \left<\v,\a_7^\vee\right> = 0$ and $\left<\v,\a_8^\vee\right> \geq 1$.  Then $\operatorname{ind}_B^G(M \t \nu)$ has factors $H^0(\a_4+\a_5+\a_6+\a_7+\a_8+\v), H^0(\a_4+\a_5+\a_6+\a_7+\v)$ and $H^0(\v)$.
\vspace{.25cm}
\item[(xix)] If $\left<\v,\a_4^\vee\right> \geq 0, \left<\v,\a_5^\vee\right> =0, \left<\v,\a_5^\vee\right> =0, \left<\v,\a_7^\vee\right> \geq 1$, and  $\left<\v,\a_8^\vee\right> =0$.  Then $\operatorname{ind}_B^G(M \t \nu)$ has factors $H^0(\a_4+\a_5+\a_6+\a_7+\a_8+\v), H^0(\a_4+\a_5+\a_6+\v)$ and $H^0(\v)$.
\vspace{.25cm}
\item[(xx)] If $\left<\v,\a_4^\vee\right> \geq 0, \left<\v,\a_5^\vee\right> =0, \left<\v,\a_6^\vee\right> \geq 1, \left<\v,\a_7^\vee\right>=0$ and $\left<\v,\a_8^\vee\right> =0$.  Then $\operatorname{ind}_B^G(M \t \nu)$has factors $H^0(\a_{4}+\a_{5}+\a_6+\a_7+\a_8+\v), H^0(\a_4+\a_5+\v)$ and $H^0(\v)$.
\vspace{.25cm}
\item[(xxi)] If $\left<\v,\a_4^\vee\right> \geq 0, \left<\v,\a_5^\vee\right> \geq 1,\left<\v,\a_6^\vee\right> = 0,\left<\v,\a_7^\vee\right> = 0,$ and $\left<\v,\a_8^\vee\right> =0$.  Then $\operatorname{ind}_B^G(M \t \nu)$ has factors $H^0(\a_{4}+\a_{5}+\a_6+\a_7+\a_8+\v), H^0(\a_4+\v)$ and $H^0(\v)$.
\vspace{.25cm}
\item[(xxii)] If $\left<\v,\a_4^\vee\right> \geq 0, \left<\v,\a_5^\vee\right> = 0,\left<\v,\a_6^\vee\right> = 0,\left<\v,\a_7^\vee\right> \geq 1,$ and $\left<\v,\a_8^\vee\right> \geq 1$.  Then $\operatorname{ind}_B^G(M \t \nu)$ has factors $H^0(\a_{4}+\a_{5}+\a_6+\a_7+\a_8+\v)$,\\ $H^0(\a_4+\a_5+\a_6+\a_7+\v), H^0(\a_4+\a_5+\a_6+\v)$ and $H^0(\v)$.
\vspace{.25cm}
\item[(xxiii)] If $\left<\v,\a_4^\vee\right> \geq 0, \left<\v,\a_5^\vee\right> = 0,\left<\v,\a_6^\vee\right> \geq 1, \left<\v,\a_7^\vee\right> = 0$ and $\left<\v,\a_8^\vee\right> \geq 1$.  Then $\operatorname{ind}_B^G(M \t \nu)$ has factors $H^0(\a_{4}+\a_{5}+\a_6+\a_7+\a_8+\v)$,\\ $H^0(\a_4+\a_5+\a_6+\a_7+\v), H^0(\a_4+\a_5+\v)$ and $H^0(\v)$.
\vspace{.25cm}
\item[(xxiv)] If $\left<\v,\a_4^\vee\right> \geq 0, \left<\v,\a_5^\vee\right> \geq 1, \left<\v,\a_6^\vee\right> = 0,\left<\v,\a_7^\vee\right> = 0,$ and $\left<\v,\a_8^\vee\right> \geq 1$.  Then $\operatorname{ind}_B^G(M \t \nu)$ has factors $H^0(\a_{4}+\a_{5}+\a_6+\a_7+\a_8+\v)$,\\ $H^0(\a_4+\a_5+\a_6+\a_7+\v), H^0(\a_4+\v)$ and $H^0(\v)$.
\vspace{.25cm}
\item[(xxv)] If $\left<\v,\a_4^\vee\right> \geq 0, \left<\v,\a_5^\vee\right> = 0, \left<\v,\a_6^\vee\right> \geq 1,\left<\v,\a_7^\vee\right> \geq 1,$ and $\left<\v,\a_8^\vee\right> =0$.  Then $\operatorname{ind}_B^G(M \t \nu)$ has factors $H^0(\a_{4}+\a_{5}+\a_6+\a_7+\a_8+\v)$,\\ $H^0(\a_4+\a_5+\a_6+\v), H^0(\a_4+\a_5+\v)$ and $H^0(\v)$.
\vspace{.25cm}
\item[(xxvi)] If $\left<\v,\a_4^\vee\right> \geq 0, \left<\v,\a_5^\vee\right> \geq 1, \left<\v,\a_6^\vee\right> = 0,\left<\v,\a_7^\vee\right> \geq 1,$ and $\left<\v,\a_8^\vee\right> = 0$.  Then $\operatorname{ind}_B^G(M \t \nu)$ has factors $H^0(\a_{4}+\a_{5}+\a_6+\a_7+\a_8+\v)$,\\ $H^0(\a_4+\a_5+\a_6+\v), H^0(\a_4+\v)$ and $H^0(\v)$.
\vspace{.25cm}
\item[(xxvii)] If $\left<\v,\a_4^\vee\right> \geq 0, \left<\v,\a_5^\vee\right> \geq 1, \left<\v,\a_6^\vee\right> \geq 1,\left<\v,\a_7^\vee\right> = 0,$ and $\left<\v,\a_8^\vee\right> =0$.  Then $\operatorname{ind}_B^G(M \t \nu)$ has factors $H^0(\a_{4}+\a_{5}+\a_6+\a_7+\a_8+\v)$, $H^0(\a_4+\a_5+\v), H^0(\a_4+\v)$ and $H^0(\v)$.
\vspace{.25cm}
\item[(xxviii)] If $\left<\v,\a_4^\vee\right> \geq 0, \left<\v,\a_5^\vee\right> =0, \left<\v,\a_6^\vee\right> \geq 1,\left<\v,\a_7^\vee\right> \geq 1,$ and $\left<\v,\a_8^\vee\right> \geq 1$.  Then $\operatorname{ind}_B^G(M \t \nu)$ has factors $H^0(\a_{4}+\a_{5}+\a_6+\a_7+\a_8+\v)$,\\ $H^0(\a_4+\a_5+\a_6+\a_7+\v),H^0(\a_4+\a_5+\a_6+\v), H^0(\a_4+\a_5+\v)$ and $H^0(\v)$.
\vspace{.25cm}
\item[(xxix)] If $\left<\v,\a_4^\vee\right> \geq 0, \left<\v,\a_5^\vee\right> \geq 1, \left<\v,\a_6^\vee\right> = 0,\left<\v,\a_7^\vee\right> \geq 1,$ and $\left<\v,\a_8^\vee\right> \geq 1$.  Then $\operatorname{ind}_B^G(M \t \nu)$ has factors $H^0(\a_{4}+\a_{5}+\a_6+\a_7+\a_8+\v)$,\\ $H^0(\a_4+\a_5+\a_6+\a_7+\v), H^0(\a_4+\a_5+\a_6+\v), H^0(\a_4+\v)$ and $H^0(\v)$.
\vspace{.25cm}
\item[(xxx)] If $\left<\v,\a_4^\vee\right> \geq 0, \left<\v,\a_5^\vee\right> \geq 1, \left<\v,\a_6^\vee\right> = \geq 1,\left<\v,\a_7^\vee\right> = 0,$ and $\left<\v,\a_8^\vee\right> \geq 1$.  Then $\operatorname{ind}_B^G(M \t \nu)$ has factors $H^0(\a_{4}+\a_{5}+\a_6+\a_7+\a_8+\v)$,\\ $H^0(\a_4+\a_5+\a_6+\v), H^0(\a_4+\a_5+\v), H^0(\a_4+\v)$ and $H^0(\v)$.
\vspace{.25cm}
\item[(xxxi)] If $\left<\v,\a_4^\vee\right> \geq 0, \left<\v,\a_5^\vee\right> \geq 1, \left<\v,\a_6^\vee\right> \geq 1,\left<\v,\a_7^\vee\right> \geq 1,$ and $\left<\v,\a_8^\vee\right> = 0$.  Then $\operatorname{ind}_B^G(M \t \nu)$ has factors $H^0(\a_{4}+\a_{5}+\a_6+\a_7+\a_8+\v)$,\\ $H^0(\a_4+\a_5+\a_6+\v),H^0(\a_4+\a_5+\v), H^0(\a_4+\v)$ and $H^0(\v)$.
\vspace{.25cm}
\item[(xxxii)] If $\left<\v,\a_4^\vee\right> \geq 0, \left<\v,\a_5^\vee\right> \geq 1, \left<\v,\a_6^\vee\right> \geq 1,\left<\v,\a_7^\vee\right> \geq 1,$ and $\left<\v,\a_8^\vee\right> \geq 1$.  Then $\operatorname{ind}_B^G(M \t \nu)$ has factors $H^0(\a_{4}+\a_{5}+\a_6+\a_7+\a_8+\v)$,\\ $H^0(\a_4+\a_5+\a_6+\a_7+\v), H^0(\a_4+\a_5+\a_6+\v), H^0(\a_4+\a_5+\v), H^0(\a_4+\v)$ and $H^0(\v)$.\\
\\
\end{itemize}
\item[(o)] If $\Phi$ is of type $F_4$ and $w=s_{\a_2}s_{\a_4}$.  Then $\left<\v,\a_1^\vee\right> \geq 0, \left<\v,\a_2^\vee\right> \geq 1, \left<\v,\a_3^\vee\right> \geq -1,$ and $\left<\v,\a_4^\vee\right> \geq 1$.  Furthermore,\\
\begin{itemize}
\item[(i)] If $\left<\v,\a_3^\vee\right>=-1$ then $\operatorname{ind}_B^G(N_{F_4}\t\v) = 0$.
\vspace{.25cm}
\item[(ii)] If $\left<\v,\a_3^\vee\right>=0,1$ or $\left<\v,\a_1^vee\right>=0$ then $\operatorname{ind}_B^G(N_{F_4}\t\v) \cong H^0(\v).$
\vspace{.25cm}
\item[(iii)] If $\left<\v,\a_3^\vee\right> \geq 2$ and $\left<\v,\a_1^vee\right> \geq 1$ then $\operatorname{ind}_B^G(N_{F_4}\t\v)$ has factors $H^0(\a_2+\v)$ and $H^0(\v)$.\\
\\
\end{itemize}
\item[(p)] If $\Phi$ is of type $F_4$ and $w=s_{\a_1}s_{\a_3}$.  Then $\left<\v,\a_1^\vee\right> \geq 1, \left<\v,\a_2^\vee\right> \geq -1, \left<\v,\a_3^\vee\right> \geq 1,$ and $\left<\v,\a_4^\vee\right> \geq 0$.  Furthermore,\\
\begin{itemize}
\item[(i)] If $\left<\v,\a_2^\vee\right>=-1, \left<\v,\a_3^\vee\right>=1$ and $\left<\v,\a_4^\vee\right>=0$ then $\operatorname{ind}_B^G(N_{F_4}\t\v) = 0$.
\vspace{.25cm}
\item[(ii)] If $\left<\v,\a_2^\vee\right>=-1, \left<\v,\a_3^\vee\right>=1$ and $\left<\v,\a_4^\vee\right> \geq 1$ then\\ $\operatorname{ind}_B^G(N_{F_4}\t\v) = H^0(\a_2+\a_3+\v).$
\vspace{.25cm}
\item[(iii)] If $\left<\v,\a_2^\vee\right>=-1, \left<\v,\a_3^\vee\right> \geq 2$ and $\left<\v,\a_4^\vee\right>=0$ then $\operatorname{ind}_B^G(N_{F_4}\t\v) = H^0(\a_2+\v)$.
\vspace{.25cm}
\item[(iv)] If $\left<\v,\a_2^\vee\right>=-1, \left<\v,\a_3^\vee\right> \geq 2$ and $\left<\v,\a_4^\vee\right> \geq 1$ then $\operatorname{ind}_B^G(N_{F_4}\t\v)$ has factors $H^0(\a_2+\a_3+\v)$ and $H^0(\a_2+\v)$.
\vspace{.25cm}
\item[(v)] If $\left<\v,\a_5^\vee\right> \geq 0, \left<\v,\a_3^\vee\right> = 1$ and $\left<\v,\a_4^\vee\right> =0$ then $\operatorname{ind}_B^G(N_{F_4}\t\v)\cong H^0(\v)$.
\vspace{.25cm}
\item[(vi)] If $\left<\v,\a_2^\vee\right> \geq 0, \left<\v,\a_3^\vee\right> =1$ and $\left<\v,\a_4^\vee\right> \geq 1$ then $\operatorname{ind}_B^G(N_{F_4}\t\v)$ has factors\\ $H^0(\a_2+\a_3+\v)$ and $H^0(\v)$.
\vspace{.25cm}
\item[(vii)] If $\left<\v,\a_2^\vee\right> \geq 0, \left<\v,\a_3^\vee\right> \geq 2$ and $\left<\v,\a_4^\vee\right> =0$ then $\operatorname{ind}_B^G(N_{F_4}\t\v)$ has factors\\ $H^0(\a_2+\v)$ and $H^0(\v)$. 
\vspace{.25cm}
\item[(viii)] If $\left<\v,\a_2^\vee\right> \geq 0, \left<\v,\a_3^\vee\right> \geq 2$ and $\left<\v,\a_4^\vee\right> \geq 1$ then $\operatorname{ind}_B^G(N_{F_4}\t\v)$ has factors\\ $H^0(\a_2+\a_3+\v), H^0(\a_2+\v)$ and $H^0(\v)$.\\  
\end{itemize}
\end{itemize}
\end{lem}
\begin{proof}
The proof is similar for all cases.\\
Consider the case $w = s_{\a_i}s_{\a_{i+2}}$.  Then $w \cdot 0 = \omega_{i-1}-2\omega_i+2\omega_{i+1}-2\omega_{i+2}+\omega_{i+3}$.  Let $\v = \sum_{i=1}^n c_i\omega_i$.  If $w \cdot 0 + 2\v$ is dominant, then $c_j \geq 0$ for $j \not\in \{i,i+1,i+2\}$, $c_{i} \geq 1, c_{i+1} \geq -1$ and $c_{i+2} \geq 1$.  To determine the structure of the induced modules, the argument follows \cite[Proposition 3.4]{BNP1}, which explains that it is necessary to determine precisely when our module factors are dominant.  The module $M \t \v$ has factors $\a_i + \v$ and $\v$.  Consider 
$$\a_i + \v = \sum_{j=1}^{i-2} c_j \omega_j + (c_{i-1}-1)\omega_{i-1} + (c_i+2)\omega_i + (c_{i+1}-1)\omega_{i+1} + \sum_{j={i+2}}^n c_j\omega_j$$
which is dominant precisely when $c_{i-1} \geq 1$ and $c_{i+1} \geq 1$.\\
\indent  If $\left<\v,\a_{i+1}^\vee\right> =-1$, then $\v$ and $\a_i+\v$ aren't dominant and $\operatorname{ind}_B^G(M \t \v)=0$.\\
\indent  If $\left<\v,\a_{i+1}^\vee\right> =0$ or $\left<\v,a_{i-1}\right>=0$, then $\a_i+\v$ isn't dominant, but $\v$ is, so using \cite{BNP1} and \cite[II.4.5]{Jan1}, then $\operatorname{ind}_B^G(\v) = H^0(\v)$; hence $\operatorname{ind}_B^G(M \t \v)=H^0(\v)$.  If $\left<\v,\a_{i+1}^\vee\right> \geq 1$ and $\left<\v,a_{i-1}\right> and \geq 1$, then both $\a_i+\v$ and $\v$ are dominant.  Using \cite[II.4.5]{Jan1} and \cite[3.4]{BNP1} then $\operatorname{ind}_B^G(\a_i+\v) = H^0(\a_i+\v)$ and $\operatorname{ind}_B^G(\v) = H^0(\v)$.  Hence $\operatorname{ind}_B^G(M\t\v)$ has factors $H^0(\a_i+\v)$ and $H^0(\v)$.
\end{proof}